\documentclass[3p,preprint]{elsarticle}
\usepackage{amsmath,amssymb,amsthm}
\usepackage{xcolor}
\usepackage{graphicx}
\graphicspath{ {figure} } % The figure path
\usepackage{array}
\usepackage{multirow}
\usepackage{arydshln}
\usepackage{booktabs}
\usepackage{subcaption}

\usepackage{setspace}
\usepackage[left]{lineno}
\usepackage{hyperref,cleveref}
\usepackage[normalem]{ulem}
\usepackage{float}

\theoremstyle{remark}
\newtheorem{remark}{Remark}

\theoremstyle{definition}
\newtheorem{definition}{Definition}

\theoremstyle{plain}
\newtheorem{prop}{Proposition}

\title{Time-Invariant Neural Operators with Applications in Solving Time-Dependent PDEs}
\author[1,2]{Zihan Zhou \fnref{fn1}}
\author[2,3]{Wenzhong Zhang \fnref{fn1}}
\author[4]{Lizuo Liu \corref{cor1}}

\fntext[fn1]{Equal contribution.}

\cortext[cor1]{Corresponding author. \textit{Email addresses}: \texttt{Lizuo.Liu@Dartmouth.edu} (LL), \texttt{zihanzhou@mail.ustc.edu.cn} (ZZ), \texttt{wenzhong@ustc.edu.cn} (WZ).}

\affiliation[1]{organization={University of Science and Technology of China},
postcode={230027},
city={Hefei},
country={P.R. China}}

\affiliation[2]{organization={Suzhou Institute for Advanced Research, University of Science and Technology of China},
postcode={215123},
city={Suzhou},
country={P.R. China}}

\affiliation[3]{organization={Key Laboratory of the Ministry of Education for Mathematical Foundations and Applications of Digital Technology},
postcode={230027},
city={Hefei},
country={P.R. China}}

\affiliation[4]{organization={Dartmouth College},
postcode={03755},
city={Hanover, NH},
country={USA}}

\begin{document}

\begin{abstract}
The deep operator network (DeepONet) is one of the basic architectures for learning nonlinear operators with neural networks.
However, for operators that describe the dynamic response of physical systems, DeepONet does not naturally respect fundamental time properties, including time causality and time invariance.
We propose time-invariant neural operator (TINO) to overcome these limitations, by creating the connection between time-delay neural network (TDNN) and DeepONet.
We then extend the proposal to solving time-dependent PDEs with initial values, which can be treated as truncated time-invariant systems, with spatial proper orthogonal decomposition (POD) implementation in the output function space.
Various numerical tests comparing TINO, neural operators with time causality, and those without time properties justify the enhanced precision of the proposed neural operator frameworks.
\end{abstract}

\maketitle

\section{Introduction}\label{sect_intro}

Learning the dynamic response of physical systems is a fundamental task in engineering and in scientific computing, with broad applications in control theory, signal processing, seismology, etc.
In many applications, the object of interest is the operator that maps a time-varying input signal, forcing term, or source field to a response function over time and space.
Two foundational temporal properties, \emph{time causality} and \emph{time invariance}, play fundamental roles to understanding and modeling these systems.
Time causality describes that the output at any time $t$ depends only on present and past input, whereas time invariance requires that a temporal shift of the input produces the corresponding shift of the output.
In signal processing, time invariance is equivalently described as \emph{shift invariance}, which indicates the commutative property that ensures consistent behavior regardless of arbitrary shifting operations.

These temporal properties are well understood in classic frameworks.
Linear time-invariant (LTI) systems admit convolution representations, while nonlinear operators are related to Volterra series, Wiener-type expansions, and state-space formulations, to name a few.
In terms of PDE solvers for finding single trajectories, numerical solvers are typically built in a time marching manner to incrementally solve the output by ascending time steps, hence strictly obeying time causality.
However, these temporal properties are not automatically inherited by data-driven operator learning methods.
Since the models are learned directly from data, it raises a nontrivial modeling problem to preserve temporal properties for prediction tasks.
Namely, time causality rejects future input to the model, and time invariance prohibits the model from explicitly depending on time.

In recent years, advancements in deep learning have made operator learning powerful and versatile for modeling dynamical systems.
Approaches such as recurrent neural networks (RNN) and their variants, Mori--Zwanzig Net \cite{fu2020learning}, OnsargarNet \cite{yu2021onsagernet} and the Laplace neural operator with a convolution kernel \cite{cao2024laplace}, Causality-DeepONet \cite{liu2024causality}, Causal-DeepONet \cite{nghiem2023causal}, Mamba and other state-space models \cite{gu2023mamba,zheng2024alias,hu2025deepomamba,cheng2025mamba,cho2026mbno,guo2026pgmno}, etc. have demonstrated significant potential in reinterpreting relationships governed by differential equations and learning temporal dynamics directly from data.
These developments indicate the effectiveness of incorporating temporal properties to the network structures on improving the accuracy and efficiency.
Nevertheless, works dedicated to time invariance are still rarely found.
The closest example in literature is the Volterra neural network (VNN) in \cite{roheda2024volterra} based on Volterra series expansion.
Despite its theoretical contribution, the discretized multiple integrals are unfortunately infeasible for numerical implementation.

This paper aims at filling the gap for learning operators that respect strict time invariance, rather than learning from data and constraints.
We construct a time-discretized, time-invariant neural operator (TINO) by integrating DeepONet with a time-delay embedding.
This formulation encodes temporal dependencies through delayed inputs while keeping the output basis independent of absolute time.
We then extend this framework to initial-value problems (IVP) for time-dependent PDEs.
Although the observation window is finite and the initial state is prescribed at a fixed time, such problems can sometimes be interpreted as truncated evolutions of an underlying time-invariant system. 
This viewpoint leads to a truncated TINO formulation with spatial POD implementation.
Extensive numerical experiments on operator fitting and on ODE/PDE benchmarks demonstrate that the proposed structures improve prediction accuracy compared with standard DeepONet variants and other baselines.

The rest of this paper is organized as follows.
In \Cref{sect_related}, we formulate time properties including time causality and time invariance, and discuss existing works on neural operators that follow time causality.
In \Cref{sect_method}, a time-invariant neural operator named TINO is derived, as well as a time-truncated version for handling the solution operator of time-invariant PDEs given initial values.
In \Cref{sect_numer}, comprehensive numerical tests are conducted to compare the proposed neural operators and those without time invariance taken into account.
A final conclusion is presented in \Cref{sect_conclusion}.

\section{Prerequisites and related works}\label{sect_related}

\subsection{Time-related properties of operators}

In this section we present a formal mathematical description of time properties of operators studied in this paper.
Let
\begin{equation*}
\mathsf{G}: \mathcal{X} \to \mathcal{Y} = \mathsf{G}[\mathcal{X}]
\end{equation*}
be a time-dependent operator given in the form
\begin{equation}\label{eq_G_demo}
    f(t, \vec{x}) \in \mathcal{X} \mapsto \mathsf{G}[f](t, \vec{y}) \in \mathcal{Y}, \quad t \in [T_{0}, T], \quad \vec{x} \in \Omega_X, \quad \vec{y} \in \Omega_Y,
\end{equation}
where $\vec{x}$ and $\vec{y}$ are spatial variables, and the time interval can be infinite, as $(-\infty, T]$.

\subsubsection{Time causality}
We say the operator $\mathsf{G}$ is \emph{time-causal} if the output values $\mathsf{G}[f](t, \cdot)$ never depend on future information.
In terms of uniqueness, we formulate the definition as follows.
\begin{definition}[Time causality of operator]
We say an operator $\mathsf{G}$ given in \cref{eq_G_demo} is time-causal, if for any $t \in [T_0, T]$, and any functions $f_1, f_2 \in \mathcal{X}$ satisfying
\begin{equation*}
    f_1(s, \vec{x}) = f_2(s, \vec{x}), \quad \forall s \in [T_0, t], \quad \forall \vec{x} \in \Omega_X,
\end{equation*}
there holds
\begin{equation*}
    \mathsf{G}[f_1](t, \vec{y}) = \mathsf{G}[f_2](t, \vec{y}), \quad \forall \vec{y} \in \Omega_Y.
\end{equation*}
\end{definition}

\begin{prop}[Truncation of time-causal operators]
If $\mathsf{G}: f(t, \vec{x}) \mapsto g(t, \vec{y})$ is a time-causal operator for $t \in [T_0, T]$, then, for any $T' \in [T_0, T]$, the operator
\begin{equation*}
    \mathsf{G}': f(t, \vec{x}) \mapsto g(t, \vec{y}), \quad t \in [T_0, T']
\end{equation*}
is also a time-causal operator.
\end{prop}

Consider a PDE with first-order in time with initial values
\begin{align*}
    \partial_t u + \mathsf{L}[u] &= f(t, \vec{x}), \quad t \in [T_0, T], \quad \vec{x} \in \Omega_X, \\
  \mathcal{I}[u(T_0, \cdot)](\vec{x}) &= v(\vec{x}), \quad \vec{x} \in \Omega_X,
  %u(t, \vec{x}) &= w(t, \vec{x}), \quad t \in (T_0, T], \quad \vec{x} \in \partial \Omega_X,
\end{align*}
where $\mathsf{L}$ is a differential operator that does not depend on time.
We omit the boundary value equations in the below discussion,
as they can usually be interpreted as differential equations with spatial masking factors.
In practice, the solution operator $\mathsf{S}: f \mapsto u$ is almost always assumed to be time-causal.

When applying finite difference schemes on ascending time steps $t_0 < t_1 < \cdots$, by letting $u^{n}(\cdot) \approx u(t_n, \cdot)$ and $(u^{n+1} - u^n)/\Delta t_n \approx \partial_t u(t_n, \cdot)$, both explicit and implicit schemes iteratively compute each $u^{n+1}$ using the information of $v(\cdot), f(t_0, \cdot), \cdots, f(t_{n+1}, \cdot)$, i.e., preserving the time causality in the approximation to the solution operator.

In the meantime, deep learning based PDE solvers usually break time causality, unless a similar time marching scheme is applied.
We take the well-known physics-informed neural network (PINN) \cite{raissi2019physics} as a heuristic example.
By using $\mathrm{NN}_u(t, \vec{x}; \vec{\theta})$ as a DNN with trainable variables $\vec{\theta}$ that approximates the PDE solution $u(t, \vec{x})$, PINN minimizes the loss function
\begin{equation}
    L(\vec{\theta}) = \frac{1}{N} \sum_{i=1}^{N} \left\lvert \partial_t \mathrm{NN}_u(t_{(i)}, \vec{x}_{(i)}; \vec{\theta}) + \mathsf{L}[\mathrm{NN}_u](t_{(i)}, \vec{x}_{(i)}; \vec{\theta}) - f(t_{(i)}, \vec{x}_{(i)}) \right\rvert^2 + L_{\mathrm{IC}}(\vec{\theta})
\end{equation}
which consists of the mean of squared PDE residuals at a set of collocation points $(t_{(i)}, \vec{x}_{(i)}) \in [T_0, T] \times \Omega_X$,and the contribution of initial value condition $L_{\mathrm{IC}}(\vec{\theta})$ which does not depend on the right-hand side $f$.
Let
\begin{equation*}
    f_i = f(t_{(i)}, \vec{x}_{(i)}), \quad i = 1, \cdots, N,
\end{equation*}
and treat the DNN parameter vector $\vec{\theta}$ as a function of the inputs $f_1, \cdots, f_N$.
The preservation of time causality requires
\begin{equation*}
    \frac{\partial \mathrm{NN}_u(t_{(i)}, \vec{x}_{(i)}; \vec{\theta})}{\partial f_j} = 0, \quad \forall t_{(i)} < t_{(j)}.
\end{equation*}
However, the governing equation of the DNN parameters $\vec{\theta}$ in PINN only enforces
\begin{equation*}
    \nabla_{\vec{\theta}} L(\vec{\theta}) = 0
\end{equation*}
with additional checks for local minimality, but never distinguishes the order of any pair $t_{(i)}$ and $t_{(j)}$. 
Therefore, PINN in general does not preserve time causality.
A workaround for the numerical side of the above issue in PINN may be found in recent works such as \cite{wang2024respecting,zhao2026casual}, where causal training strategies have been proposed to sequentially minimize the residual loss at ascending time steps.
It remains challenging to keep strict time causality of an entire PDE solver.

\subsubsection{Time invariance}

Time invariance refers to a system yielding the same output over time.
Mathematically, it is equivalent to the commutativity between the operator $\mathsf{G}$ and shift operators.
Since the shifted input and output functions have different time domains, in this paper, we adopt \emph{zero padding} for the shift operator $\mathsf{T}^{\tau}$ for $\tau \in (0, T - T_0)$ as follows:
\begin{equation}
    \mathsf{T}^\tau[h](t) = \begin{cases}
        0, & t \in [ T_0 , T_0 + \tau), \\
        h(t - \tau), & t \in [T_0 + \tau, T].
    \end{cases}
\end{equation}
The padding can be neglected if $T_0 = -\infty$.
\begin{definition}[Time invariance of operator]
We say an operator $\mathsf{G}$ given in \cref{eq_G_demo} is time-invariant, if $\mathsf{G}$ is time-causal, and that for any $\tau \in (0, T - T_0)$ and any $f \in \mathcal{X}$ satisfying $f^\tau = \mathsf{T}^\tau[f] \in \mathcal{X}$, there holds
\begin{equation}
    \mathsf{G}[f^\tau](t, \vec{y}) = \mathsf{G}[f](t - \tau, \vec{y}), \quad \forall t \in [T_0 + \tau, T], \quad \forall \vec{y} \in \Omega_Y.
\end{equation}
\end{definition}

\begin{remark}
The integration operator $h(t) \mapsto \int_{T_0}^{t} h(s) ds$ is time-invariant, because it is time-causal, and that for any $t \in [T_0 + \tau, T]$,
\begin{equation*}
    \int_{T_0}^{t} h^\tau(s) \,\mathrm{d}s = \int_{T_0}^{T_0 + \tau} 0 \,\mathrm{d}s + \int_{T_0 + \tau}^{t} h(s - \tau) \,\mathrm{d}s = \int_{T_0}^{t - \tau} h(s) \,\mathrm{d}s.
\end{equation*}
It is clear that the padding has an impact on the conclusion when the time domain is bounded.
\end{remark}

For time-only continuous time-invariant operators $f(t) \mapsto \mathsf{G}[f](t)$, the Volterra series \cite{sandberg2003volterra,sandberg1982expansions,volterra1887sopra,volterra1931theory} expands the operator in a series expansion of operators similar to the Taylor series
\begin{equation}\label{eq_Volterra_Series}
\mathsf{G} = \mathsf{G}^0 + \mathsf{G}^1 + \mathsf{G}^2 + \cdots,
\end{equation}
where
\begin{equation}
\mathsf{G}^0[f](t) = W^0
\end{equation}
is constant, and 
\begin{equation}\label{eq_Volterra_int}
\mathsf{G}^{n}[f](t) = \int_{[T_0, t]^n} W^n(t - s_1, \cdots, t - s_n) f(s_1) \cdots f(s_n) \,\mathrm{d}s_1 \cdots \,\mathrm{d}s_n, \quad n \ge 1
\end{equation}
are higher orders of the series expansion, with an $n$-variable integral kernel function $W^n$.
A linear continuous time-invariant operator requires only the leading terms of $\mathsf{G}^0$ and $\mathsf{G}^1$, the latter of which is a linear convolution, see e.g. \cite{licurse,li2022approximation} for a strict characterization.
For system identification with random input data, people also adopt the Wiener series \cite{wiener1966nonlinear}, which is an orthogonalized rearrangement of the Volterra series.

Typical difficulties of numerically implementating the Volterra series include the actual convergence of the series expansion, which is not often easily assessed beforehand, and the approximation of the multiple integrals.
Although it is possible to rewrite any truncated Volterra series into a nested formulation similar to the Horner's scheme on polynomials (see e.g. \cite{osowski1994multilayer}), the cost for training such an operator in deep learning is still expensive.
In \cite{roheda2024volterra} where a Volterra series based neural network is proposed, a cascade of second-order Volterra series has been used instead of directly handling higher-order series.

\subsection{Deep neural operators}

The learning of nonlinear operators has been made possible by the capability of parametrizing operators by DNNs.
Bluntly speaking, any DNN based mapping from discretized input functions to discretized output functions can serve as a deep neural operator.
This includes popular network models such as U-Net \cite{ronneberger2015u}, transformer \cite{vaswani2017attention,shih2025transformers}, etc.
However, the deep neural operator networks including DeepONet \cite{jin2022mionet,lu2021learning} and FNO \cite{lifourier} play a central role in the study of operator learning, as they consist of intrinsic structures dedicated to continuous operators. 

\subsubsection{Deep Operator Network (DeepONet)}\label{sect_don}

The DeepONet \cite{jin2022mionet,lu2021learning} is a deep neural operator network that approximates nonlinear continuous operators using data sampled from fixed locations (known as sensors).
The framework of DeepONet originates in the universal approximation theory of continuous functionals and continuous operators \cite{chen1993approximations,chen1995universal}.
DeepONet approximates a continuous operator $\mathsf{G}: f(\vec{x}) \in \mathcal{X} \mapsto \mathsf{G}[f](\vec{y}) \in \mathsf{G}[\mathcal{X}]$ using the formulation
\begin{equation}\label{eq_deepont}
\mathsf{G}[f](\vec{y}) \approx \overrightarrow{\mathrm{NN}}_{\mathrm{branch}}\left( f(\vec{x}_1), \cdots, f(\vec{x}_M) \right) \cdot \overrightarrow{\mathrm{NN}}_{\mathrm{trunk}}(\vec{y}) = \sum_{j=1}^{J} \mathrm{NN}^{(j)}_{\mathrm{branch}}\left( f(\vec{x}_1), \cdots, f(\vec{x}_M) \right) \mathrm{NN}^{(j)}_{\mathrm{trunk}}(\vec{y}),
\end{equation}
where $\mathcal{X}$ is a function space, $\overrightarrow{\mathrm{NN}}_{\mathrm{branch}} = (\mathrm{NN}^{(1)}_{\mathrm{branch}}, \cdots, \mathrm{NN}^{(J)}_{\mathrm{branch}})$ and $\overrightarrow{\mathrm{NN}}_{\mathrm{trunk}} = (\mathrm{NN}^{(1)}_{\mathrm{trunk}}, \cdots, \mathrm{NN}^{(J)}_{\mathrm{trunk}})$ are DNNs that are called the branch net and the trunk net, respectively.
The points $\vec{x}_1, \cdots, \vec{x}_M$ are fixed points chosen from the domain, known as sensors.

Once a DeepONet is trained for the approximation of an operator $\mathsf{G}$, it approximates every output function in $\mathsf{G}[\mathcal{X}]$ using the same trunk net, effectively assuming that the output function space $\mathsf{G}[\mathcal{X}]$ is linearly approximated by the $J$ basis functions $\mathrm{NN}^{(1)}_{\mathrm{trunk}}(\vec{y}), \cdots, \mathrm{NN}^{(J)}_{\mathrm{trunk}}(\vec{y})$.
Based upon this fact, \cite{lu2022comprehensive} proposed POD-DeepONet as an alternative version to DeepONet to improve the approximation quality of $\mathsf{G}[\mathcal{X}]$, where the basis functions of $\mathsf{G}[\mathcal{X}]$ are chosen by proper orthogonal decomposition (POD) from training data, instead of a trunk DNN trained together with the branch net.
Namely, let $\{\mathsf{G}[f_i](\vec{y})\}_{i=1}^{N_{\mathrm{data}}}$ be the snapshot set of target output functions from training data.
The POD-DeepONet approximates $\mathsf{G}$ using a mean-plus-correction expansion
\begin{equation}\label{eq_pod_deepont}
\mathsf{G}[f](\vec{y}) \approx \phi_0(\vec{y}) + \sum_{j=1}^{J}\mathrm{NN}^{(j)}_{\mathrm{branch}}\left( f(\vec{x}_1), \cdots, f(\vec{x}_M) \right) \phi_j(\vec{y}),
\end{equation}
where
\begin{equation*}
    \phi_0(\vec{y}) = \frac{1}{N}\sum_{i=1}^{N_{\mathrm{data}}}\mathsf{G}[f_i](\vec{y}),
\end{equation*}
and the basis functions $\phi_1(\vec{y}), \cdots, \phi_J(\vec{y})$ are extracted via POD, such as the singular value decomposition (SVD) on the set of functions $\{\mathsf{G}[f_i](\vec{y}) - \phi_0(\vec{y})\}_{i=1}^{N_{\mathrm{data}}}$ discretized on a sequence of collocation points.

\subsubsection{Fourier Neural Operator (FNO)}

The Fourier Neural Operator (FNO) \cite{lifourier} is known as another parallel series of deep neural operator network. FNO approximates a continuous operator $\mathsf{G}: f(\vec{x}) \in \mathcal{X} \mapsto \mathsf{G}[f](\vec{x}) \in \mathsf{G}[\mathcal{X}]$ by composing a sequence of Fourier layers via the Fast Fourier Transform (FFT). Each Fourier layer applies a nonlinear transformation to a combination of local linear mappings and global spectral convolutions. The update at the $\ell$-th layer is given by
\begin{equation*}
    v_{\ell+1}(\vec{x}) = \sigma\left( W_\ell v_\ell(\vec{x}) + c_\ell + \mathcal{F}^{-1}\left( R_\ell \cdot \mathcal{F}(v_\ell) \right)(\vec{x}) \right),
\end{equation*}
The affine term $W_\ell v_\ell(\vec{x}) + c_\ell$ combined with the pointwise activation $\sigma$ acts as a local channel-mixer, aiming at capturing fine-grained spatial features and performing nonlinear transformations independently at each grid point.
The parallel application of FFT $\mathcal{F}$ shifts the hidden representation into the frequency domain.
The element-wise multiplication with the trainable matrix $R_\ell$ applies a trainable spectral filter.
Finally, the inverse Fourier transform $\mathcal{F}^{-1}$ maps the filtered spectrum back to physical space.

The Fourier transformation is deliberately employed to convert domain-wide spatial correlations into frequency modes, thereby exposing long-range patterns that are otherwise costly to model in physical space. Moreover, since high-frequency information decays faster in the frequency domain, FNO is less susceptible to interference from high-frequency noise during training.

FNO's FFT is inherently global and a single layer can  information from all points in the domain. As a result, it naturally violates causal ordering unless explicitly constrained.

\subsubsection{Causality-DeepONet}
In \cite{liu2024causality}, a causality enhanced DeepONet, namely Causality-DeepONet, was proposed for effective prediction of the seismic response of buildings where a thorough simulation is rather expensive.
The branch nets of Causality-DeepONet adopt time sensors in a chronological order.
With temporal sensors
\begin{equation*}
s_i = T_0 + i h, \quad i = 0, \cdots, m = \frac{T - T_0}{h},
\end{equation*}
the operator network takes the form
\begin{equation}\label{eq_causality_don}
\mathsf{G}[f](t) \approx \sum_{j=1}^{J} \sum_{k=1}^{K} c_k^j \sigma_{\mathrm{b}} \left( \sum_{i=1}^{\left\lceil \frac{t - T_0}{h}\right\rceil} \xi_{k, m - \left\lceil \frac{t - T_0}{h}\right\rceil + i}^{j} f(s_i) \right) \sigma_{\mathrm{t}}(\omega_j t + \zeta_j), \quad t \in [T_0, T],
\end{equation}
where $c_k^j$, $\xi_{k, i}^{j}$, $\omega_k$ and $\zeta_k$ are trainable network parameters, independent of the input variable $t$, and $\sigma_{\mathrm{b}}$ and $\sigma_{\mathrm{t}}$ are activation functions for the branch net and the trunk net, respectively.

\subsubsection{Time-causal DeepONet}

With a straightforward causal masking (c.f. \cite{vaswani2017attention}), DeepONet can achieve time causality.
Below we formulate such neural operator, mainly for reference by the rest of this paper.
The proposal is spiritually close to Causal-DeepONet \cite{nghiem2023causal}.
Given a grid of spatiotemporal sensors
\begin{equation*}
    (t_i, \vec{x}_m), \quad i=1,\cdots,I, \quad m=1,\cdots,M,
\end{equation*}
causal masks are embedded in the branch nets as follows:
\begin{equation}\label{eq_TC_DON}
    \mathsf{G}[f](t, \vec{y}) \approx \overrightarrow{\mathrm{NN}}_{\mathrm{branch}} \left( \vec{f}_{\le t}(\vec{x}_1), \cdots, \vec{f}_{\le t}(\vec{x}_M) \right) \cdot \overrightarrow{\mathrm{NN}}_{\mathrm{trunk}}(t, \vec{y}),
\end{equation}
where each masked vector
\begin{equation*}
    \vec{f}_{\le t} (\vec{x}_m) = \left( f_{\le t}(t_1, \vec{x}_m), \cdots, f_{\le t}(t_I, \vec{x}_m) \right), \quad 
    f_{\le t}(s, \vec{x}) = \begin{cases}
        f(s, \vec{x}), & T_0 \le s \le t, \\
        0, & \mathrm{otherwise}.
    \end{cases}
\end{equation*}
We name the neural operator in \cref{eq_TC_DON} \textbf{time-causal DeepONet (TC-DeepONet)}.

Note that the branch nets of time-causal DeepONet in \cref{eq_TC_DON} (as well as \cite{nghiem2023causal}) differ from Causality-DeepONet, due to the different arrangement of the $\xi_{k,i}^{j}$ network parameters, and the bias parameter of the input layer.

\subsubsection{Neural operators on state-space models}

The state-space model (SSM) has been paid close attention to because of the recent development of Mamba architecture \cite{gu2023mamba} and its subsequent works in scientific machine learning such as \cite{zheng2024alias,hu2025deepomamba,cheng2025mamba,cho2026mbno,guo2026pgmno}, to name a few.
The SSM assumes there is a hidden state $u(t)$ with dynamics
\begin{align*}
u'(t) &= \phi(u(t), f(t), t), \\
\mathsf{G}[f](t) &= \psi(u(t), f(t), t),
\end{align*}
where $\phi$ and $\psi$ are functions that do not depend on the input function $f(t)$ or the state $u(t)$.
Time invariance may be achieved when $\psi$ and $\psi$ do not explicitly depend on the input variable $t$.
\cite{paduart2010identification} applies polynomials to approximate the nonlinear functions; meanwhile \cite{suykens1995nonlinear} uses shallow neural networks.
In Mamba, for performance concerns, the SSM is simplified so that the hidden state $u(t)$ keeps linear dynamics.

\section{Time-invariant neural operators}\label{sect_method}

\subsection{Time-delay neural network (TDNN)}

We begin with a continuous time-invariant operator $\mathsf{G}$ that maps 1-D functions of time on $t \in [T_0, T]$ to 1-D functions of time on the same interval.
Suppose the Volterra series expansion of $\mathsf{G}$ (c.f. \cref{eq_Volterra_Series}) converges.
With a time grid $t_i = T_0 + i(T- T_0) / I$, $i=0,\cdots,I$, denoting $f_i = f(t_i)$, we consider the discretization of each integral
\begin{align}\label{eq_td_deriv}
\begin{split}
G^n(t_i) &= \int_{[T_0, t_i]^n} W^n(t_i - s_1, \cdots, t_i - s_n) f(s_1) \cdots f(s_n) \,\mathrm{d}s_1 \cdots \,\mathrm{d}s_n \\
&\approx \sum_{i_1, \cdots, i_n = 0}^{i - 1} \widehat{W}^n_{i - i_1, \cdots, i - i_n} f_{i_1} \cdots f_{i_n} \\
&= \sum_{j_1, \cdots, j_n = 1}^{i} \widehat{W}^n_{j_1,\cdots,j_n} f_{i - j_1} \cdots f_{i - j_n} \\
&= \sum_{j_1, \cdots, j_n = 1}^{i} \widehat{W}^n_{j_1,\cdots,j_n} ( \vec{f} * \vec{e}_{j_1} )_i \cdots ( \vec{f} * \vec{e}_{j_n} )_i,
\end{split}
\end{align}
where
\begin{equation*}
    \vec{f} = (f_0, \cdots, f_I),
\end{equation*}
and vectors $\vec{e}_1, \cdots, \vec{e}_I$ are unit vectors in $\mathbb{R}^{I+1}$, such that
\begin{equation*}
    ( \vec{f} * \vec{e}_{j} )_i = \begin{cases}
        f_{i - j}, & j \le i, \\
        0, & \mathrm{otherwise}.
    \end{cases}
\end{equation*}
Indeed, the above discretization of the integral $G^n(t_i)$ is an $n$-th order \emph{homogeneous} polynomial of $(\vec{f} \ast \vec{e}_{1})_i, \cdots, (\vec{f} \ast \vec{e}_{i})_i$.
By approximating polynomials with DNNs and adding up for all $n \ge 0$, the resulting network can be written as a function
\begin{equation}
    \mathsf{G}[f](t_i) \approx \mathrm{NN}_{\mathrm{TD}}\left( (\vec{f} \ast \vec{e}_{1})_i, \cdots, (\vec{f} \ast \vec{e}_{I})_i \right),
\end{equation}
where the network $\mathrm{NN}_{\mathrm{TD}}$ is known as the time-delay neural network (TDNN) \cite{waibel1989phoneme}.
The above heuristic derivation had already been discussed in the literature, see e.g. \cite{marmarelis1994relation, wray1994calculation}.
Compared to the Volterra series with polynomial backbone, the DNN formulation further extends the scope of representation.

\subsection{Time-invariant neural operators (TINO)}

We then extend the above discussion to the case where the inputs are 1-D temporal functions on $[0, T]$, while the output functions have an additional spatial variable $\vec{y} \in \Omega_Y$.
An extension to the original Volterra series \cref{eq_Volterra_Series} becomes parametrized by $\vec{y}$, namely,
\begin{equation}\label{eq_Volterra_ty}
    \mathsf{G}[f](t,\vec{y}) = W^0(\vec{y}) + \sum_{n=1}^{\infty} \int_{[T_0,t]^n} W^n(t-s_1, \cdots, t-s_n; \vec{y}) f(t - s_1) \cdots f(t - s_n) \,\mathrm{d}s_1 \cdots \,\mathrm{d}s_n.
\end{equation}
The corresponding temporal discretization on the $n$-th integral term becomes
\begin{equation*}
    G^n(t_i; \vec{y}) \approx \sum_{j_1, \cdots, j_n = 1}^{i} \widehat{W}^n_{j_1,\cdots,j_n}(\vec{y}) ( \vec{f} * \vec{e}_{j_1} )_i \cdots ( \vec{f} * \vec{e}_{j_n} )_i.
\end{equation*}
That is, the discretized integral kernels $\widehat{W}^n$ carry spatial information of $\vec{y}$, and temporal information is implicitly recorded in the convolutions.

Recall that in DeepONet, the trunk net uses $J$ basis functions to linearly approximate any output function of the operator.
Let these \emph{spatial} basis functions be denoted by $\omega_1(\vec{y}), \cdots, \omega_J(\vec{y})$, and use them as the (approximated) linear basis in the above expansion.
The coefficients are again homogeneous $n$-th order polynomials of $(\vec{f} \ast \vec{e}_{1})_i, \cdots, (\vec{f} \ast \vec{e}_{i})_i$.
By further replacing polynomials and basis functions by DNNs, and adding up for all $n \ge 0$, we arrive at a neural operator formulation
\begin{equation}\label{eq_TINO_ty}
    \mathsf{G}[f](t_i, \vec{y}) \approx \sum_{j=1}^{J} {\mathrm{NN}}^{(j)}_{\mathrm{branch}}\left( (\vec{f} \ast \vec{e}_{1})_i, \cdots, (\vec{f} \ast \vec{e}_{I})_i \right) \mathrm{NN}^{(j)}_{\mathrm{trunk}}(\vec{y}).
\end{equation}
In \cref{eq_TINO_ty}, each ${\mathrm{NN}}^{(j)}_{\mathrm{branch}}$ refers to a component of the output of a branch net like in DeepONet, while the trunk net counterpart with basis functions $\mathrm{NN}^{(j)}_{\mathrm{trunk}}$ does \emph{not} depend on time $t$.

When the input function has additional spatial variable $\vec{x} \in \Omega_X$, we employ spatial sensors at $\vec{x}_1, \cdots, \vec{x}_M$, and replace each entry of $\vec{f}$ in \cref{eq_TINO_ty} by a vector on these spatial sensors.
Namely,
\begin{equation}\label{eq_TINO_txy}
\mathsf{G}[f](t_i, \vec{y}) \approx \sum_{j=1}^{J} {\mathrm{NN}}^{(j)}_{\mathrm{branch}}\left( \left\{ (\vec{f}^m * \vec{e}_{i'})_i \right\}_{1 \le m \le M, 1 \le i' \le I}  \right) \mathrm{NN}^{(j)}_{\mathrm{trunk}}(\vec{y}),
\end{equation}
where $\vec{f}^m = (f(t_0, \vec{x}_m), \cdots, f(t_I, \vec{x}_m))$.

We refer to \cref{eq_TINO_txy} as a (time discretized) \textbf{time-invariant neural operator (TINO)}, and treat \cref{eq_TINO_ty} as its special case.

The ${\mathrm{NN}}^{(j)}_{\mathrm{branch}}$ functions may be joined into a fully-connected neural network.
We emphasize that, when doing so, the input layer can be further reformulated as a convolution layer with trainable variables.

\begin{figure}[h]
    \centering
    \includegraphics[keepaspectratio, clip, trim={0 0 0 0}, pagebox=cropbox, width=\linewidth]{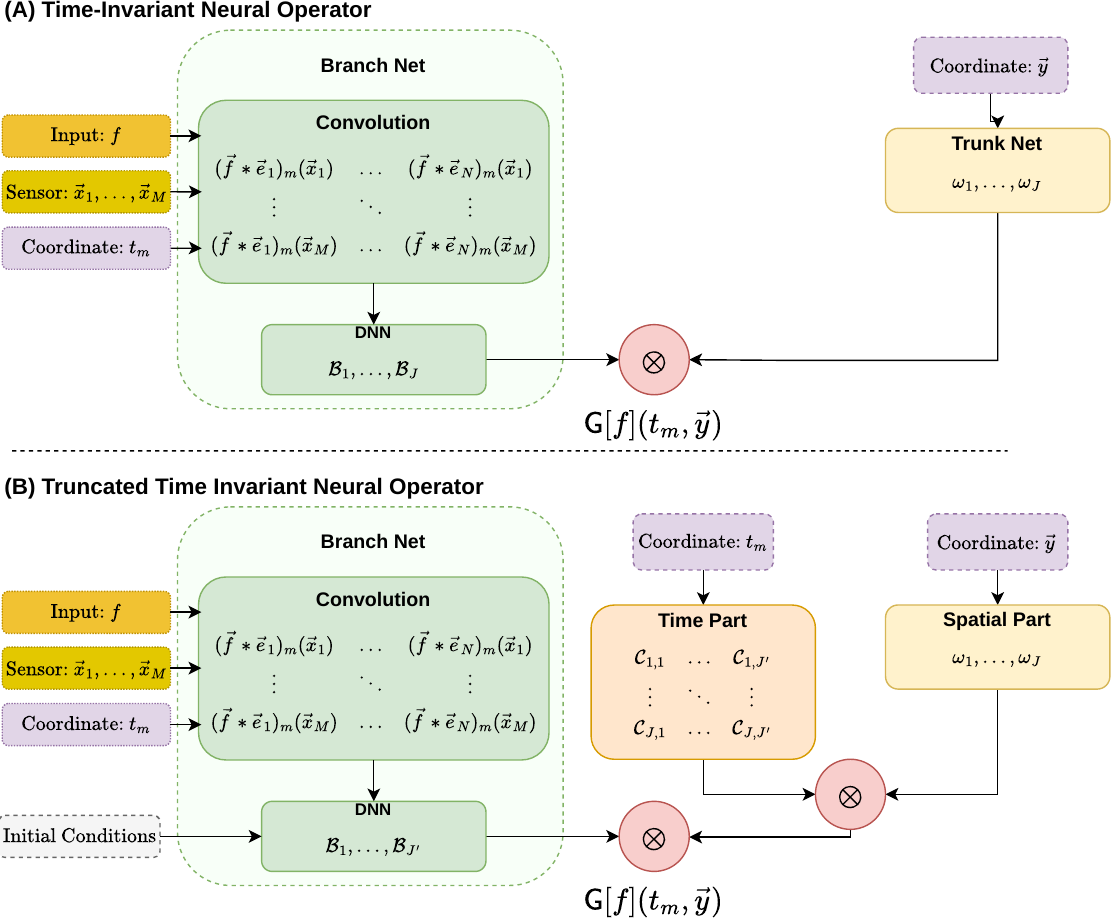}
    \caption{Architectures of TINO and TrTINO}
    \label{TINO and TrTINO}
\end{figure}

\subsubsection{POD-TINO}

The POD versions of the above TINOs (\cref{eq_TINO_ty,eq_TINO_txy}) are given by performing \emph{spatial} POD to replace the trunk net functions $\mathrm{NN}^{(j)}_{\mathrm{trunk}}(\vec{y})$.

\subsection{Truncated time-invariant neural operators (TrTINO)}

The dynamics of a time-invariant physical system may have started from an initial time $t = T_- < T_0$ before our observation at $T_0$.
Suppose that instead of the past, an initial condition (IC) as supplementing information at $t=T_0$ is known.
Here, we derive the time truncated version of TINO for this situation.

For simplicity, we consider such an operator $\mathsf{G}$ that can be extended to a mapping from a temporal function $f(t)$ defined on $t \in [T_-, T]$ to a spatiotemporal function $\mathsf{G}[f](t, \vec{y})$, with a converging Volterra series expansion
\begin{equation}
    \mathsf{G}[f](t, \vec{y}) = W^0(\vec{y}) + \sum_{n=1}^{\infty} \int_{[T_-, t]^n} W^n(t - s_1, \cdots, t - s_n; \vec{y}) f(s_1) \cdots f(s_n) \,\mathrm{d}s_1 \cdots \,\mathrm{d}s_n.
\end{equation}
From the $n$-th integral, denoted by $G^n(t; \vec{y})$, we split the integral domain $[T_-, t]^n$ into $2^n$ rectangular pieces by a tensor product of the division of interval $[T_-, t] = [T_-, T_0) \cup [T_0, t]$ on each dimension.
Such decomposition can be denoted by
\begin{equation*}
G^n(t; \vec{y}) = \sum_{\Lambda \subset \{1,2,\cdots,n\}} G^{n, \Lambda}(t; \vec{y}),    
\end{equation*}
where
\begin{equation*}
    G^{n, \Lambda}(t; \vec{y}) = \int_{D_1 \times \cdots \times D_n} W^n(t-s_1, \cdots, t - s_n; \vec{y}) f(s_1) \cdots f(s_n) \,\mathrm{d}s_1 \cdots \,\mathrm{d}s_n, \quad D_i = \begin{cases}
        [T_0, t], & i \in \Lambda, \\
        [T_-, T_0), & \mathrm{otherwise}.
    \end{cases}
\end{equation*}
Since for any $i \notin \Lambda$, the information of $f(s_i)$ in the integrand of $G^{n,\Lambda}(t; \vec{y})$ is not accessible, we integrate them and get
\begin{equation*}
    G^{n, \Lambda}(t; \vec{y}) = \int_{D_{l_1} \times \cdots \times D_{l_r}} W^{n,\Lambda}(t-s_{l_1}, \cdots, t-s_{l_r}; t ; \vec{y}) f(s_{l_1}) \cdots f(s_{l_r}) \,\mathrm{d}s_{l_1}\cdots \,\mathrm{d}s_{l_r}, \quad \{l_1,\cdots,l_r\} = \Lambda,
\end{equation*}
where the function $W^{n,\Lambda}$ contains \emph{all} the information of $f$ from the past $[T_-, T_0)$.
When the integral $G^{n, \Lambda}$ is further discretized on a time grid $t_i = T_0 + i(T - T_0)/I$, similar to the derivation in \cref{eq_td_deriv},
\begin{equation}
\begin{aligned}
     G^{n, \Lambda}(t_i; \vec{y}) &\approx \sum_{i_1,\cdots,i_r = 0}^{i-1} \widehat{W}^{n,\Lambda}_{i-i_1,\cdots,i-i_r}(t_i;\vec{y}) f_{i_1} \cdots f_{i_r}\\ 
     &= \sum_{j_1,\cdots,j_r = 1}^{m} \widehat{W}^{n,\Lambda}_{j_1,\cdots,j_n}(t_i;\vec{y}) (\vec{f} \ast \vec{e}_{j_1})_i \cdots (\vec{f} \ast \vec{e}_{j_r})_i,
\end{aligned}
\end{equation}
which is again a polynomial of $(\vec{f} \ast \vec{e}_{1})_i, \cdots, (\vec{f} \ast \vec{e}_{I})_i$, with coefficients depending on $t$, and containing information from the past input of $f$ on $[T_-, T_0)$.
Then, similar to the proposal of TINO, with the same spatial basis functions $\omega_1(\vec{y}), \cdots, \omega_J(\vec{y})$, we can split the temporal and spatial parts
\begin{equation*}
    \widehat{W}^{n,\Lambda}_{j_1,\cdots,j_n}(t;\vec{y}) \approx \sum_{j=1}^{J} c_{j;j_1,\cdots,j_n}^{n,\Lambda}(t) \omega_j(\vec{y}).
\end{equation*}
By summing up for any $\Lambda, n$ and replacing with DNN approximation, we get an (intermediate) approximation
\begin{equation*}
    \mathsf{G}[f](t_i, \vec{y}) \approx \sum_{j=1}^{J} \mathrm{NN}^{(j)}\left( (\vec{f} \ast \vec{e}_{1})_i, \cdots, (\vec{f} \ast \vec{e}_{I})_i; t_i \right) \omega_j(\vec{y}).
\end{equation*}
Now that $\mathrm{NN}^{(j)}$ has both information of the input function $f$ from observation (by the convolutions) and from the past (by variable $t$), where the latter is usually provided by the initial conditions, and noting that $f$ and $t$ typically have different units, we further split the networks $\mathrm{NN}^{(j)}$ into a functional part multiplied by a time part, leading to a neural operator formulation
\begin{equation}\label{eq_TrTINO_ty}
    \mathsf{G}[f](t_i, \vec{y}) \approx \sum_{j'=1}^{J'} \sum_{j=1}^{J} \mathrm{NN}_{\mathrm{branch}}^{(j')} \left( (\vec{f} \ast \vec{e}_{1})_i, \cdots, (\vec{f} \ast \vec{e}_{I})_i; \mathrm{IC} \right) \mathrm{NN}_{\mathrm{t}}^{(j',j)}(t) \mathrm{NN}_{\mathrm{branch}}^{(j)}(\vec{y}),
\end{equation}
where ``IC'' refers to the initial value information that is passed to the branch net.

When the input function $f$ also takes spatial variable $\vec{x} \in \Omega_X$, the corresponding proposal is accordingly using sensors at $\vec{x}_1, \cdots, \vec{x}_M$, by
\begin{equation}\label{eq_TrTINO_txy}
    \mathsf{G}[f](t_i, \vec{y}) \approx \sum_{j=1}^{J} {\mathrm{NN}}^{(j)}_{\mathrm{branch}}\left( \left\{ (\vec{f}^m * \vec{e}_{i'})_i \right\}_{1 \le m \le M, 1 \le i' \le I}; \text{IC}  \right) \mathrm{NN}_{\mathrm{t}}^{(j',j)}(t) \mathrm{NN}_{\mathrm{branch}}^{(j)}(\vec{y}),
\end{equation}
effectively separating temporal and spatial variables $t$, $\vec{x}$ and $\vec{y}$ into distinct components.

We name the above neural operator \cref{eq_TrTINO_txy} \textbf{truncated time-invariant neural operator (TrTINO)}, and treat \cref{eq_TrTINO_ty} as its special case.
When the output functions do not depend on spatial variable $\vec{y}$, we set $J = 1$ and let the branch net be $1$, a trivial constant.

\subsubsection{TrTINO with spatial POD}

The POD versions of the above TrTINOs \cref{eq_TrTINO_ty,eq_TrTINO_txy}, named \textbf{SPOD-TrTINO}, are given by performing \emph{spatial POD} to replace the trunk net.
Compared to a spatiotemporal POD, the spatial version is analogous to the formulation of TrTINO \cref{eq_TrTINO_txy}.
Since the POD is practically taken by heavy SVD algorithms, reducing the temporal dimension improves overall performance, meanwhile keeping a time network $\overrightarrow{\mathrm{NN}}_{\mathrm{t}}$ is indeed not expensive in training, as suggested by U-Net and other architectures with time embedding mechanism.

In our implementation, for evaluating the spatial POD, we pick spatial functions from the dataset of $\mathsf{G}[\cdot](\tau, \cdot)$ with
\begin{equation*}
    \tau = T_0 + \frac{1}{4}(T - T_0),
\end{equation*}
i.e. at the first quartile of the time domain.
The SVD on these functions vectorized at collocation points is truncated by discarding singular values $\sigma_i$ that satisfy
\begin{equation}\label{eq_SVD_eta}
    \frac{\sigma_i}{\sigma_{\max}} < 10^{-6}.
\end{equation}
where $\sigma_{\max}$ is the largest singular value.

We also apply the above spatial POD implementation to DeepONet and TC-DeepONet, respectively.
The resulting formulations are given as follows, where $\phi_0(\vec{y})$ and $\phi_1(\vec{y}), \cdots, \phi_J(\vec{y})$ are derived in the same manner.

\textbf{SPOD-TrTINO} (c.f. \cref{eq_TrTINO_txy}):
\begin{equation}\label{eq_spod_trtino}
    \mathsf{G}[f](t_i, \vec{y}) \approx \phi_0(\vec{y}) + \sum_{j=1}^{J} {\mathrm{NN}}^{(j)}_{\mathrm{branch}}\left( \left\{ (\vec{f}^m * \vec{e}_{i'})_i \right\}_{1 \le m \le M, 1 \le i' \le I}; \text{IC}  \right) \mathrm{NN}_{\mathrm{t}}^{(j',j)}(t) \phi_j(\vec{y}).
\end{equation}

\textbf{SPOD-DeepONet} (c.f. \cref{eq_deepont}):
\begin{equation}\label{eq_spod_don}
    \mathsf{G}[f](t, \vec{y}) \approx \phi_0(\vec{y}) + \sum_{j'=1}^{J'} \sum_{j=1}^{J} {\mathrm{NN}}^{(j')}_{\mathrm{branch}} \left( \vec{f}(\vec{x}_1), \cdots, \vec{f}(\vec{x}_M) \right) \mathrm{NN}_{\mathrm{t}}^{(j',j)}(t) \phi_j(\vec{y}),
\end{equation}
where each $\vec{f}(\vec{x_m}) = \left( f(t_1, \vec{x}_m), \cdots, f(t_I, \vec{x}_m) \right)$.

\textbf{TC-SPOD-DeepONet} (c.f. \cref{eq_TC_DON}):
\begin{equation}\label{eq_TC_SPOD_DON}
\mathsf{G}[f](t, \vec{y}) \approx \phi_0(\vec{y}) + \sum_{j'=1}^{J'} \sum_{j=1}^{J} {\mathrm{NN}}^{(j')}_{\mathrm{branch}} \left( \vec{f}_{\le t}(\vec{x}_1), \cdots, \vec{f}_{\le t}(\vec{x}_M) \right) \mathrm{NN}_{\mathrm{t}}^{(j',j)}(t) \phi_j(\vec{y}).
\end{equation}

\subsection{Remarks}

The truncated TINO operators does not keep the time invariance properties, due to the involvement of initial conditions at the fixed time $t=T_0$.
When treated as a truncated version of a TINO operator, in terms of network structure, the major difference between TrTINO and DeepONet, or between their POD variants, consist of the one-layer temporal convolution on the input, as well as a separation of space and time in the design of the trunk net.
Readers may also refer to \cite{wu2024capturing} for a DeepONet variant with branch nets implemented by convolutional neural networks.
\section{Numerical Tests}\label{sect_numer}

In this section, we conduct comprehensive numerical tests for comparison of the models discussed above, including:
\begin{itemize}
    \item DeepONet (DON, \cref{eq_deepont}), its POD version (POD-DON, \cref{eq_pod_deepont}) and its spatial POD version (SPOD-DON, \cref{eq_spod_don});
    \item Time-causal DeepONet (TC-DON, \cref{eq_TC_DON}), and its spatial POD version (TC-SPOD-DON, \cref{eq_TC_SPOD_DON});
    \item Time-invariant neural operator (TINO, \cref{eq_TINO_ty,eq_TINO_txy}), its time truncated version (TrTINO, \cref{eq_TrTINO_ty,eq_TrTINO_txy}), and the version with spatial POD (SPOD-TrTINO, \cref{eq_spod_trtino});
    \item FNO and U-Net;
\end{itemize}
where FNO and U-Net are not specially modified for time properties, and require output functions having the same domain as input functions.
We categorize the test problems into $3$ groups:
\begin{itemize}
    \item The fitting of time-invariant operators, in \Cref{sect_fitting_ti_op};
    \item Solving time-dependent ODEs and ODE systems with initial values, in \Cref{sect_sol_ode};
    \item Solving time-dependent PDEs with initial values, in \Cref{sect_sol_pde}.
\end{itemize}
In the last group, we also include a test case where the solution operator is time-causal but not time-invariant.

\Cref{figure_numer_result} shows the relative $L^2$ errors across all numerical experiments under different methods.

\subsection{Experimental setup}
To ensure a fair comparison, we standardize the configurations across all models. For each benchmark, we generate $1000$ training and $200$ testing function pairs using identical sampling procedures (except a special treatment in \Cref{subsubsect_navier_stokes_icvar}).
All the spatial and temporal grids are uniform grids.
The models are trained and tested on a single NVIDIA A40 GPU with 48GB video RAM. In the case of ODEs and 1-D PDEs, the training batch size is set to 200, and for PDEs of higher dimensions, it is reduced to 100(except a special treatment in \Cref{subsubsect_navier_stokes_icvar}).

Detailed network parameters and data generation procedures are provided in \cref{app_para_set}.

During the training phase, all neural operators are optimized by minimizing the Mean Squared Error (MSE) between the predicted solutions and the ground truth over the training dataset. Let $N_{\text{train}}$ denote the total number of training samples, each discretized on $N$ grid points. The training loss function $\mathcal{L}_{\text{train}}$ is formulated as:
\begin{equation} \label{eq_loss_func}
    \mathcal{L}_{\text{train}} = \frac{1}{N_{\text{train}}} \sum_{k=1}^{N_{\text{train}}} \left( \frac{1}{N} \sum_{i=1}^{N} \left( u_i^{(k)} - \hat{u}_i^{(k)} \right)^2 \right),
\end{equation}
where $u_i^{(k)}$ and $\hat{u}_i^{(k)}$ represent the ground truth and the predicted value at the $i$-th grid point of the $k$-th training sample, respectively.

After training, we evaluate the models using the Mean Squared Error (MSE) and its Standard Deviation (SD) across the entire test set. For the $k$-th test sample (out of $N_{\text{test}}$ total test samples), the sample-wise MSE is computed as:
\begin{equation} \label{eq_single_mse}
    \text{MSE}^{(k)} = \frac{1}{N} \sum_{i=1}^{N} \left( u_i^{(k)} - \hat{u}_i^{(k)} \right)^2.
\end{equation}
The reported results are the mean and standard deviation of $\text{MSE}^{(k)}$ over all $N_{\text{test}}$ test samples, calculated as:
\begin{align} \label{eq_mse_and_sd}
    \overline{\text{MSE}} &= \frac{1}{N_{\text{test}}} \sum_{k=1}^{N_{\text{test}}} \text{MSE}^{(k)}, \\
    \text{SD} &= \sqrt{ \frac{1}{N_{\text{test}} - 1} \sum_{k=1}^{N_{\text{test}}} \left( \text{MSE}^{(k)} - \overline{\text{MSE}} \right)^2 }.
\end{align}

Furthermore, to ensure a fair and rigorous comparison, we also evaluate the performance of all models using standard error metrics. For a single test sample discretized on $N$ grid points, the mean Relative $L^2$ Error of $N_{\text{test}}$ total test samples is defined as:
\begin{equation} \label{eq_rela_l2}
    \text{Relative } L^2 \text{ Error} = \frac{1}{N_{\text{test}}}\sum_{k=1}^{N_{\text{test}}}\frac{\| u - \hat{u} \|_{L^2}}{\| u \|_{L^2}} =  \frac{1}{N_{\text{test}}}\sum_{k=1}^{N_{\text{test}}}\sqrt{ \frac{\sum_{i=1}^{N} (u_i^{(k)} - \hat{u}_i^{(k)})^2}{\sum_{i=1}^{N}( u_i^{(k)})^2} }.
\end{equation}
This metric is used for the overall error distribution analysis (e.g., \cref{figure_numer_result}).

\subsection{Fitting time-invariant operators}\label{sect_fitting_ti_op}

\subsubsection{Time-only operator} \label{subsubsect_fitting_time}
We focus  on the fitting of the following operator:
\begin{equation}
    \mathsf{G}: f(t) \in C[0,T] \mapsto f(t)^2 + \cos\left((10-t)f(t)\right) + \int_{0}^{t}f(s)\,\mathrm{d}s \in C[0, T].
\end{equation}

In this experiment, we set $T=10$.
Input functions are restricted by the parameterization
\begin{equation*}
    f(t) = A \sin(b t + c)
\end{equation*}
with $A \in [0.5, 2]$, $b \in [0.1, 2]$, $c \in [0, 2\pi]$.
The temporal grid has $n_{\mathrm{t}} = 1000$ points.

\subsubsection{Time-dependent operator in 1-D space}\label{subsubsect_fitting_time_space}

We focus on the fitting of the following operator in one spatial and one temporal dimension:
\begin{equation}
    \mathsf{G}: f(t,x) \in C([0,T]\times [0,L]) \mapsto \frac{1}{2}\int_0^t \int_0^y K\left(\frac{z}{2},\frac{y}{2}\right) u\left(s,\frac{z}{2}\right)^2 \,\mathrm{d}z\,\mathrm{d}s \in C([0,T]\times [0,2L]).
\end{equation}

In this experiment, we set $T=4$ and $L=1$.
The temporal grid has $n_\mathrm{t} = 200$ points.
The spatial grids have $n_{\mathrm{in}} = 128$ points for input and $n_{\mathrm{out}} = 128$ for output.
The integral kernel function $K(x,y) = xy + 1$.
Input functions are restricted by the parametrization
\begin{equation*}
    f(t, x) = A\sin(bt+c) e^{-(x-\mu)^2/(2\sigma^2)}
\end{equation*}
with $A \in [0.1, 1]$, $b \in [0.2, 1]$, $c \in [0, 2\pi]$, $\mu \in [L/4, 3L/4]$, $\sigma \in [0.1L, 0.4L]$.

\subsubsection{Results}
\begin{table}[h]
\centering
\caption{Error of neural operators on purely temporal and spatio-temporal fitting tasks. Results are reported as the mean squared error (MSE) $\pm$ standard deviation (SD).}
\label{tab:mse_results_fitting}
\begin{tabular}{lcc}
\toprule
 & \Cref{subsubsect_fitting_time}
 & \Cref{subsubsect_fitting_time_space} \\ 
%\midrule
 & Fitting (time) & Fitting (time and 1-D space) \\
\midrule
FNO & \textbf{2.20e-2 $\pm$ 7.01e-2}  & - \\

\hdashline
DON & 1.20e-1 $\pm$ 8.43e-2 & 3.92e-5 $\pm$ 9.01e-5 \\
POD-DON & 1.31e-1 $\pm$ 7.16e-2 & 1.63e-4 $\pm$ 6.25e-4 \\
SPOD-DON & - & 1.75e-5 $\pm$ 6.75e-5 \\
TC-DON & 1.81e-1 $\pm$ 9.77e-2 & \underline{1.16e-5 $\pm$ 3.31e-5}  \\
TC-SPOD-DON & - & 4.00e-5 $\pm$ 1.18e-4 \\
\hdashline
\textbf{TINO} & \underline{6.85e-2 $\pm$ 6.18e-2}  & \textbf{5.09e-6 $\pm$  1.28e-5}  \\
\bottomrule
\end{tabular}
\end{table}

We do not evaluate the SPOD variants or U-Net in the purely temporal case, since the SPOD variants are designed for spatial problems and UNet is not suitable for one-dimensional data. In the spatio-temporal case, FNO and UNet are not evaluated because their input and output grids are defined on different domains.

\begin{figure}[h]
    \centering
    \begin{subfigure}{0.3\textwidth}
        \centering
        \includegraphics[width=\linewidth]{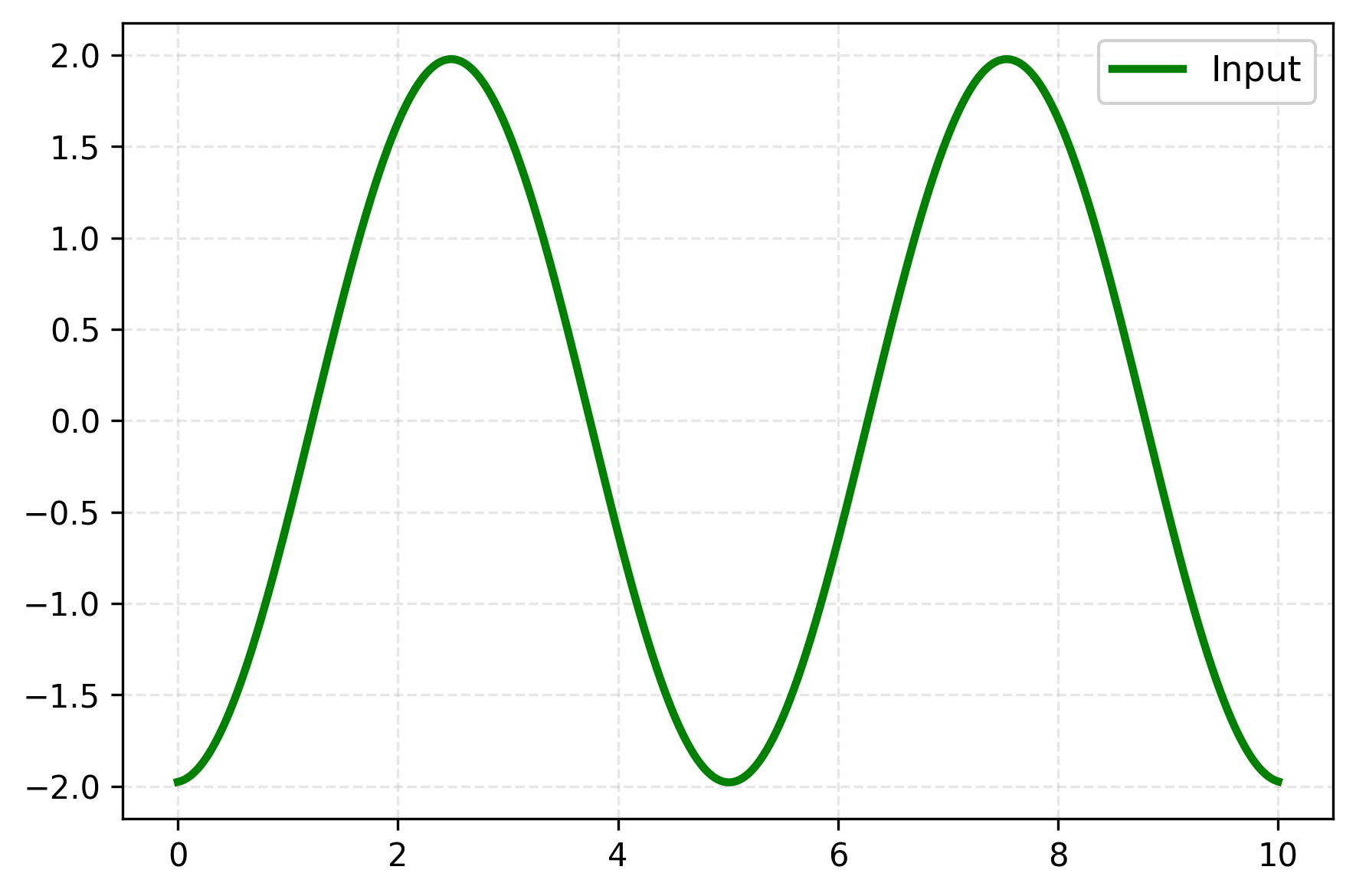}
        \caption{Input function}
        \label{fig:fitting_time_input_func}
    \end{subfigure}
    \hfill
    \begin{subfigure}{0.3\textwidth}
        \centering
        \includegraphics[width=\linewidth]{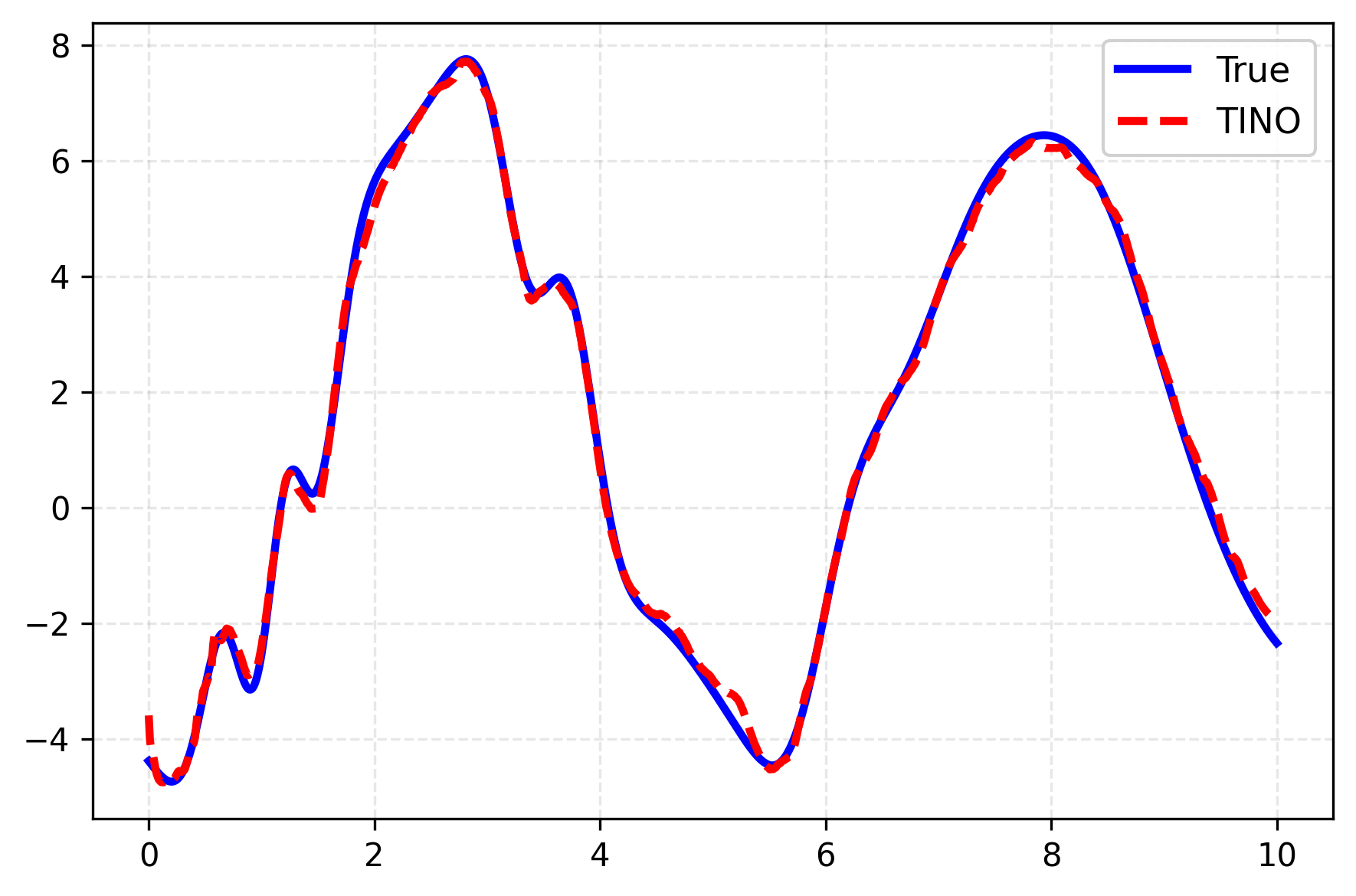}
        \caption{Ground truth vs TINO prediction}
        \label{fig:fitting_time_comparison}
    \end{subfigure}
    \hfill
    \begin{subfigure}{0.3\textwidth}
        \centering
        \includegraphics[width=\linewidth]{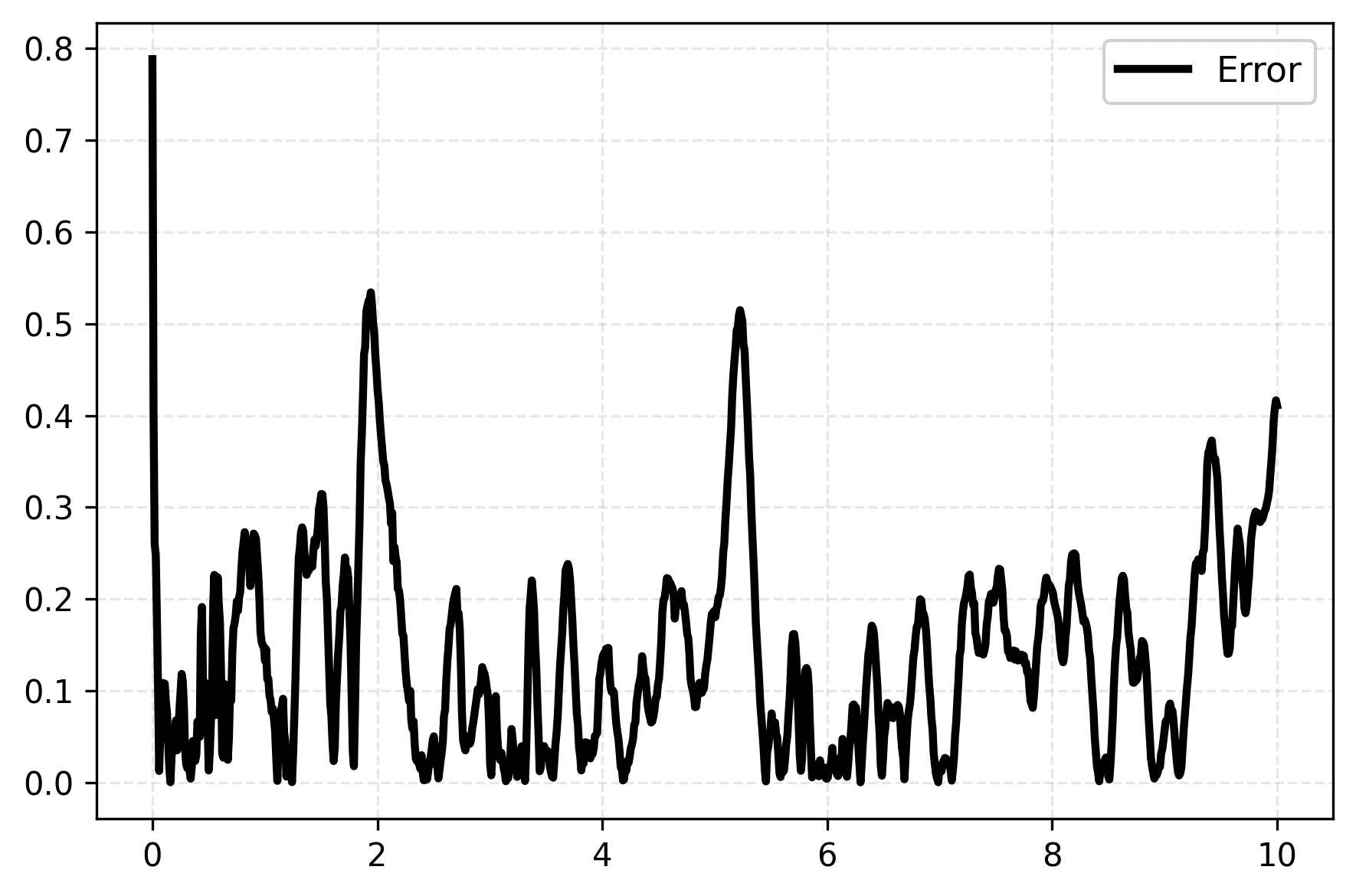}
        \caption{Absolute error}
        \label{fig:fitting_time_error}
    \end{subfigure}
    \caption{Example output of TINO for fitting time-only operator (\Cref{subsubsect_fitting_time}).}
    \label{fig:fitting_time_three_subfigures}
\end{figure}

\begin{figure}[h]
	\centering
	\begin{subfigure}{0.2\textwidth}
		\centering
		\includegraphics[width=\linewidth]{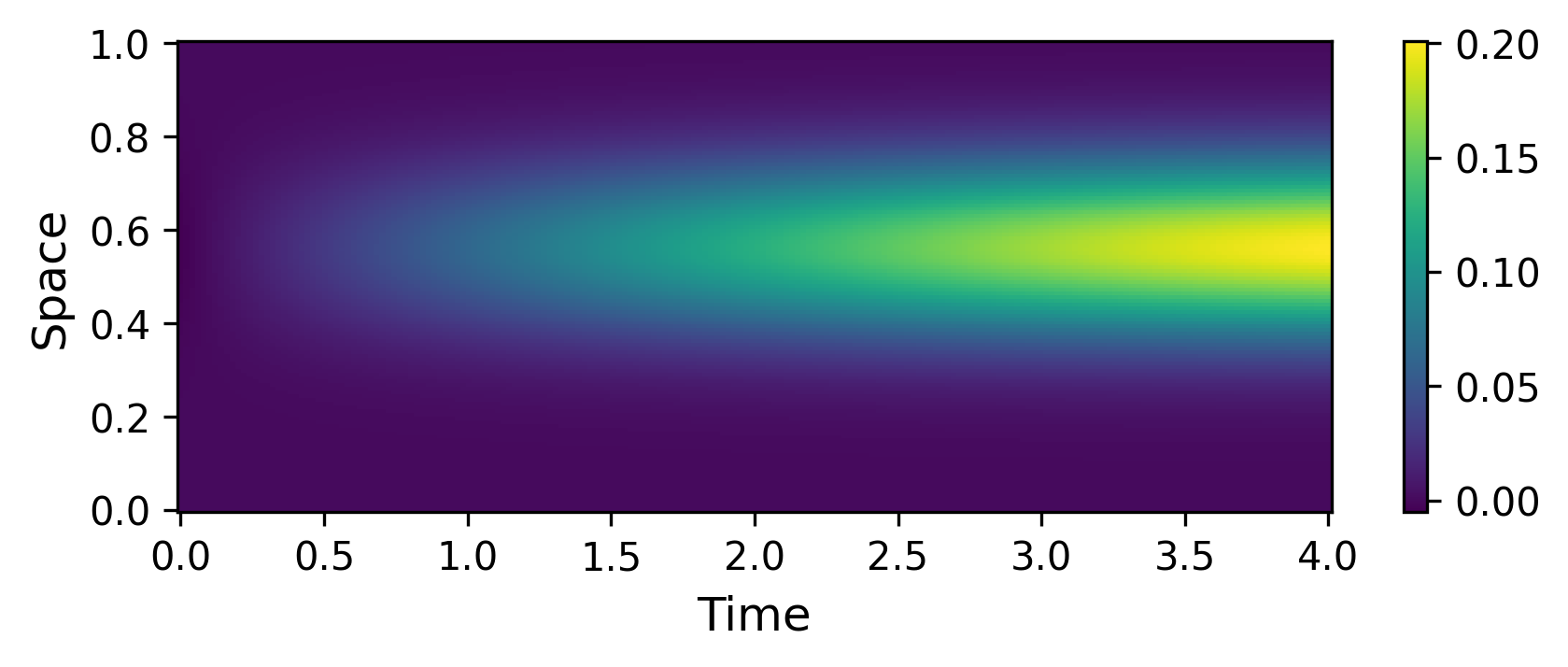}
		\caption{Input function}
		\label{fig:fitting_space_and_time_input_func}
	\end{subfigure}
	\hfill
	\begin{subfigure}{0.2\textwidth}
		\centering
		\includegraphics[width=\linewidth]{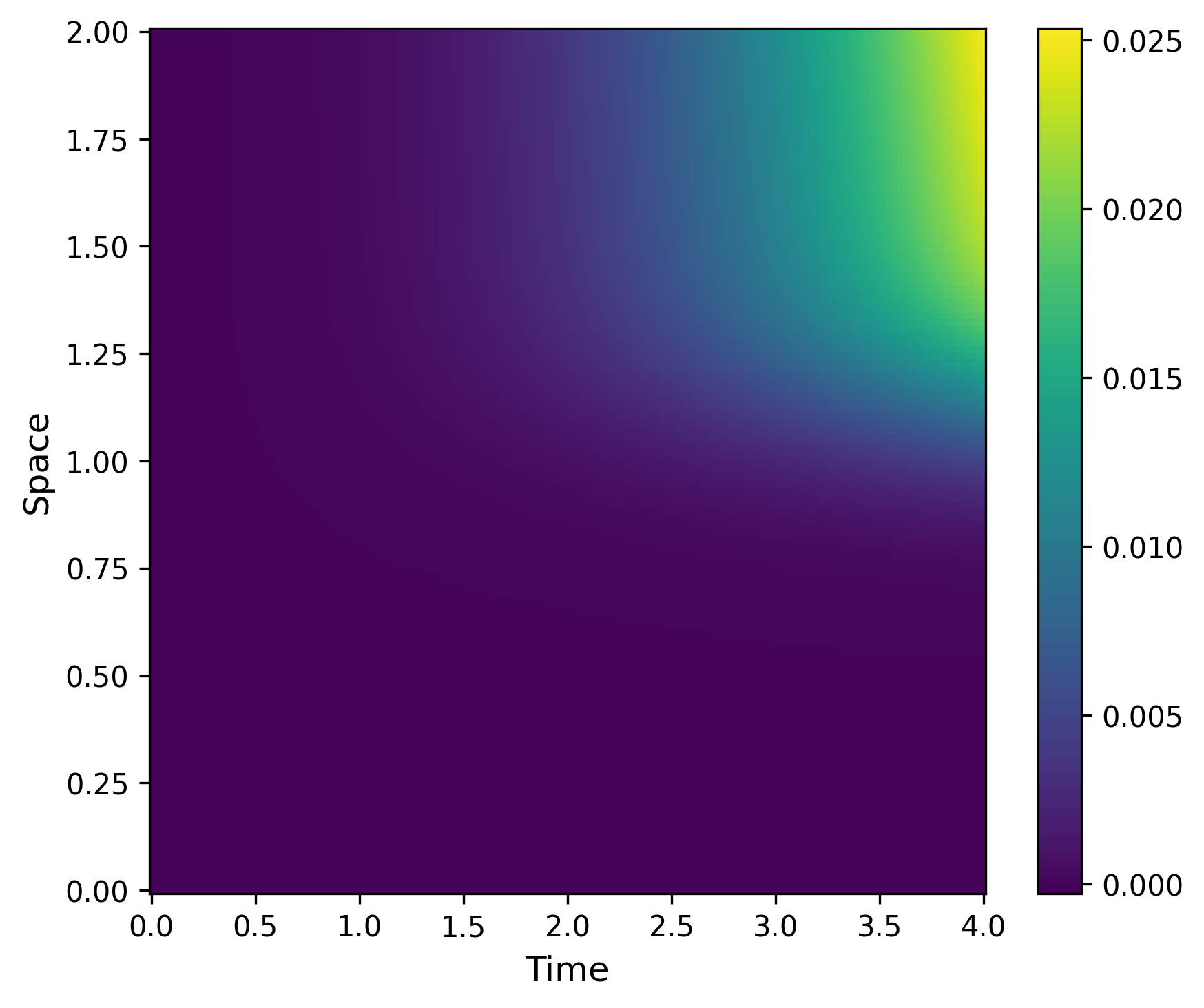}
		\caption{Ground truth}
		\label{fig:fitting_space_and_time_input_func}
	\end{subfigure}
	\hfill
	\begin{subfigure}{0.2\textwidth}
		\centering
		\includegraphics[width=\linewidth]{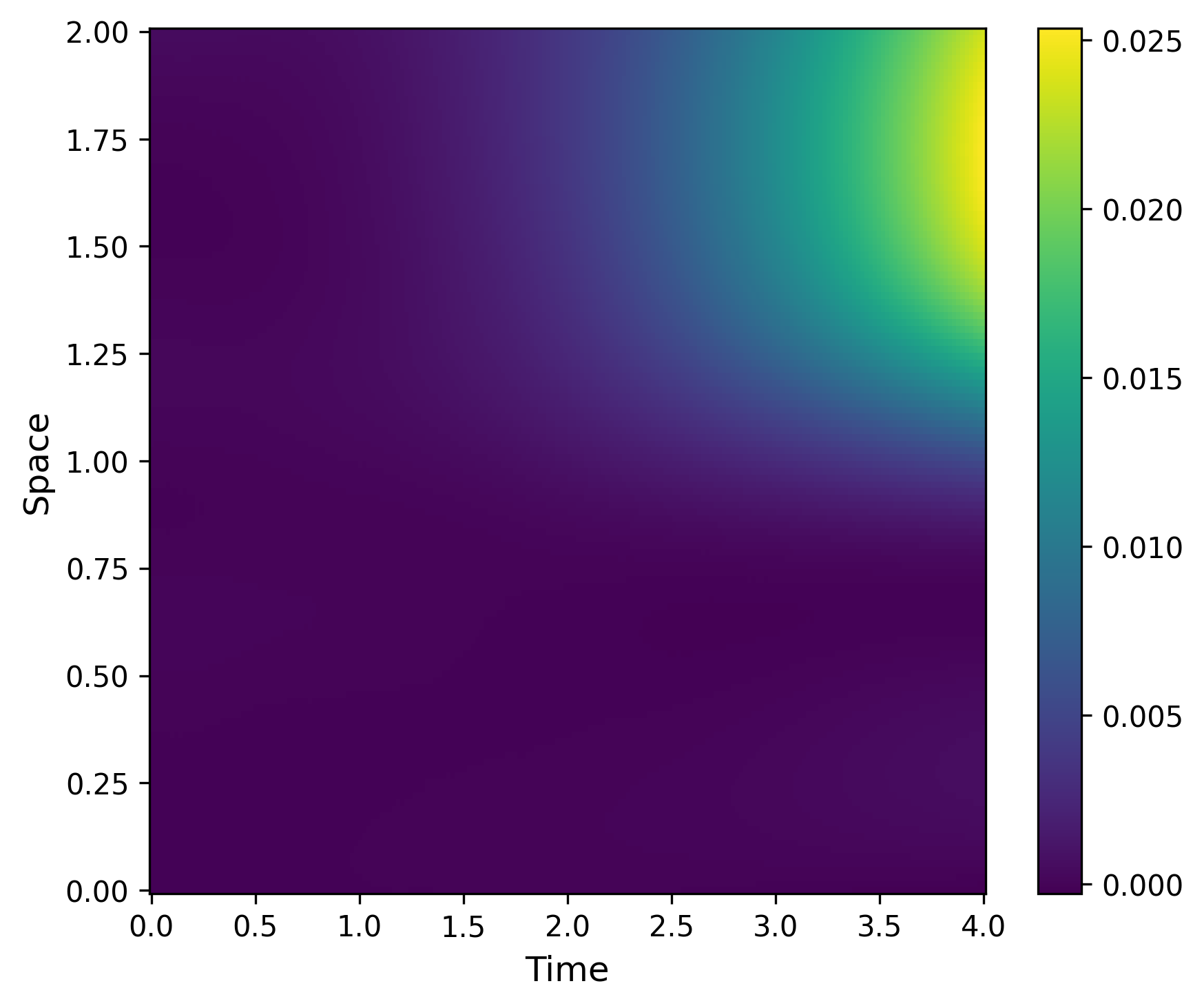}
		\caption{TINO Prediction}
		\label{fig:fitting_space_and_time_input_func}
	\end{subfigure}
	\hfill
	\begin{subfigure}{0.2\textwidth}
		\centering
		\includegraphics[width=\linewidth]{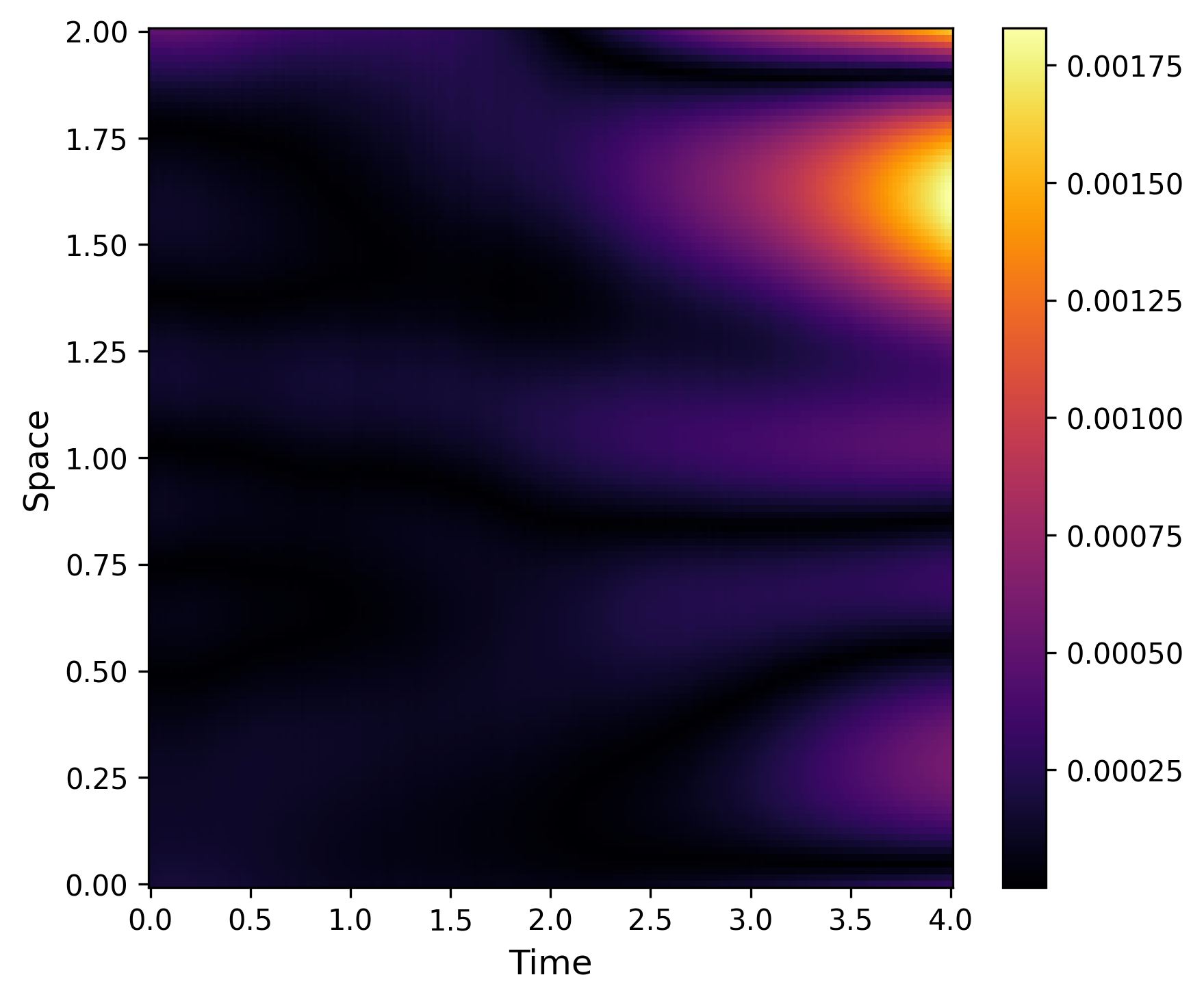}
		\caption{Absolute error}
		\label{fig:fitting_space_and_time_input_func}
	\end{subfigure}
	\caption{Example output of TINO for fitting time-dependent operator in 1-D space (\Cref{subsubsect_fitting_time_space}).}
	\label{fig:fitting_space_and_time_four_subfigures}
\end{figure}

We conclude the test results in \Cref{tab:mse_results_fitting}.
TINO achieves the best accuracy on the spatio-temporal fitting task, outperforming all DeepONet variants by approximately one order of magnitude. For the purely temporal case, FNO performs best while TINO remains competitive. 
Example outputs of TINO are illustrated in \Cref{fig:fitting_time_three_subfigures,fig:fitting_space_and_time_four_subfigures}.

\subsection{Solving time-dependent ODEs and ODE systems}\label{sect_sol_ode}

{In this section, we focus on learning the solution operators for ODEs/ODE systems from external source terms \(f(t)\) to solution \(u(t)\) with fixed initial conditions. The generation of training and testing data follows the same procedure discussed in \cref{subsubsect_fitting_time}.

\subsubsection{Duffing oscillator}\label{subsubsect_duffing_osc}
We consider the following second-order ODE governing the Duffing oscillator:
\begin{align}
\begin{split}
    \frac{\mathrm{d}^2u}{\mathrm{d}t^2} + c\frac{\mathrm{d}u}{\mathrm{d}t} + u + u^3 &= f(t), \quad t \in [0, T],\\
    u(0) = u'(0) &= 0.
\end{split}
\end{align}
Here, $f(t)$ represents the external force over time, while $u$, $\frac{\mathrm{d}u}{\mathrm{d}t}$, and $\frac{\mathrm{d}^2u}{\mathrm{d}t^2}$ denote displacement, velocity, and acceleration, respectively.
We set the damping coefficient $c=1$, and learn the operator $f(t) \in C[0,T] \mapsto u(t) \in C^2[0,T]$.

\subsubsection{Gradient flow on a double-well potential} \label{subsubsect_gradient_flow}
We consider the following forced gradient-flow equation
\begin{align}
\begin{split}
    \frac{\mathrm{d} u}{\mathrm{d} t} + F'(u) &= f(t), \quad t \in [0, T], \\
    u(0) &= 0,
\end{split}
\end{align}
where $f(t)$ is the external force, $u(t)$ is the phase field variable, and $F(u)=(u^2-1)^2/4$ is a double-well potential. We learn the operator $f(t) \in C[0,T] \mapsto u(t) \in C[0,T]$.

\subsubsection{Lorenz system} \label{subsubsect_lorenz_sys}
The Lorenz system considered has the following form:
\begin{equation}
    \begin{split}
        \frac{\mathrm{d}u}{\mathrm{d}t} &= 10(y-u), \\
        \frac{\mathrm{d}y}{\mathrm{d}t} &= u(\rho-z) - y, \\
        \frac{\mathrm{d}z}{\mathrm{d}t} &= uy - \frac{8}{3}z - f(t), \\
        u(0) &= 1, \quad y(0) = z(0) = 0.
    \end{split}
\end{equation}
Here, $t \in [0, T]$, $f(t)$ is the external force, $u(t)$ is proportional to the rate of convection, and $y(t), z(t)$ represent horizontal and vertical temperature variations, respectively \cite{cao2024laplace, shi2021analysis}.
We set $\rho = 10$ and learn the operator $f(t) \in C[0,T] \mapsto u(t) \in C[0,T]$.

\subsubsection{Results}

\begin{table}[h]
\centering
\caption{Error of neural operators on learning solution operators of ODEs. Results are reported as the mean squared error (MSE) $\pm$ standard deviation (SD).} 
\label{tab:mse_results_ode}
\begin{tabular}{lccc}
\toprule
&\Cref{subsubsect_duffing_osc} &\Cref{subsubsect_gradient_flow} &\Cref{subsubsect_lorenz_sys}\\
%\midrule
 & Duffing oscillator & gradient flow & Lorenz system \\
\midrule
FNO & \textbf{8.18e-6 $\pm$ 1.41e-5}  & \textbf{2.14e-2 $\pm$ 1.10e-1}  & \underline{1.71e-6 $\pm$ 2.16e-6}  \\
\hdashline
DON & 1.48e-3 $\pm$ 2.43e-3 & 7.25e-2 $\pm$ 1.30e-1 & 2.61e-3 $\pm$ 3.09e-3 \\
POD-DON & 1.97e-3 $\pm$ 3.42e-3 & 1.18e-1 $\pm$ 1.58e-1 & \textbf{1.02e-6 $\pm$ 2.14e-6}  \\
TC-DON & 1.48e-3 $\pm$ 2.54e-3 & \underline{3.56e-2 $\pm$ 9.89e-2}  & 1.71e-3 $\pm$ 1.04e-3 \\
\hdashline
\textbf{TrTINO} & \underline{3.93e-5 $\pm$ 6.32e-5}  & 3.67e-2 $\pm$ 1.05e-1 & 6.32e-4 $\pm$ 7.73e-4 \\
\bottomrule
\end{tabular}
\end{table}

Again, we skip the evaluation for SPOD variants and UNet for similar reasons discussed in \cref{sect_fitting_ti_op}.

We summarize the test results for ODEs and ODE systems in \Cref{tab:mse_results_ode}.
For the Duffing oscillator problem, TrTINO achieves the second-best performance, substantially outperforming standard DeepONet variants.
FNO excels on this smooth temporal problem.
On the gradient flow task, TC-DeepONet and TrTINO show comparable performance, both significantly better than vanilla DeepONet.
For the Lorenz system, POD-DeepONet achieves the best results, due to the low-rank structure being sufficiently well captured by POD.
TrTINO maintains reasonable accuracy for the ODE system.

Visualizations of representative cases of TrTINO are illustrated in \cref{fig:ode_duffing_oscillato_three_subfigures,fig:ode_double_well_potential_three_subfigures,fig:ode_lorenz_system_three_subfigures}.

\begin{figure}[h]
    \centering
    \begin{subfigure}{0.3\textwidth}
        \centering
        \includegraphics[width=\linewidth]{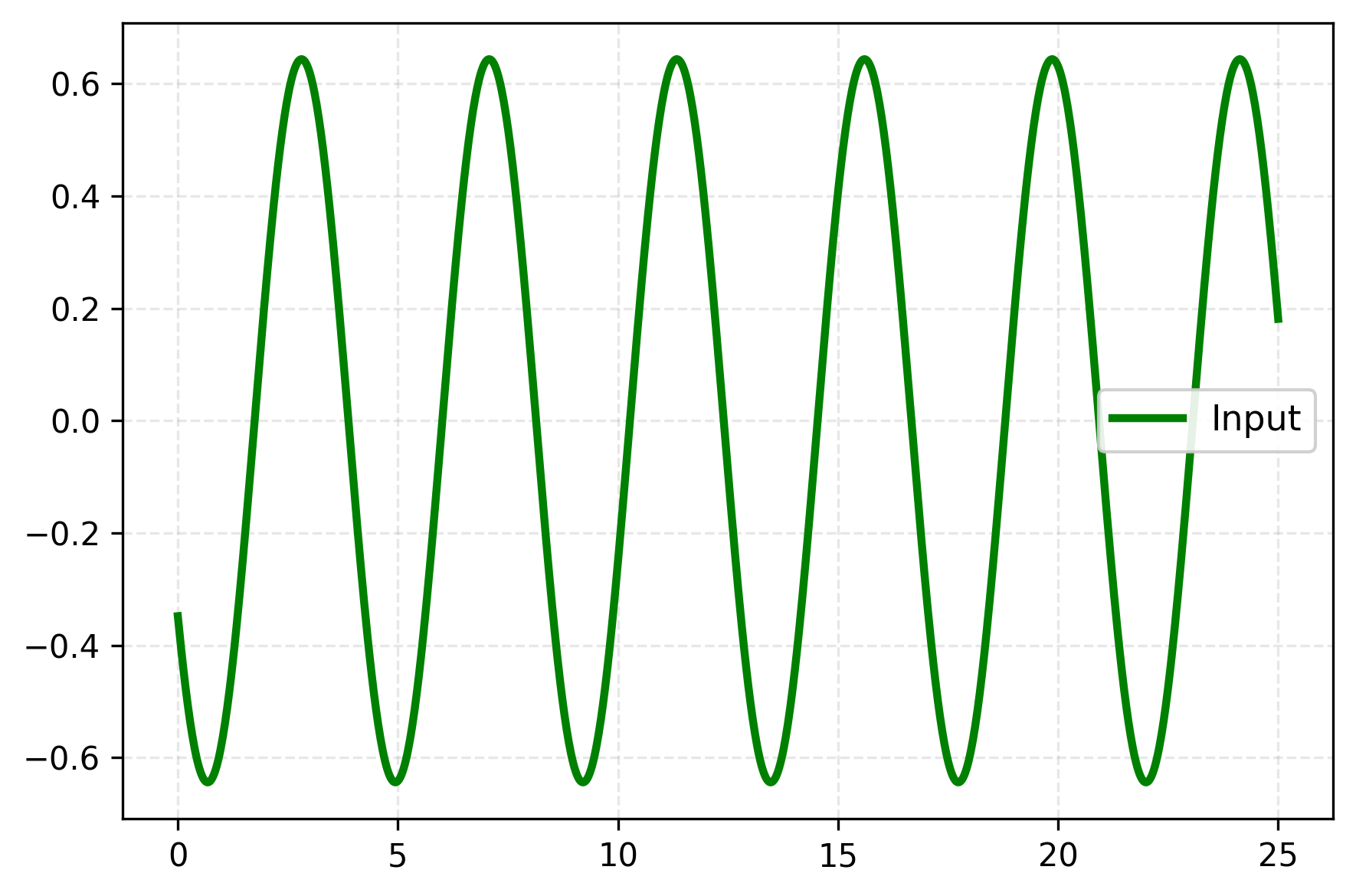}
        \caption{Input function}
        \label{fig:ode_duffing_oscillato_input_func}
    \end{subfigure}
    \hfill
    \begin{subfigure}{0.3\textwidth}
        \centering
        \includegraphics[width=\linewidth]{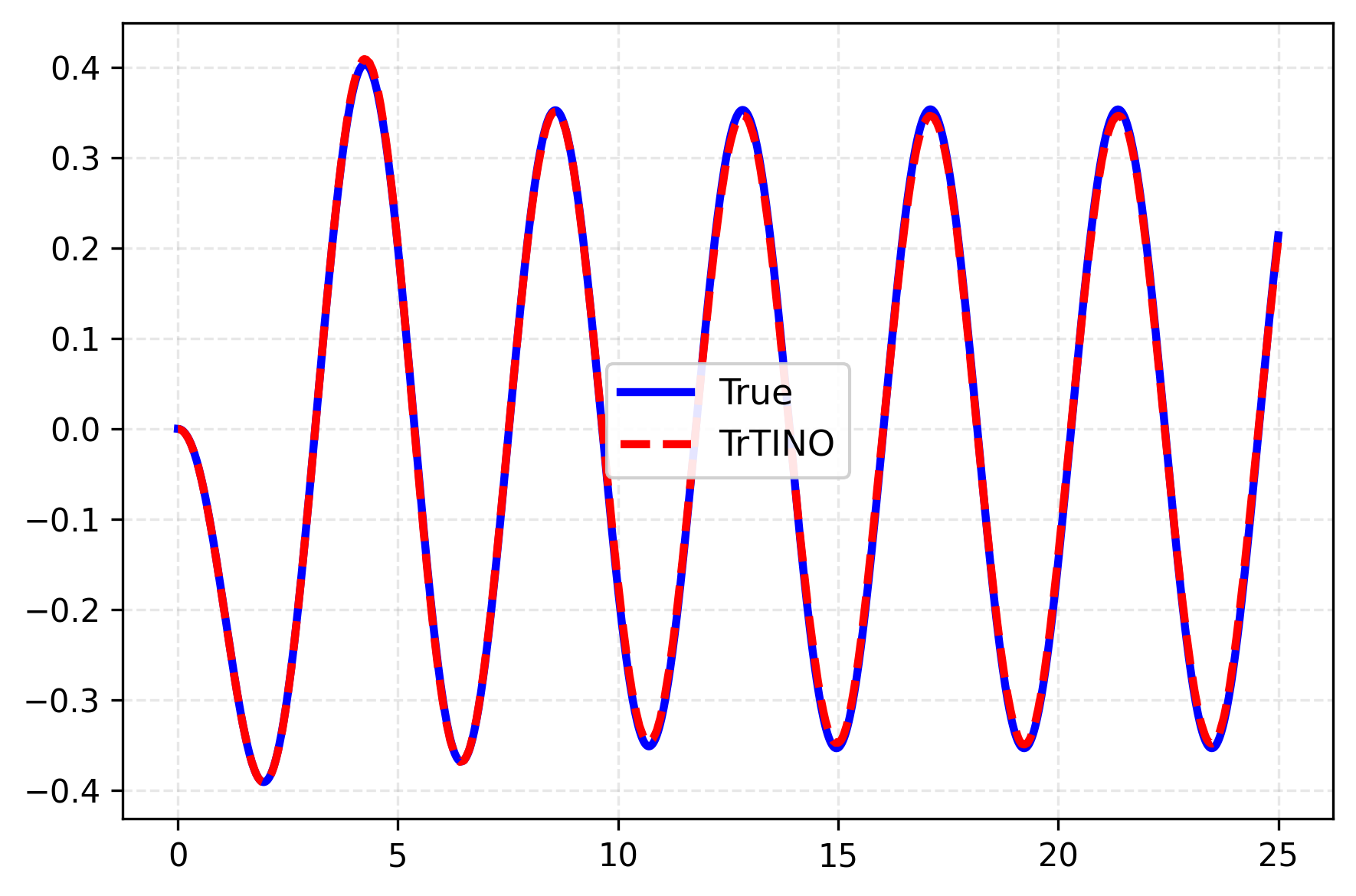}
        \caption{Ground truth vs. TrTINO prediction}
        \label{fig:ode_duffing_oscillato_comparison}
    \end{subfigure}
    \hfill
    \begin{subfigure}{0.3\textwidth}
        \centering
        \includegraphics[width=\linewidth]{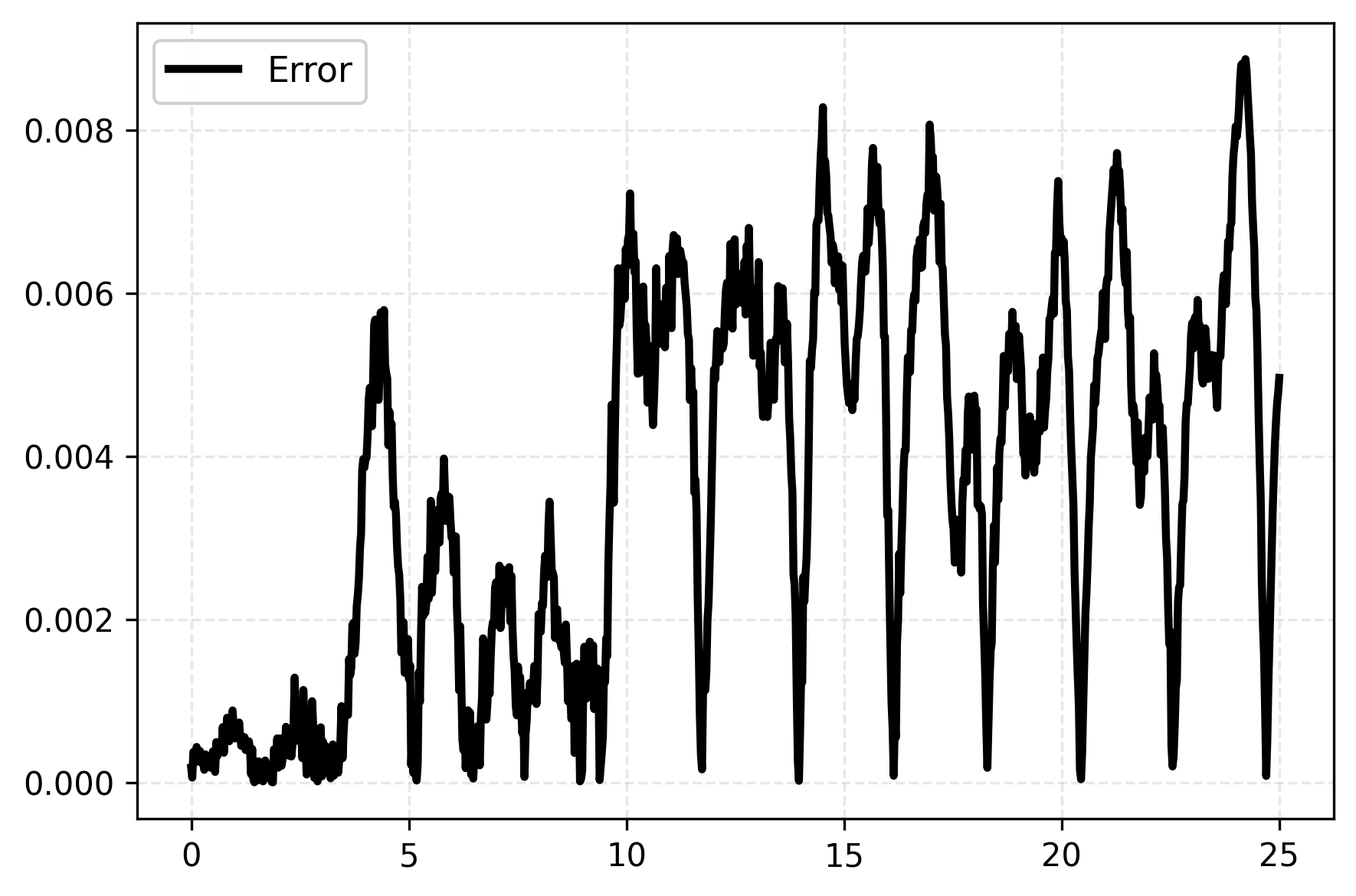}
        \caption{Absolute error}
        \label{fig:ode_duffing_oscillato_error}
    \end{subfigure}
    \caption{Example output of TrTINO for learning the solution of Duffing oscillator (\Cref{subsubsect_duffing_osc}).}
    \label{fig:ode_duffing_oscillato_three_subfigures}
\end{figure}

\begin{figure}[h]
    \centering
    \begin{subfigure}{0.3\textwidth}
        \centering
        \includegraphics[width=\linewidth]{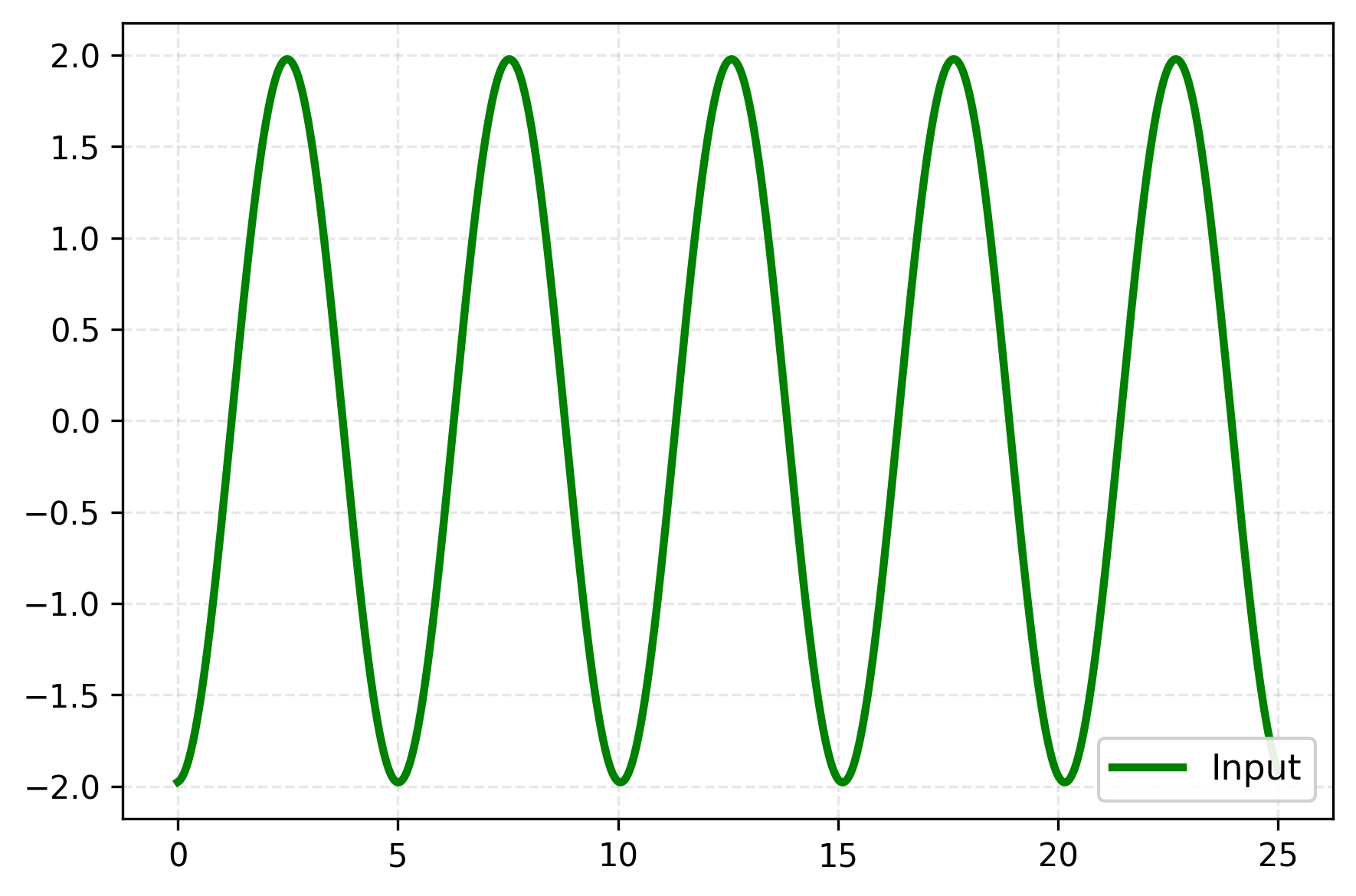}
        \caption{Input function}
        \label{fig:ode_double_well_potential_input_func}
    \end{subfigure}
    \hfill
    \begin{subfigure}{0.3\textwidth}
        \centering
        \includegraphics[width=\linewidth]{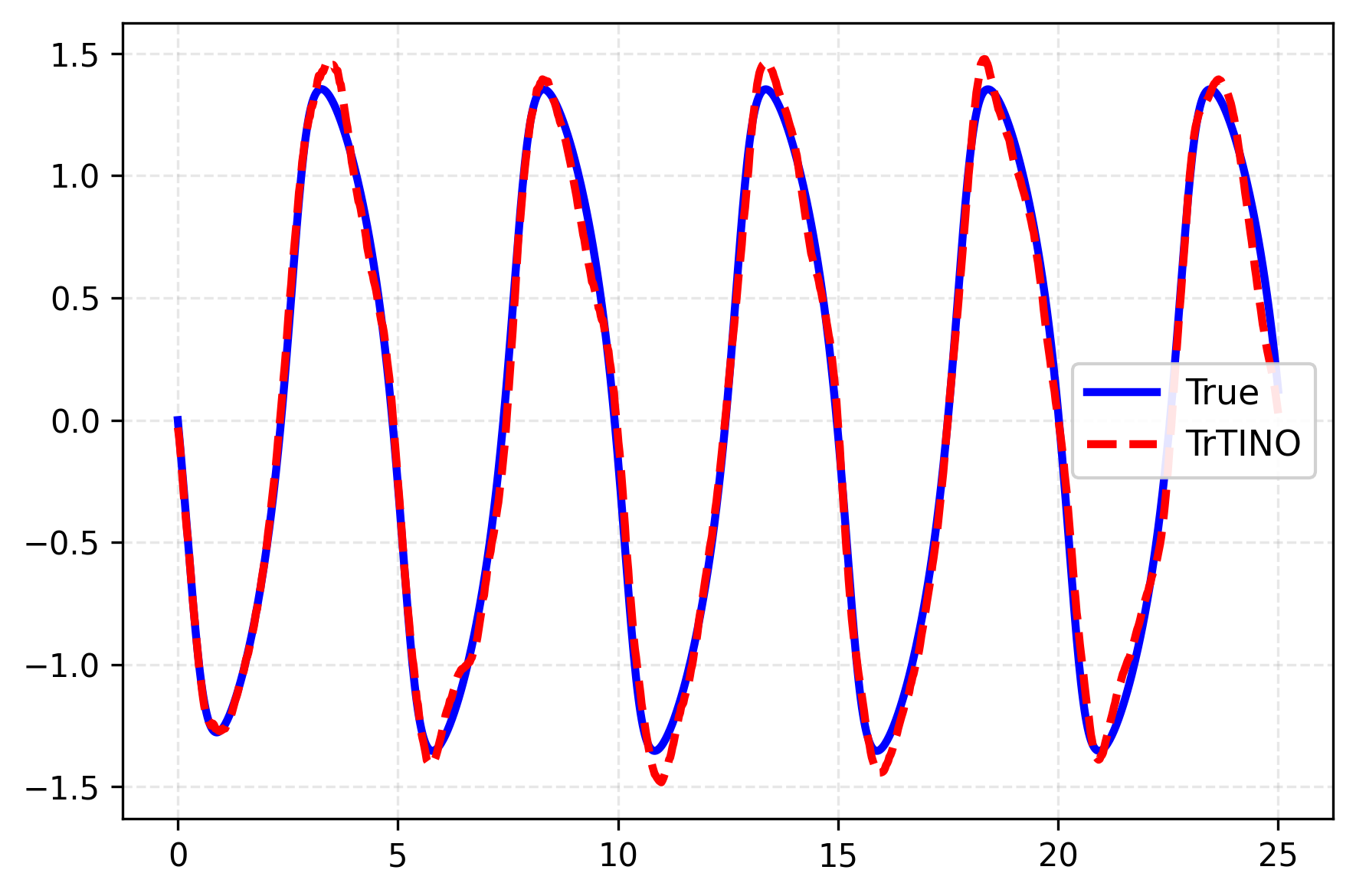}
        \caption{Ground truth vs. TrTINO prediction}
        \label{fig:ode_double_well_potential_comparison}
    \end{subfigure}
    \hfill
    \begin{subfigure}{0.3\textwidth}
        \centering
        \includegraphics[width=\linewidth]{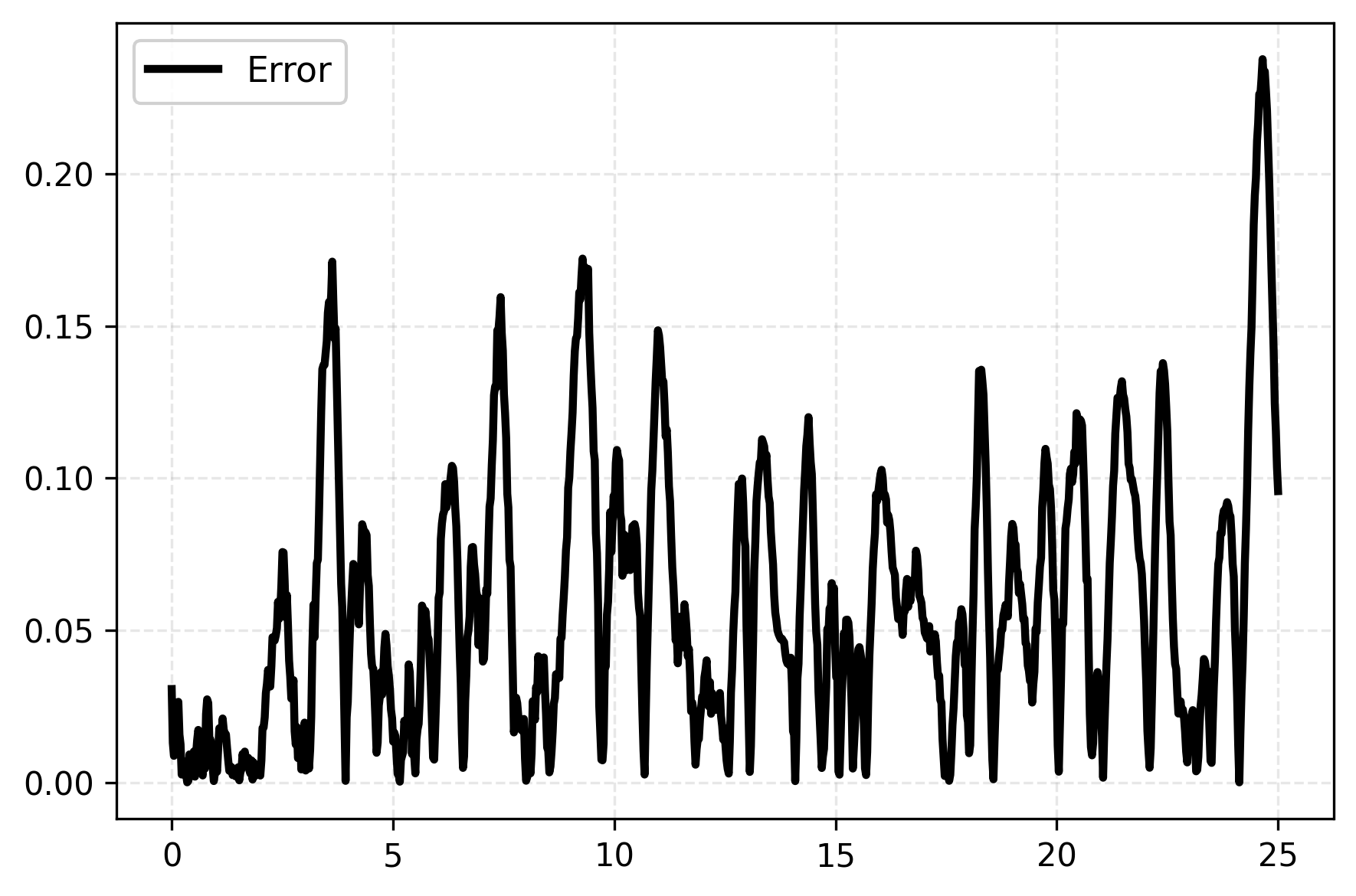}
        \caption{Absolute error}
        \label{fig:ode_double_well_potential_error}
    \end{subfigure}
    \caption{Example output of TrTINO for learning the solution of gradient flow on a double-well potential (\Cref{subsubsect_gradient_flow}).}
    \label{fig:ode_double_well_potential_three_subfigures}
\end{figure}

\begin{figure}[h]
    \centering
    \begin{subfigure}{0.3\textwidth}
        \centering
        \includegraphics[width=\linewidth]{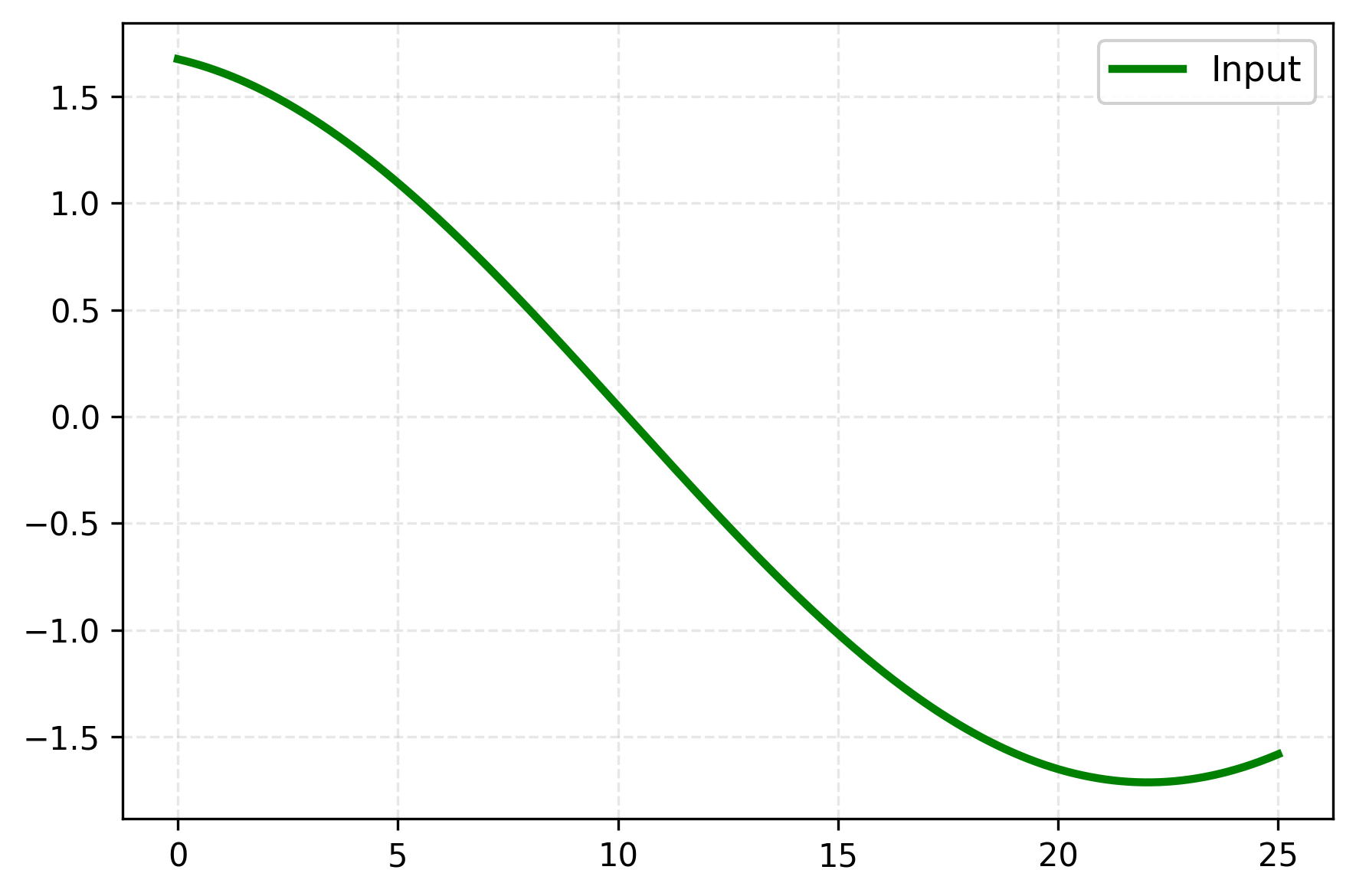}
        \caption{Input function}
        \label{fig:ode_lorenz_system_input_func}
    \end{subfigure}
    \hfill
    \begin{subfigure}{0.3\textwidth}
        \centering
        \includegraphics[width=\linewidth]{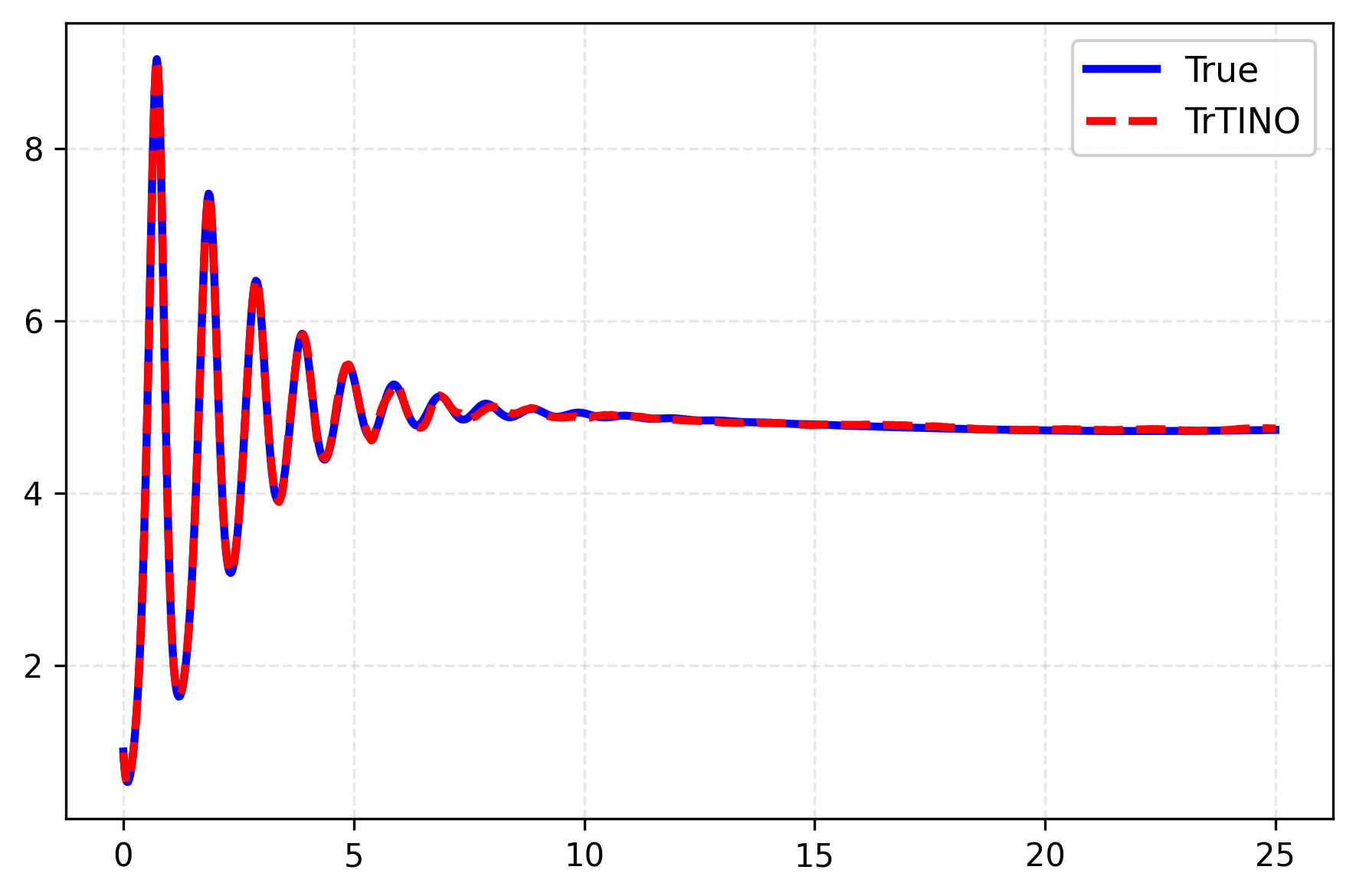}
        \caption{Ground truth vs. TrTINO prediction}
        \label{fig:ode_lorenz_system_comparison}
    \end{subfigure}
    \hfill
    \begin{subfigure}{0.3\textwidth}
        \centering
        \includegraphics[width=\linewidth]{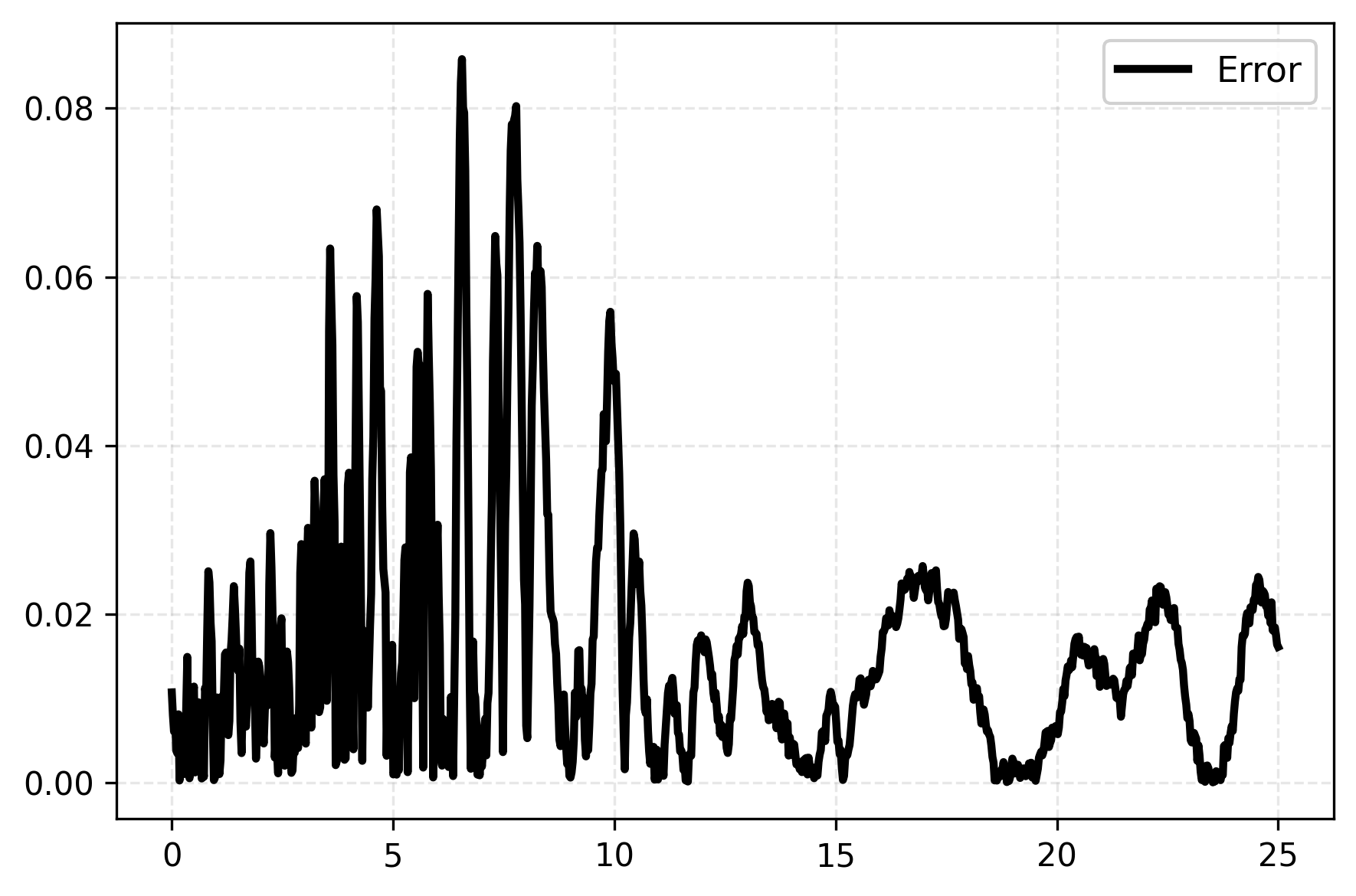}
        \caption{Absolute error}
        \label{fig:ode_lorenz_system_error}
    \end{subfigure}
    \caption{Example output of TrTINO for learning the solution operator of Lorenz system (\Cref{subsubsect_lorenz_sys}).}
    \label{fig:ode_lorenz_system_three_subfigures}
\end{figure}

\subsection{Solving time-dependent PDEs}\label{sect_sol_pde}
We extend our evaluation to time-dependent PDEs with initial conditions. We benchmark SPOD-TrTINO and comparative models by learning the solution operators in the form $f(t,\vec{x}) \in C([0,T]\times \Omega) \mapsto u(t,\vec{x}) \in C([0,T]\times \Omega)$.
For PDEs with 1-D in space, we use $n_{\mathrm{t}} = 200$ temporal samples and $n_1=128$ spatial samples for picking sensors.
For PDEs with 2-D in space, we use $n_{\mathrm{t}} = 100$ temporal samples and $n_1 \times n_2= 32 \times 32$ spatial samples.
For PDEs with 3-D in space, we use $n_{\mathrm{t}} = 50$ temporal samples and $n_1 \times n_2 \times n_3=16 \times 16 \times 16$ spatial samples due to limited memory.

\subsubsection{1-D Allen--Cahn equation} \label{subsubsect_allen_cahn}
We consider the 1-D Allen--Cahn equation for $u \in C^1([0, T] \times [0, L])$:
\begin{equation}
    \begin{split}
        \frac{\partial u}{\partial t} - \varepsilon^2 \left(\frac{\partial u}{\partial x}\right)^2 &= - F'(u) + f(t,x), \quad (t, x) \in [0, T] \times [0, L], \\
        u(0,x) &= u_0(x), \quad x \in [0, L], \\
        \frac{\partial u(t, x)}{\partial t}&=0, \quad t \in [0, T], \quad x \in \{0, L\}.
    \end{split}
\end{equation}
Here, $f(t,x)$ is the external force, $u(t,x)$ is the phase field variable, $\varepsilon$ is the interfacial width parameter,
$F(u)=(u^2-1)^2/4$ is the double-well potential.
We learn the solution operator $f(t,x) \mapsto u(t,x)$.

In this experiment, we set $T=4$, $L=1$, $\varepsilon=0.01$.
Input functions are restricted by the parametrization
\begin{equation*}
    f(t,x ) = A\sin(bt+c) e^{-(x-\mu)^2/(2\sigma^2)}
\end{equation*}
with $A \in [0.1, 1]$, $b \in [0.2, 1]$, $c \in [0, 2\pi]$, $\mu \in [L/4, 3L/4]$, $\sigma \in [0.1L, 0.4L]$.
For the initial condition, we set $u_0(x) = 0.5$ for any $x \in [0, L]$.

\subsubsection{1-D viscous Burgers' equation} \label{subsubsect_burgers}
We consider the one-dimensional viscous Burgers' equation:
\begin{equation}
    \begin{split}
        \frac{\partial u}{\partial t} + \frac{1}{2} \frac{\partial u^2}{\partial x} &= \nu \frac{\partial^2 u}{\partial x^2} + f(t, x), \quad (t, x) \in [0, T] \times [0, L] \\
        u(0,x) &= u_0(x), \quad x \in [0, L]
    \end{split}
\end{equation}
subject to periodic boundary conditions.
Here, $f(t,x)$ is the external force, $u(t,x)$ is the fluid velocity, and $\nu$ is the viscosity coefficient.
We learn the solution operator $f(t,x) \mapsto u(t,x)$ for the cases of $\nu=0.1$ and $\nu=0.01$, respectively.

In this experiment, we set $T=4$ and $L=1$. 
The input functions are parametrized by Fourier series
\begin{equation} \label{equ_burgers_input}
	f(t,x) = A\sin(bt+c) \sum_{n=1}^{10} \frac{a_n \sin(2\pi n x) + b_n \cos(2\pi n x)}{n^2},
\end{equation}
where $A \in [0.1, 1]$, $b \in [0.2, 1]$, $c \in [0, 2\pi]$, $a_n, b_n \in [-1, 1]$.
The initial condition considered has the form
\begin{equation} \label{equ_burgers_ic}
    u_0(x)=\sum_{n=1}^{10} \frac{\sin{(2n\pi x)} + \cos{(2n\pi x)}}{n^2}.
\end{equation}

\subsubsection{1-D convection-diffusion equation} \label{subsubsect_conv_diff}
We consider the 1-D convection-diffusion equation
\begin{equation}
    \begin{split}
        \frac{\partial u}{\partial t} + v\frac{\partial u}{\partial x} &= D \frac{\partial^2 u}{\partial x^2} + f(t,x), \quad (t, x) \in [0, T] \times [0, L], \\
        u(0,x) &= u_0(x), \quad x \in [0, L], \\
        \frac{\partial u}{\partial t}&=0, \quad t \in [0, T], \quad x \in \{0, L\}.
    \end{split}
\end{equation} 
Here, $f(t,x)$ is the external force, $u(t,x)$ represents the concentration, $v$ is the advection velocity, and $D$ is the diffusion coefficient.
We learn the solution operator $f(t,x) \mapsto u(t,x)$.

In this experiment, we set $T=4$, $L=1$ with $D=0.01$, $v=0.1$.
Input functions are parametrized as $f(t,x) = (A\sin(bt+c))^2 \cdot e^{-(x-\mu)^2/(2\sigma^2)}$, similar to the 1-D Allen--Cahn equation (\cref{subsubsect_allen_cahn}) test case.
The initial condition is given by $u_0(x) = 3\exp\left(-(x-L/2)^2/(0.1L^2)\right)$.

\subsubsection{2-D Navier--Stokes Equation} \label{subsubsect_navier_stokes}
We consider the 2-D Navier--Stokes equation for a viscous, incompressible fluid
\begin{equation}
    \begin{split}
    \frac{\partial w}{\partial t} + \vec{u} \cdot \nabla w &= \nu \Delta w + f(t,\vec{x}) \\
    \nabla \cdot \vec{u} &= 0 \\
    w(0,\vec{x}) &= w_0(\vec{x})        
    \end{split}
\end{equation}
for $(t, \vec{x}) \in [0, T] \times [0, L]^2$ with periodic boundary conditions.
Here, $f(t,\vec{x})$ is the vorticity forcing, $\vec{u}(t,\vec{x})$ is the velocity field, $w(t,\vec{x}) = \nabla \times \vec{u}$ is the vorticity, $w_0$ is the initial vorticity, $\nu$ is the constant viscosity coefficient.
We learn the solution operator $f(t,\vec{x}) \mapsto w(t,\vec{x})$ in the domain.

In this experiment, we set $T=2$, $L=1$ and $\nu=0.01$.
Input functions are generated using a 2-D Fourier series
\begin{equation} \label{equ_ns_input}
	f(t,\vec{x}) = A\sin(bt+c) \sum_{n=1}^{10} \frac{1}{n^2} \left( a_{n,1} \sin(2\pi n x_1) + b_{n,1} \cos(2\pi n x_1) + a_{n,2} \sin(2\pi n x_2) + b_{n,2} \cos(2\pi n x_2) \right)
\end{equation}
where $\vec{x} = (x_1, x_2)$, $A \in [0.1, 1]$, $b \in [0.2, 1]$, $c \in [0, 2\pi]$, and $a_{n,i}, b_{n,i} \in [-1, 1]$, $n=1,\dots,10$, $i=1,2$.
The initial vorticity
\begin{equation} \label{equ_ns_ic}
    w(0,\vec{x}) = \sum_{n=1}^{10} \frac{\sin(2\pi n x_1) + \cos(2\pi n x_1) + \sin(2\pi n x_2) + \cos(2\pi n x_2)}{n^2}.
\end{equation}

\subsubsection{2-D Navier--Stokes equation with variational initial condition} \label{subsubsect_navier_stokes_icvar}

We also include a test case based on \Cref{subsubsect_navier_stokes} but with variable initial vorticity $w_0(\vec{x})$.
The goal of the test is to learn the solution operator that maps both external force terms and  initial conditions to solutions
\begin{equation*}
    \left( f, w_0 \right) \mapsto w.
\end{equation*}
The test setup reuses the training data from the previous function pairs $(f(t, \vec{x}), w(t, \vec{x})$ for $t \in [0, T]$, but clipped into shorter time intervals
\begin{equation*}
    \left[\frac{k-1}{10}T, \frac{k+1}{10}T\right], \quad k=1,\cdots,9,
\end{equation*}
of size $0.2T$.
The test results are collected on all these $9$ intervals.

\subsubsection{Quasi 3-D advection-diffusion system}\label{subsubsect_3d_adc_diff}
We consider a system governing the evolution of a scalar field $w$ (generalized vorticity) under a 2-D velocity field extended to a 3-D domain:
\begin{equation}
    \begin{split}
    \frac{\partial w}{\partial t} + \vec{u} \cdot \nabla_{12} w &= \nu \Delta_{123} w + f(t, \vec{x}), \\
    \nabla_{12} \cdot \vec{u} &= 0, \\
    w &= \frac{\partial u_2}{\partial x_1} - \frac{\partial u_1}{\partial x_2}, \\
    w(0, \vec{x}) &= w_0(\vec{x})
    \end{split}
\end{equation}
for $t \in [0, T]$ and $\vec{x} = (x_1, x_2, x_3) \in [0, L]^3$ with periodic boundary conditions.
The subscripts in $\nabla$ and $\Delta$ refer to the partial derivatives on the corresponding entries of $\vec{x}$.
Here, $f(t,\vec{x})$ is the external force, the scalar field $w(t,\vec{x})$ diffuses in all three spatial dimensions, while the velocity $\vec{u}(t,\vec{x}) = (u_1, u_2, 0)$ remains planar.
This models a stack of interacting 2-D fluid layers coupled solely via viscosity in the perpendicular direction.
We learn the solution operator $f(t,\vec{x}) \mapsto w(t,\vec{x})$.

In this experiment, we set $T=1$, $L=1$ and $\nu=0.001$.
The input functions are generated as
\begin{equation} \label{equ_3d_input_scalar}
    \begin{split}
        f(t,\vec{x}) = A\sin(bt+c) \sum_{n=1}^{10} \sum_{i=1}^{3} \frac{1}{n^2} \left( a_{n,i} \sin(2\pi n x_i) + b_{n,i} \cos(2\pi n x_i) \right)
    \end{split}
\end{equation}
where $A \in [0.1, 1]$, $b \in [0.2, 1]$, $c \in [0, 2\pi]$, and the spatial coefficients $a_{n,i}, b_{n,i} \in [-1, 1]$, $n=1,\dots,10$, $i=1,2,3$.
The initial value function $w_0(\vec{x})$ is constructed as a multi-modal superposition of Fourier terms
\begin{equation} \label{equ_3d_ic_scalar}
    w_0(\vec{x}) = \sum_{n=1}^{10} \sum_{i=1}^{3} \frac{1}{n^2} \left( \sin(2\pi n x_i) + \cos(2\pi n x_i) \right).
\end{equation}

\subsubsection{Heat equation in a triangular domain} \label{subsubsect_heat_equ}
We consider 2-D heat equation with a source term
\begin{equation}
\begin{split}
    \frac{\partial u}{\partial t} - \alpha \Delta u &=  f(t,\vec{x}), \\
    u(0,\vec{x})&=0
\end{split}
\end{equation}
for $t, \vec{x} \in [0, T] \times \Omega$, with Neumann boundary conditions $\frac{\partial u}{\partial n}=0$ on $\partial \Omega$.
Here, $f(t,\vec{x})$ is the external force, $u(t,\vec{x})$ is the temperature field, and $\alpha$ is the thermal diffusivity.
We set $\Omega$ as an equilateral triangle with vertices
\begin{equation*}
    (0, 0), \quad (L, 0), \quad (L/2, \sqrt{3}L / 2).
\end{equation*}
We learn the operator $f(t,\vec{x})\mapsto u(t,\vec{x})$.

In this experiment, we set $T=1$, $L=1$ and $\alpha=0.005$.
Input functions are generated as
\begin{equation} \label{equ_heat_input}
    f(t,\vec{x}) = (A \sin(bt+c))^2 \exp\left( -\frac{(x_1-x_{1\mathrm{c}})^2 + (x_2-x_{2\mathrm{c}})^2}{2\sigma^2} \right)
\end{equation}
where $A \in [0.1, 1]$, $b \in [0.2, 1]$, $c \in [0, 2\pi]$, and the source parameters are sampled with $x_{1\mathrm{c}}, x_{2\mathrm{c}} \in [0.3, 0.7]$ and $\sigma \in [0.2, 0.3]$.

With the above setup, during the initial stage of evolution, the spatial profile of the temperature field is primarily dominated by the external heat source $f(t,\vec{x})$, as the diffusion effect $\alpha \Delta u$ is negligible when the temperature gradients are small.

For FNO and U-Net which do not naturally support solving on the triangular domain, 
we embed the domain into a rectangular region $[0, L]^2$, following the strategy in \cite{lu2022comprehensive}.

\subsubsection{1-D advection-reaction equation} \label{subsubsect_adv_react}
We consider the 1-D time-dependent advection-reaction equation
\begin{equation}
    \begin{split}
        \frac{\partial u}{\partial t} + v(t)\frac{\partial u}{\partial x} + \sin (0.1u) &= f(t, x), \quad (t, x) \in [0, T] \times [0, L], \\
        u(0,x) &= u_0(x), \quad x \in [0, L],
    \end{split}
\end{equation}
with periodic boundary conditions.
Here, $f(t,x)$ is the external source term, $u(t,x)$ is the state variable, and $R(u) = -\sin(0.1 u)$ is a nonlinear reaction term.
We learn the solution operator $f(t,x) \mapsto u(t,x)$.

We emphasize that the above system is time-causal but in general \emph{not} time-invariant, due to the time-dependent advection velocity $v(t)$.

In this experiment, we set $T=4$, $L=1$, $v(t)=\cos t$.
The input functions and initial value functions are chosen to be identical to those specified for the Burgers' equation in \Cref{subsubsect_burgers}, following \Cref{equ_burgers_input} and \Cref{equ_burgers_ic}, respectively.

\subsubsection{Results}

\begin{table}[h]
\centering
\caption{Error of neural operators on learning solution operators of time-dependent 1-D PDEs. Results are reported as the mean squared error (MSE) $\pm$ standard deviation (SD).}
\label{tab:mse_comparison_part1}
\resizebox{\columnwidth}{!}{
\begin{tabular}{lcccc}
\toprule
& \Cref{subsubsect_allen_cahn} & \Cref{subsubsect_burgers} & \Cref{subsubsect_burgers} & \Cref{subsubsect_conv_diff} \\
%\midrule
& 1-D Allen--Cahn & 1-D Burgers' ($\nu=0.1$) & 1-D Burgers' ($\nu=0.01$) & 1-D conv. diff. \\
\midrule
FNO & \underline{4.07e-7 $\pm$ 6.98e-7} & \underline{3.00e-7 $\pm$ 2.52e-7} & \textbf{2.03e-5 $\pm$ 3.24e-5}  & 6.07e-6 $\pm$ 7.82e-6 \\
U-Net & 4.19e-6 $\pm$ 3.67e-6 & 2.44e-5 $\pm$ 7.00e-5 & 3.05e-4 $\pm$ 4.58e-4 & 3.15e-4 $\pm$ 1.03e-3 \\
\hdashline
DON  & 2.59e-5 $\pm$ 5.68e-5 & 2.84e-4 $\pm$ 5.31e-5 & 3.29e-3 $\pm$ 2.48e-3 & 1.97e-4 $\pm$ 1.57e-4 \\
POD-DON  & 8.65e-7 $\pm$ 2.71e-6 & 1.52e-6 $\pm$ 2.05e-6 & 6.02e-4 $\pm$ 1.06e-3 & 5.59e-4 $\pm$ 8.97e-4 \\
SPOD-DON  & 7.80e-7 $\pm$ 2.33e-6 & 1.07e-6 $\pm$ 3.27e-6 & 3.45e-4 $\pm$ 7.04e-4 & 2.35e-5 $\pm$ 4.69e-5 \\
TC-DON  & 1.53e-5 $\pm$ 4.02e-5 & 7.85e-5 $\pm$ 7.61e-5 & 1.98e-3 $\pm$ 1.98e-3 & 5.57e-5 $\pm$ 5.33e-5 \\
TC-SPOD-DON & 4.32e-7 $\pm$ 1.33e-6 & 1.35e-6 $\pm$ 2.82e-6 & 1.49e-4 $\pm$ 3.10e-4 & \underline{2.84e-6 $\pm$ 8.33e-6}  \\
\hdashline
\textbf{SPOD-TrTINO} & \textbf{9.90e-8 $\pm$ 2.92e-7}  & \textbf{6.86e-8 $\pm$ 1.08e-7}  & \underline{1.04e-4$\pm$ 2.08e-4}  & \textbf{7.63e-7 $\pm$ 9.05e-7}  \\
\bottomrule
\end{tabular}
}
\end{table}

\begin{table}[h]
\centering
\caption{Error of neural operators on learning solution operators of time-dependent 2-D and 3-D PDEs. Results are reported as the mean squared error (MSE) $\pm$ standard deviation (SD).}
\label{tab:mse_comparison_part2}
\resizebox{\columnwidth}{!}{
\begin{tabular}{lcccc}
\toprule
& \Cref{subsubsect_navier_stokes} & \Cref{subsubsect_navier_stokes_icvar} & \Cref{subsubsect_3d_adc_diff} & \Cref{subsubsect_heat_equ}\\
%\midrule
& 2-D NS & 2-D NS (var. IC) & 3-D adv. diff. & 2-D heat \\
\midrule
FNO & \underline{1.16e-3 $\pm$ 1.96e-3} & \textbf{7.83e-06 $\pm$ 1.92e-05} & \underline{4.29e-5 $\pm$ 4.36e-5} & 1.38e-6 $\pm$ 2.30e-6 \\
U-Net & 1.11e-2 $\pm$ 2.30e-2 & 6.65e-04 $\pm$ 2.60e-03 & 7.68e-1 $\pm$ 1.24e-2 & 5.03e-5 $\pm$ 1.86e-4 \\
\hdashline
DON  & 2.76e-2 $\pm$ 9.84e-3 & 7.75e-03 $\pm$ 7.32e-03 & 1.42e-2 $\pm$ 1.47e-2 & 2.27e-5 $\pm$ 3.95e-5 \\
%POD-DON  & -- & -- & -- & -- \\
SPOD-DON  & 2.58e-2 $\pm$ 1.49e-2 & 2.35e-04 $\pm$ 9.60e-04 & 6.88e-4 $\pm$ 9.54e-4 & 2.79e-7 $\pm$ 5.71e-7 \\
TC-DON  & 2.48e-2 $\pm$ 7.67e-3 & 6.81e-03 $\pm$ 6.59e-03 & 5.92e-3 $\pm$ 3.46e-3 & 1.24e-5 $\pm$ 2.42e-5 \\
TC-SPOD-DON & 1.29e-3 $\pm$ 2.49e-3 & 2.18e-04 $\pm$ 9.16e-04 & 2.10e-3 $\pm$ 2.58e-3 & \underline{9.34e-8 $\pm$ 2.75e-7} \\
\hdashline
\textbf{SPOD-TrTINO} & \textbf{9.13e-4 $\pm$ 1.92e-3} & \underline{8.59e-05 $\pm$ 4.94e-04} & \textbf{2.95e-5 $\pm$ 1.26e-4} & \textbf{5.67e-8 $\pm$ 1.39e-7} \\
\bottomrule
\end{tabular}
} % resizebox
\end{table}

\begin{table}[h]
\centering
\caption{Error of neural operators on learning the solution operator of time-dependent 1-D advection-reaction equation, which is time-causal but not time-invariant. Results are reported as the mean squared error (MSE) $\pm$ standard deviation (SD).}
\label{tab:mse_comparison_part3}
\begin{tabular}{lc}
\toprule
& \Cref{subsubsect_adv_react} \\
%\midrule
& 1-D adv. react. \\
\midrule
FNO & \textbf{4.23e-6 $\pm$ 5.96e-6}  \\
U-Net & 4.15e-4 $\pm$ 9.34e-4 \\
\hdashline
DON  & 3.84e-3 $\pm$ 1.90e-3 \\
POD-DON  & \underline{7.78e-6 $\pm$ 1.43e-5}  \\
SPOD-DON  & 1.53e-4 $\pm$ 2.27e-4 \\
TC-DON  & 9.63e-4 $\pm$ 1.33e-3 \\
TC-SPOD-DON & 8.49e-5 $\pm$ 1.13e-4 \\
\hdashline
\textbf{SPOD-TrTINO} & 2.17e-4 $\pm$ 2.90e-4 \\
\bottomrule
\end{tabular}
\end{table}

We summarize numerical results in \Cref{tab:mse_comparison_part1} for 1-D PDEs, in \Cref{tab:mse_comparison_part2} for 2-D and 3-D PDEs, and in \Cref{tab:mse_comparison_part3} for the time-causal test case.
We note that the test problem in \Cref{subsubsect_navier_stokes_icvar} takes more training data for a prediction on a shorter time interval, which explains the higher accuracy compared to the one in \Cref{subsubsect_navier_stokes}.

When the solution operator is time-invariant in the underlying problem, the overall accuracy of TC variants of DeepONet has superior accuracy compared to vanilla DeepONet versions, and SPOD-TrTINO contributes further improvement.
In some of the test cases, SPOD-TrTINO has comparable precision when compared to a reference FNO implementation, which takes longer per-epoch time to train despite having less network parameters.

When the dynamics is time-causal but not time-invariant, as in \Cref{subsubsect_adv_react}, the results shown in \Cref{tab:mse_comparison_part3} show clearly different trend compared to other 1-D PDE test cases shown by \Cref{tab:mse_comparison_part1}.
In this test case, SPOD-TrTINO yields the worst accuracy among POD-based neural operators, penalized by the wrong assumption of time invariance.

For the 2-D Navier-Stokes equation, all models perform better in the setup with the variational initial condition. As the prediction horizon is reduced to 0.2T for each sub-interval, the cumulative errors associated with long-term temporal integration are drastically reduced. Moreover, conditioning the model on the exact intermediate states transforms the task from learning a global trajectory into learning a short-horizon, local flow map, thereby passing the error accumulation and initial transient challenges present in the full-horizon setup.

We visualize some predictions of SPOD-TrTINO in \Cref{fig:pde_allen_cahn_four_subfigures,fig:pde_burgers_nu01_four_subfigures,fig:pde_burgers_nu001_four_subfigures,fig:pde_convection_diffusion_four_subfigures,fig:ns_results_matrix,fig:heat_results_matrix}.

\begin{figure}[h!]
    \centering
    \begin{subfigure}{0.2\textwidth}
        \centering
        \includegraphics[width=\linewidth]{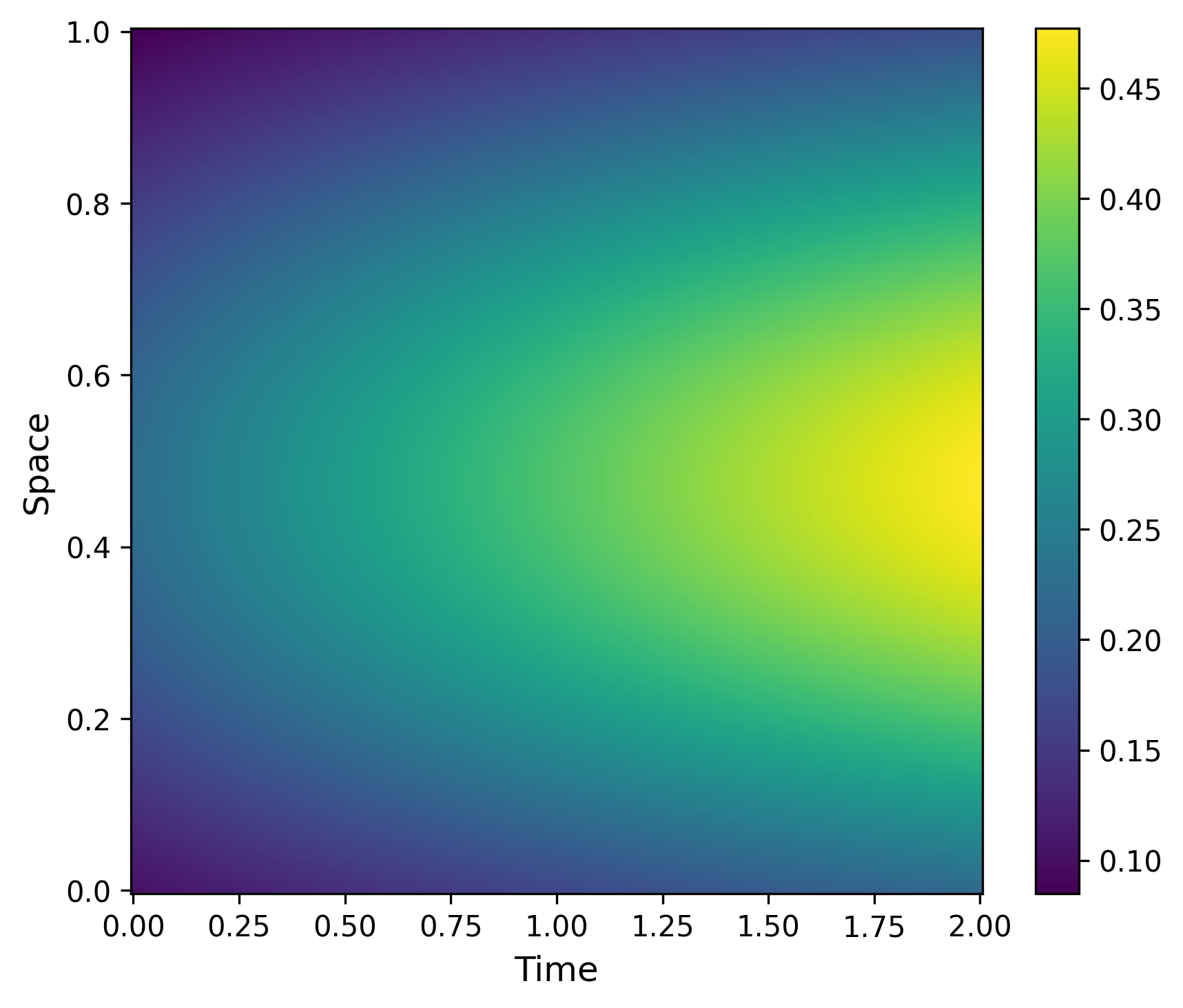}
        \caption{Input function}
        \label{fig:pde_allen_cahn_input_func}
    \end{subfigure}
    \hfill
    \begin{subfigure}{0.2\textwidth}
        \centering
        \includegraphics[width=\linewidth]{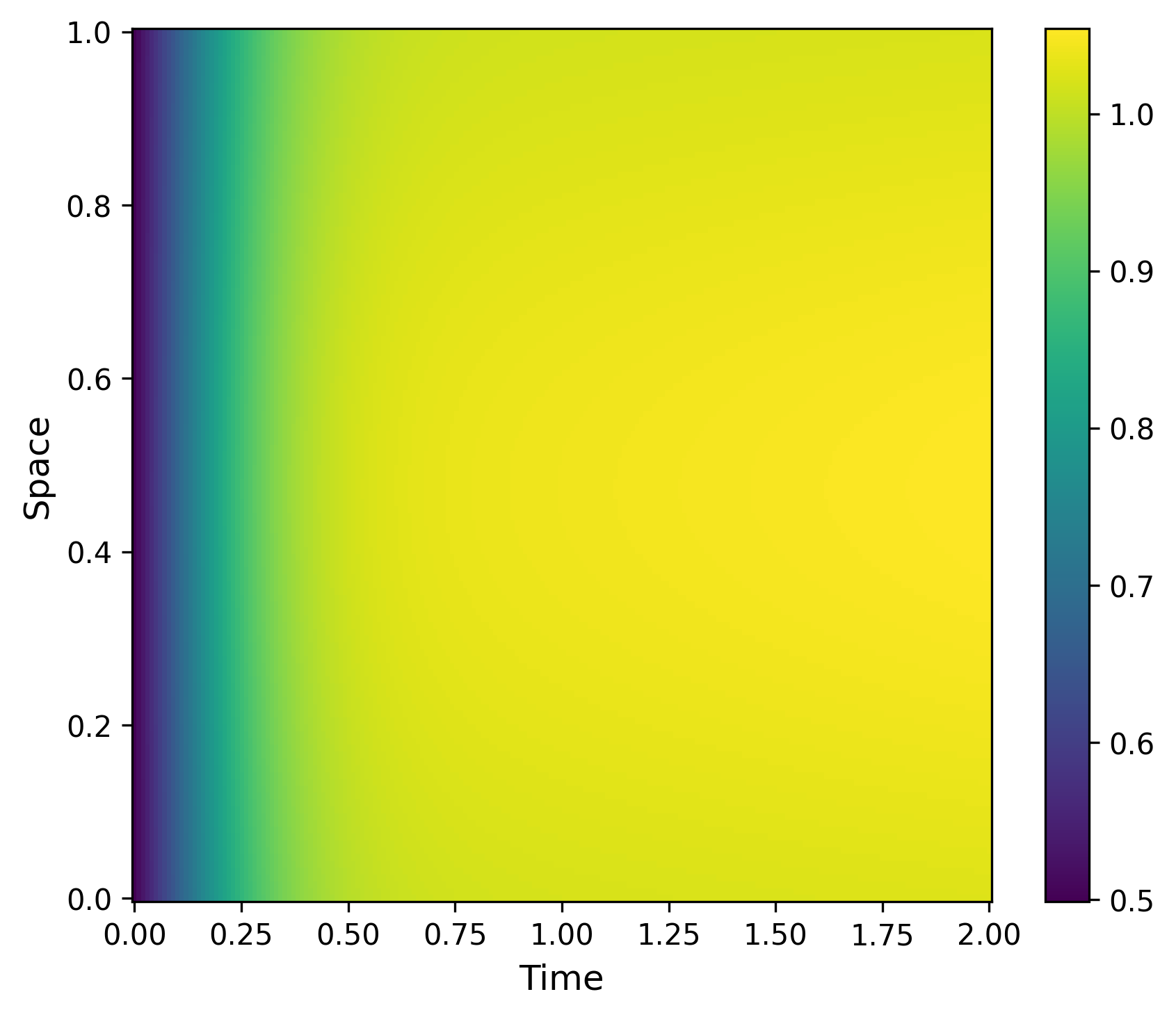}
        \caption{Ground truth}
        \label{fig:pde_allen_cahn_ground_truth}
    \end{subfigure}
    \hfill
    \begin{subfigure}{0.2\textwidth}
        \centering
        \includegraphics[width=\linewidth]{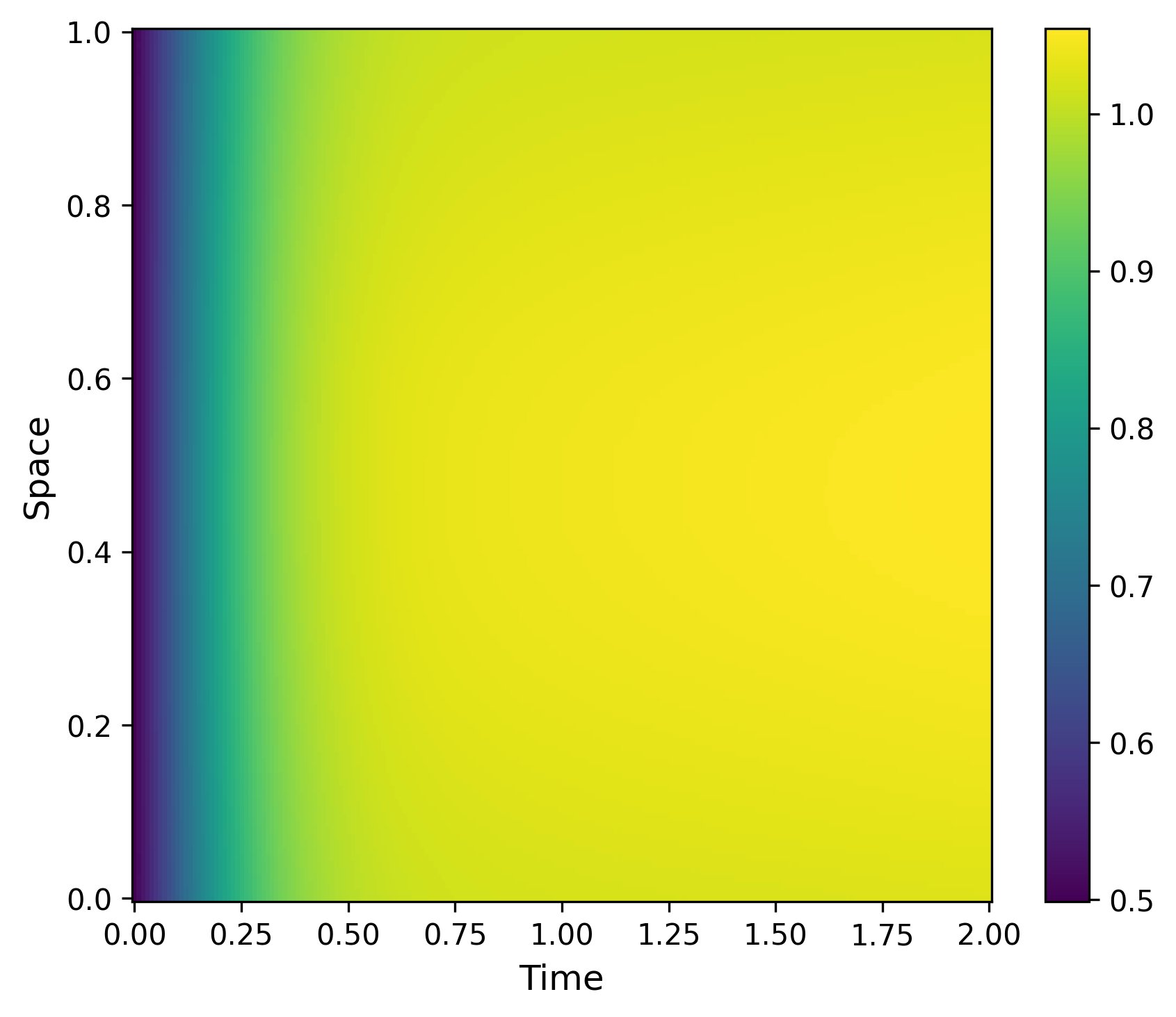}
        \caption{Prediction}
        \label{fig:pde_allen_cahn_result}
    \end{subfigure}
    \hfill
    \begin{subfigure}{0.2\textwidth}
        \centering
        \includegraphics[width=\linewidth]{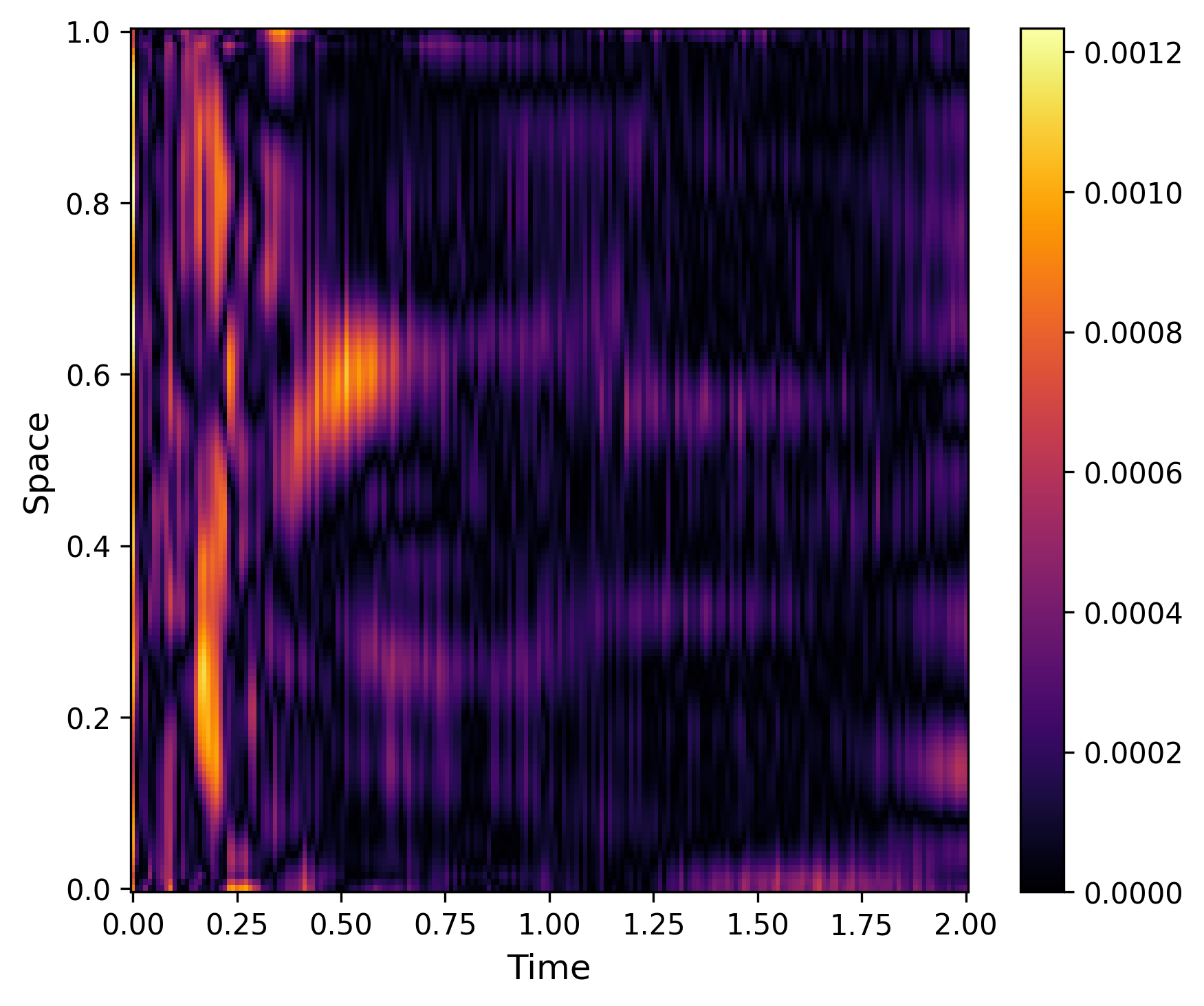}
        \caption{Absolute error}
        \label{fig:pde_allen_cahn_error}
    \end{subfigure}
    \caption{Example output of SPOD-TrTINO for learning the solution of Allen--Cahn equation (\Cref{subsubsect_allen_cahn}).}
    \label{fig:pde_allen_cahn_four_subfigures}
\end{figure}

\begin{figure}[h]
    \centering
    \begin{subfigure}{0.2\textwidth}
        \centering
        \includegraphics[width=\linewidth]{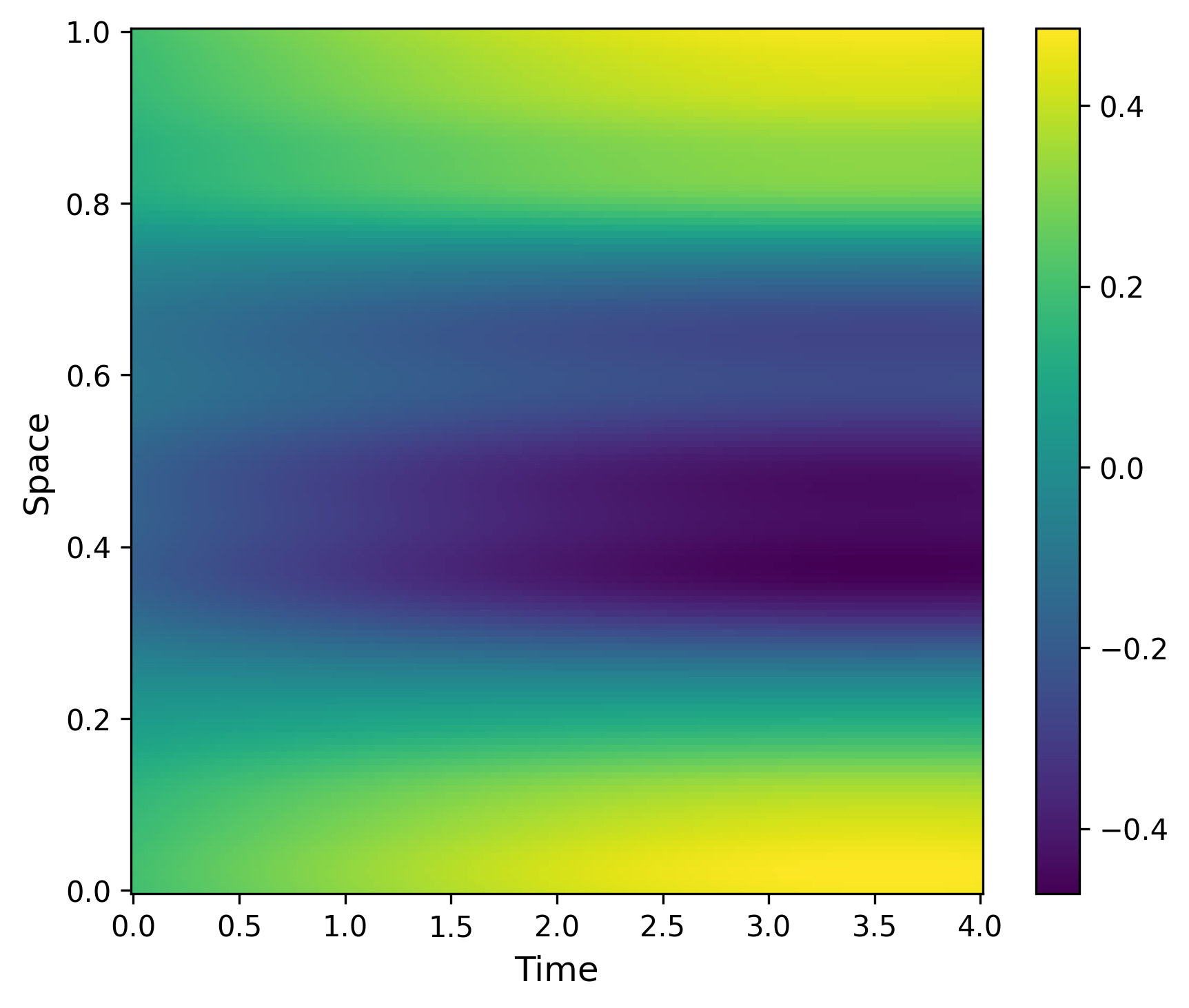}
        \caption{Input function}
        \label{fig:pde_burgers_nu01_input_func}
    \end{subfigure}
    \hfill
    \begin{subfigure}{0.2\textwidth}
        \centering
        \includegraphics[width=\linewidth]{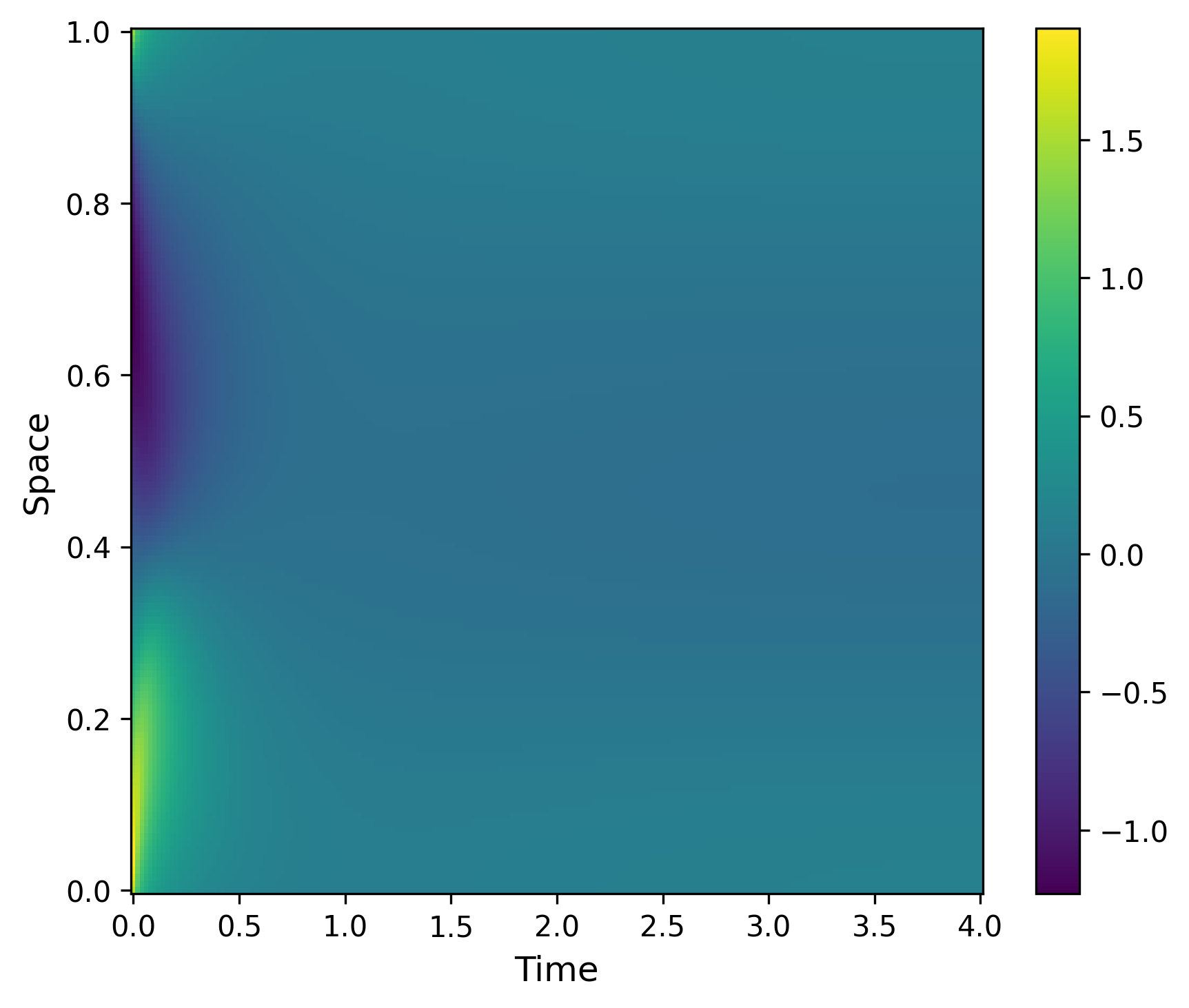}
        \caption{Ground truth}
        \label{fig:pde_burgers_nu01_ground_truth}
    \end{subfigure}
    \hfill
    \begin{subfigure}{0.2\textwidth}
        \centering
        \includegraphics[width=\linewidth]{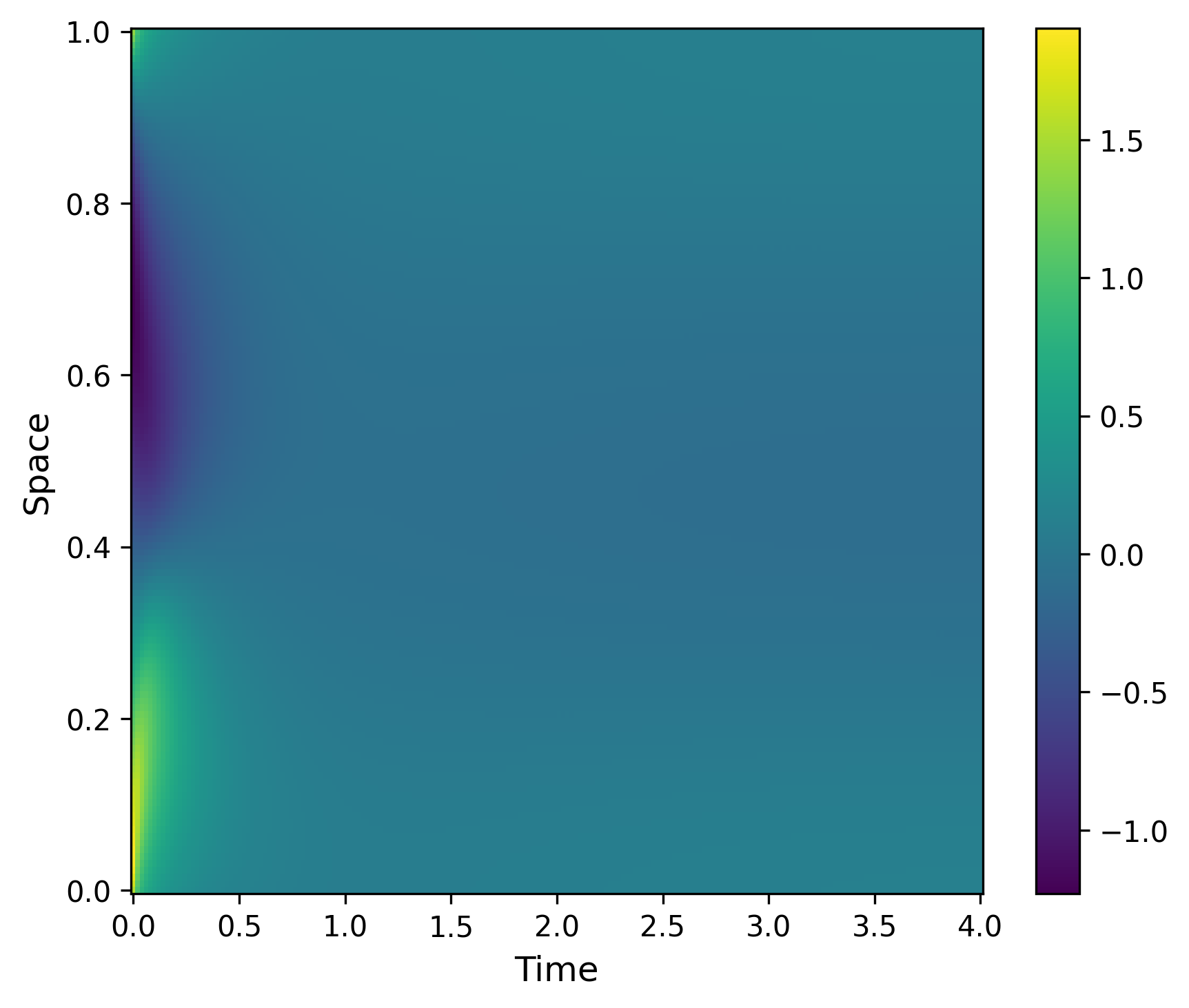}
        \caption{Prediction}
        \label{fig:pde_burgers_nu01_result}
    \end{subfigure}
    \hfill
    \begin{subfigure}{0.2\textwidth}
        \centering
        \includegraphics[width=\linewidth]{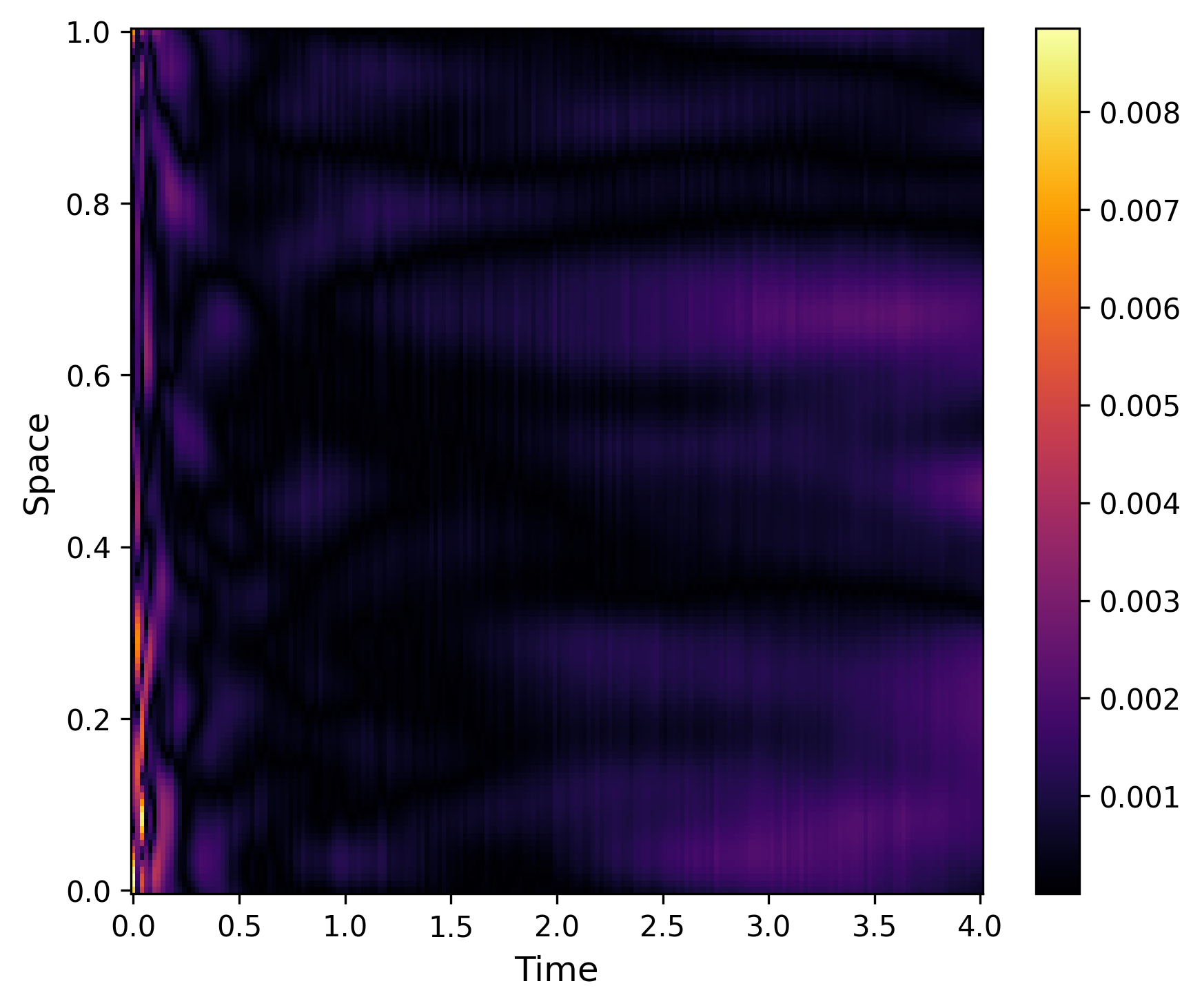}
        \caption{Absolute error}
        \label{fig:pde_burgers_nu01_error}
    \end{subfigure}
    \caption{Example output of SPOD-TrTINO for learning the solution of Burgers' equation with $\nu=0.1$ (\Cref{subsubsect_burgers}).}
    \label{fig:pde_burgers_nu01_four_subfigures}
\end{figure}

\begin{figure}[h]
    \centering
    \begin{subfigure}{0.2\textwidth}
        \centering
        \includegraphics[width=\linewidth]{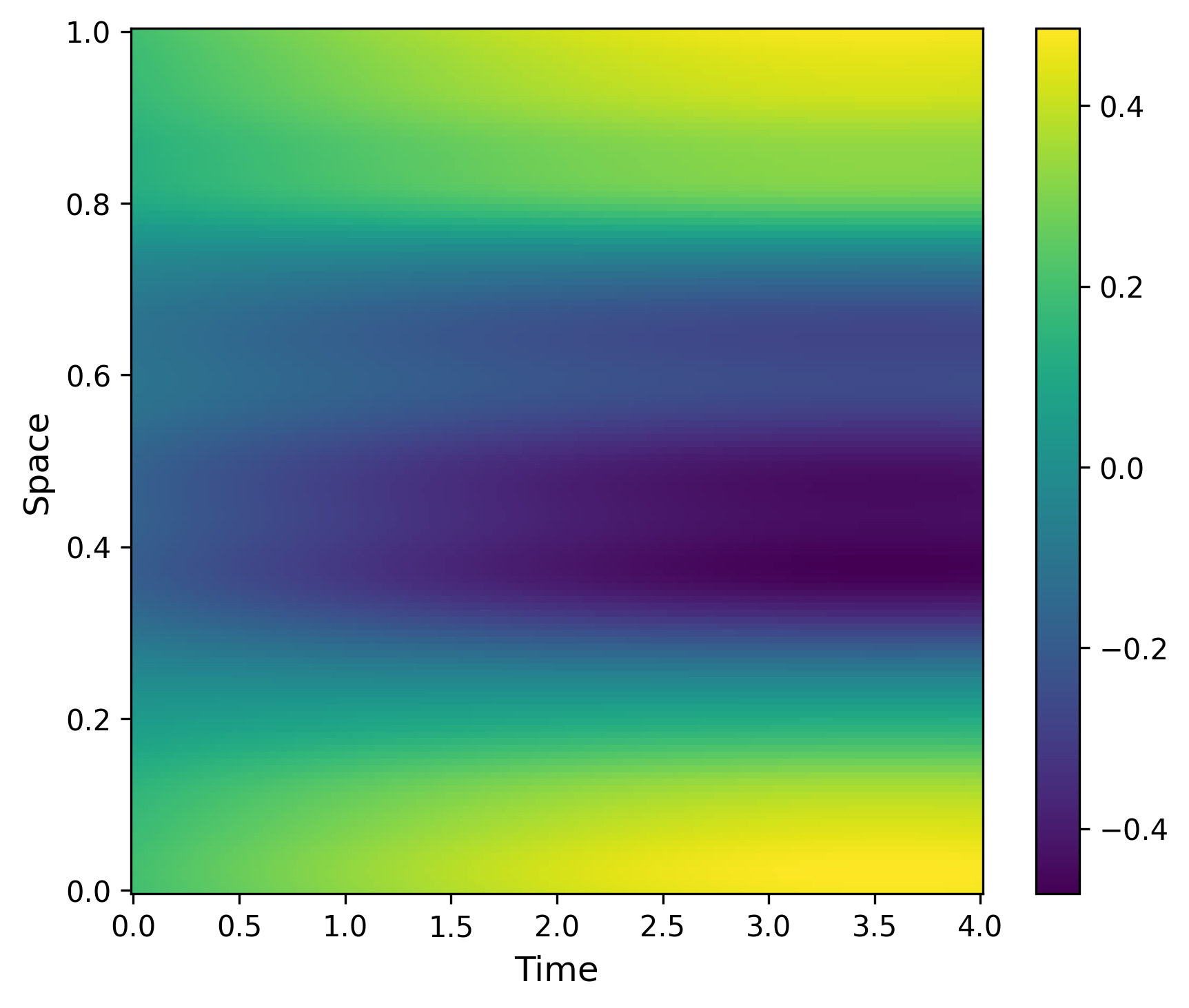}
        \caption{Input function}
        \label{fig:pde_burgers_nu001_input_func}
    \end{subfigure}
    \hfill
    \begin{subfigure}{0.2\textwidth}
        \centering
        \includegraphics[width=\linewidth]{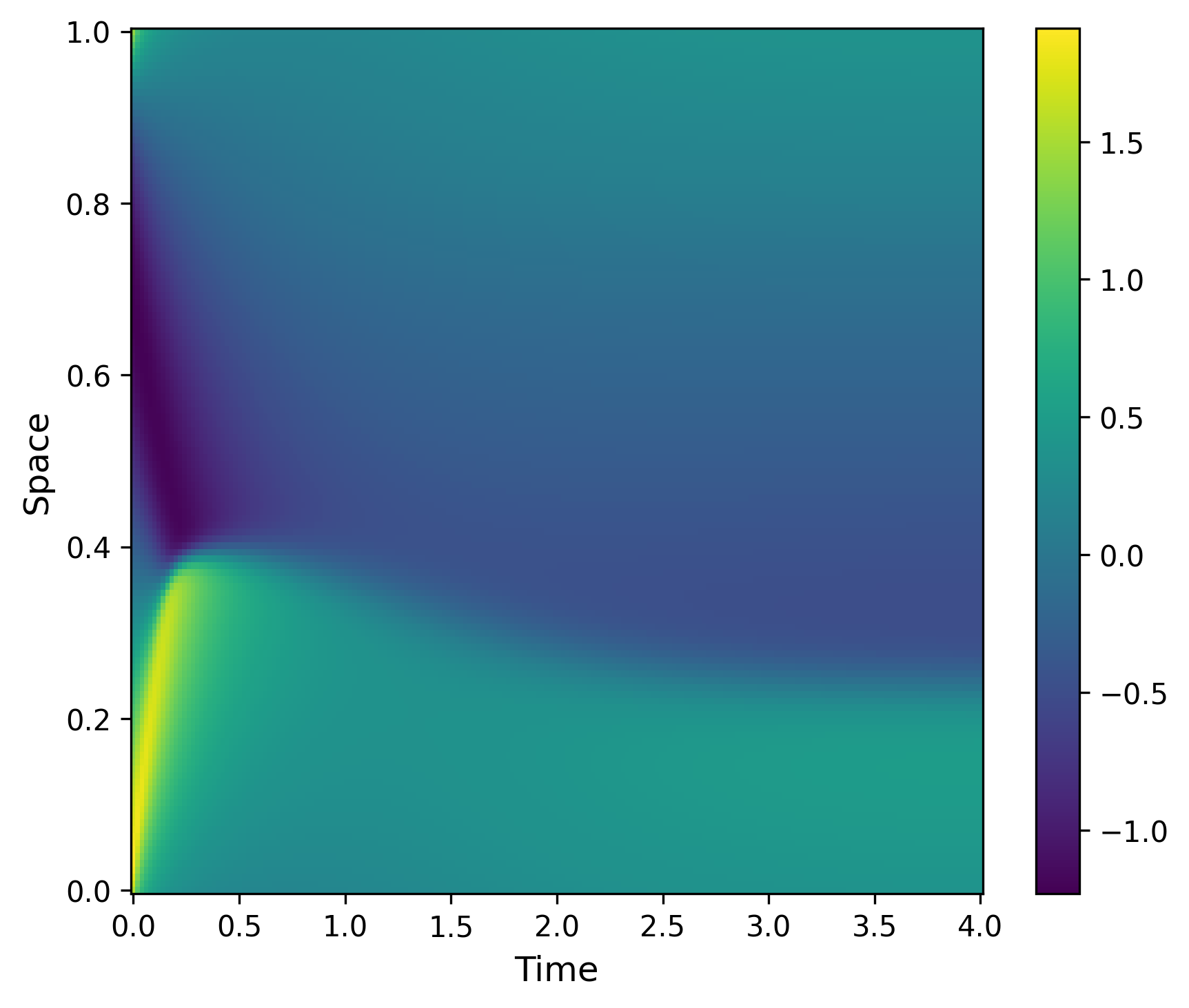}
        \caption{Ground truth}
        \label{fig:pde_burgers_nu001_ground_truth}
    \end{subfigure}
    \hfill
    \begin{subfigure}{0.2\textwidth}
        \centering
        \includegraphics[width=\linewidth]{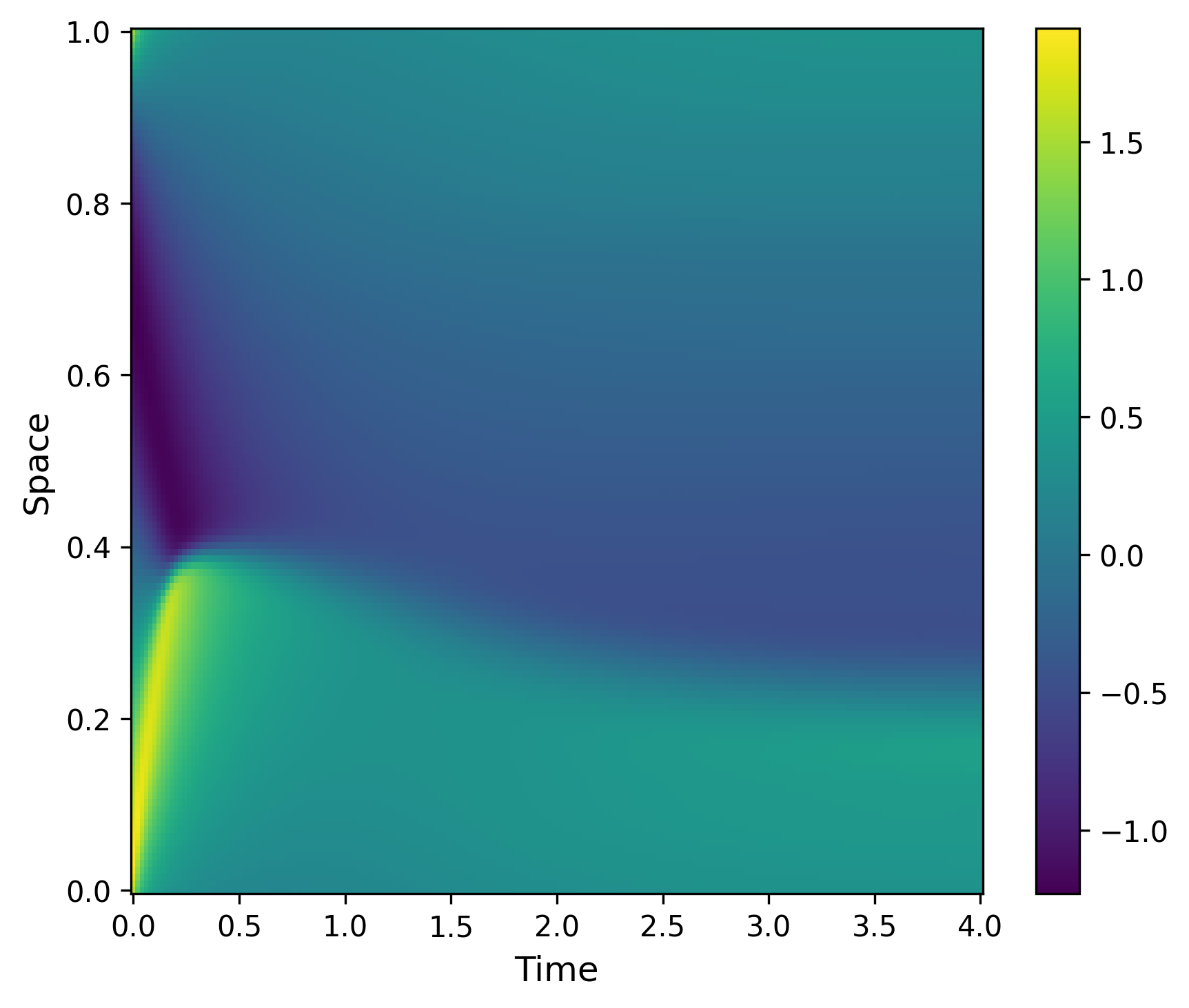}
        \caption{Prediction}
        \label{fig:pde_burgers_nu001_result}
    \end{subfigure}
    \hfill
    \begin{subfigure}{0.2\textwidth}
        \centering
        \includegraphics[width=\linewidth]{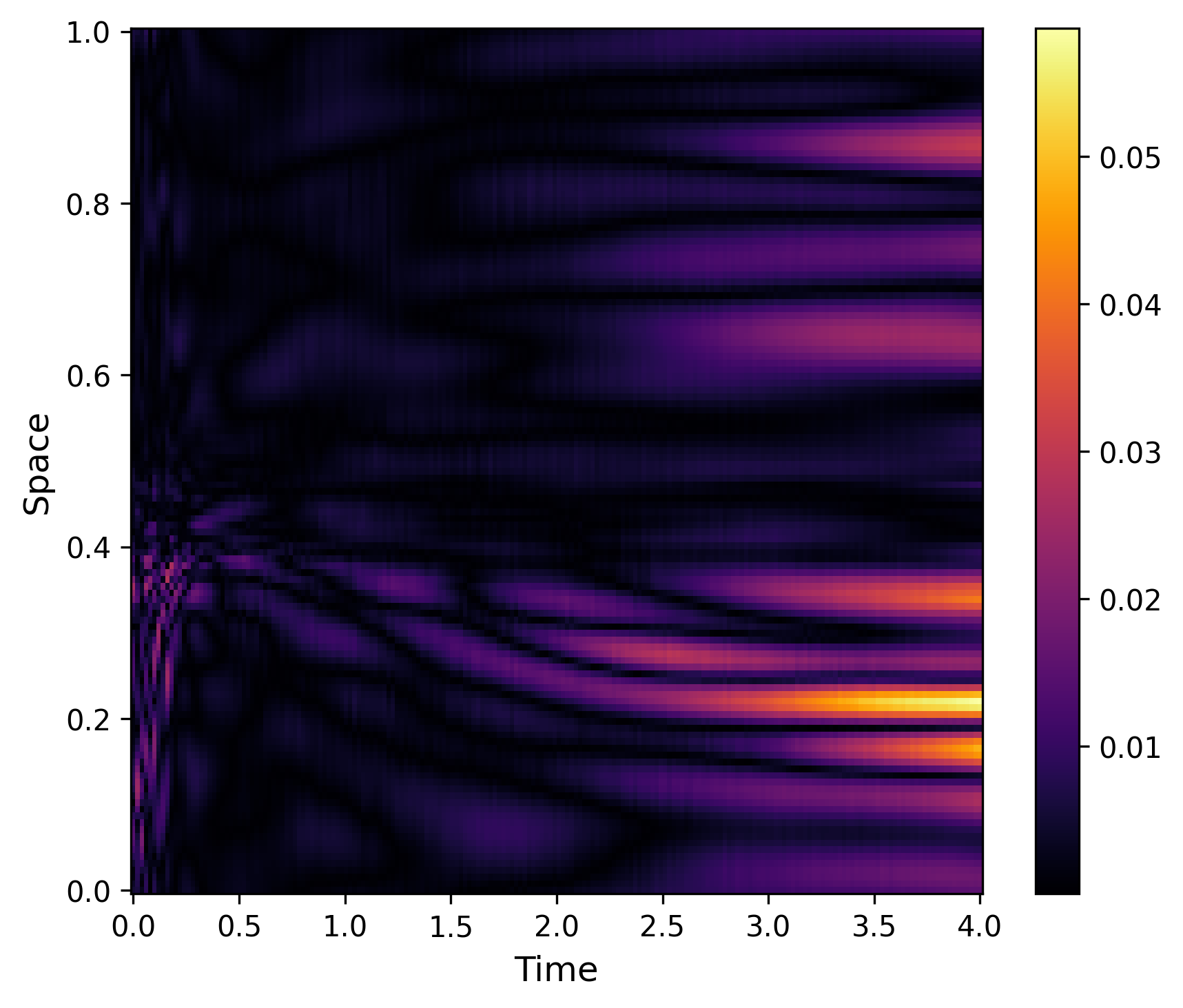}
        \caption{Absolute error}
        \label{fig:pde_burgers_nu001_error}
    \end{subfigure}
    \caption{Example output of SPOD-TrTINO for learning the solution of Burgers' equation with $\nu=0.01$ (\Cref{subsubsect_burgers}).}
    \label{fig:pde_burgers_nu001_four_subfigures}
\end{figure}

\begin{figure}[h]
    \centering
    \begin{subfigure}{0.2\textwidth}
        \centering
        \includegraphics[width=\linewidth]{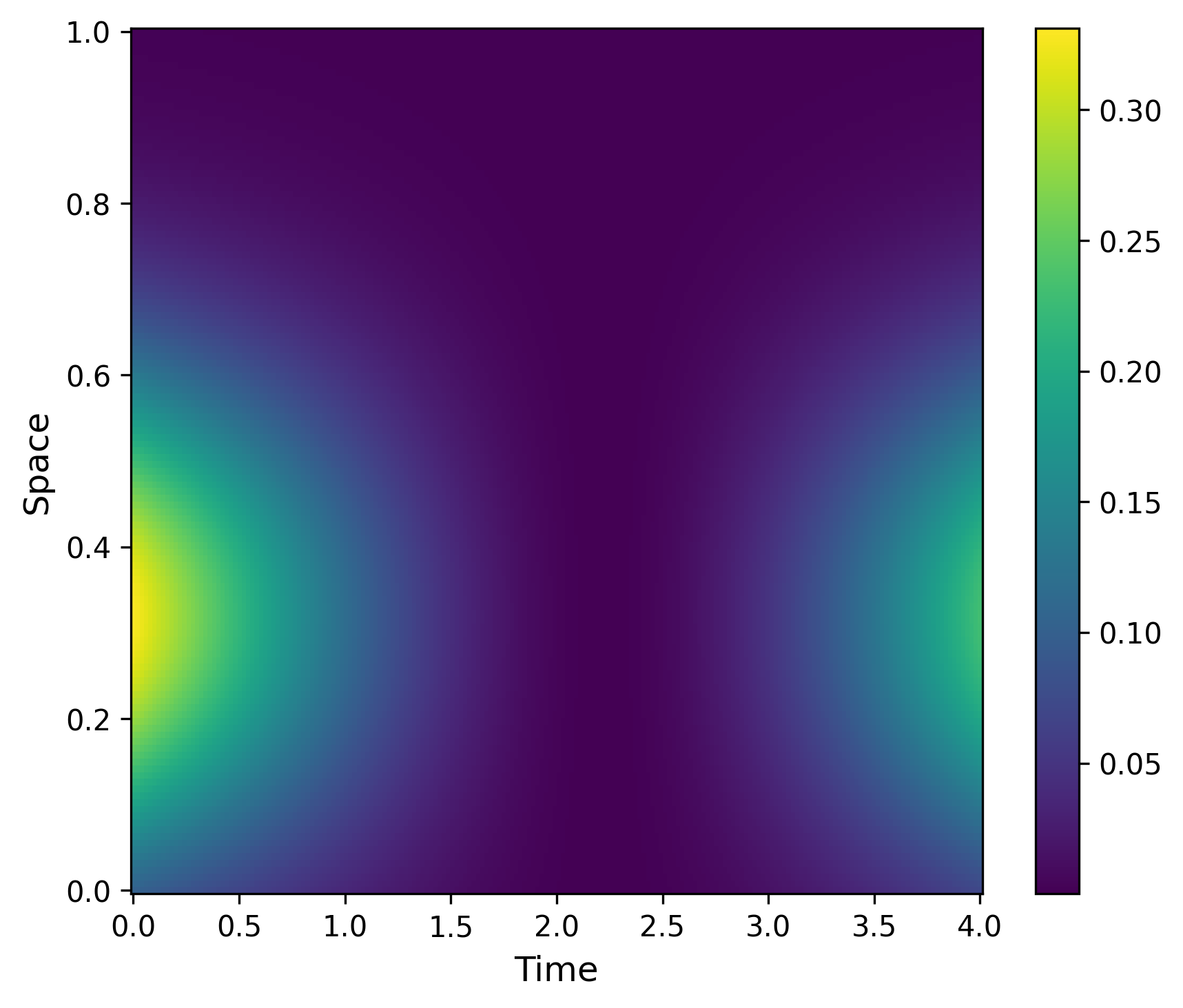}
        \caption{Input function}
        \label{fig:pde_convection_diffusion_input_func}
    \end{subfigure}
    \hfill
    \begin{subfigure}{0.2\textwidth}
        \centering
        \includegraphics[width=\linewidth]{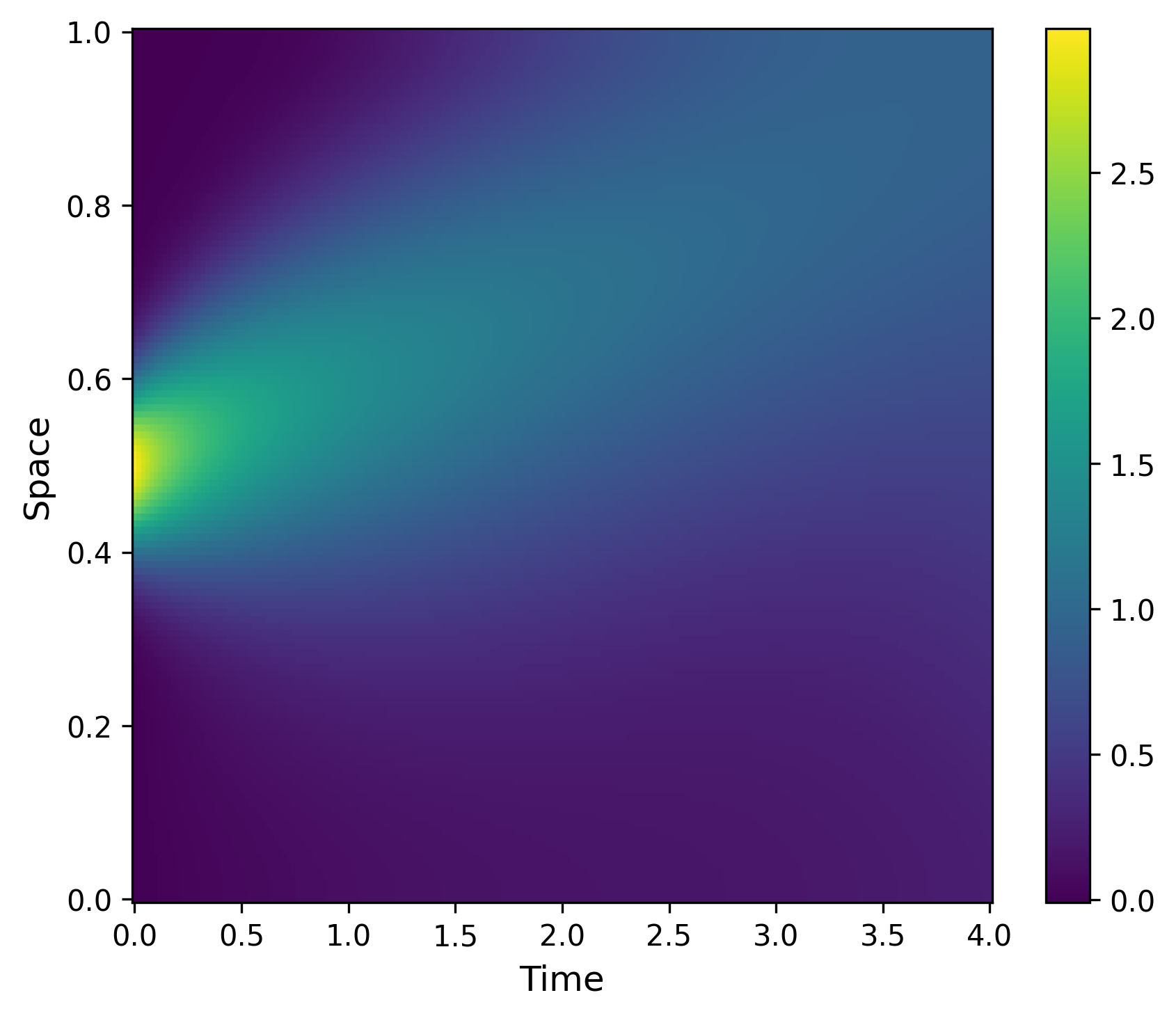}
        \caption{Ground truth}
        \label{fig:pde_convection_diffusion_ground_truth}
    \end{subfigure}
    \hfill
    \begin{subfigure}{0.2\textwidth}
        \centering
        \includegraphics[width=\linewidth]{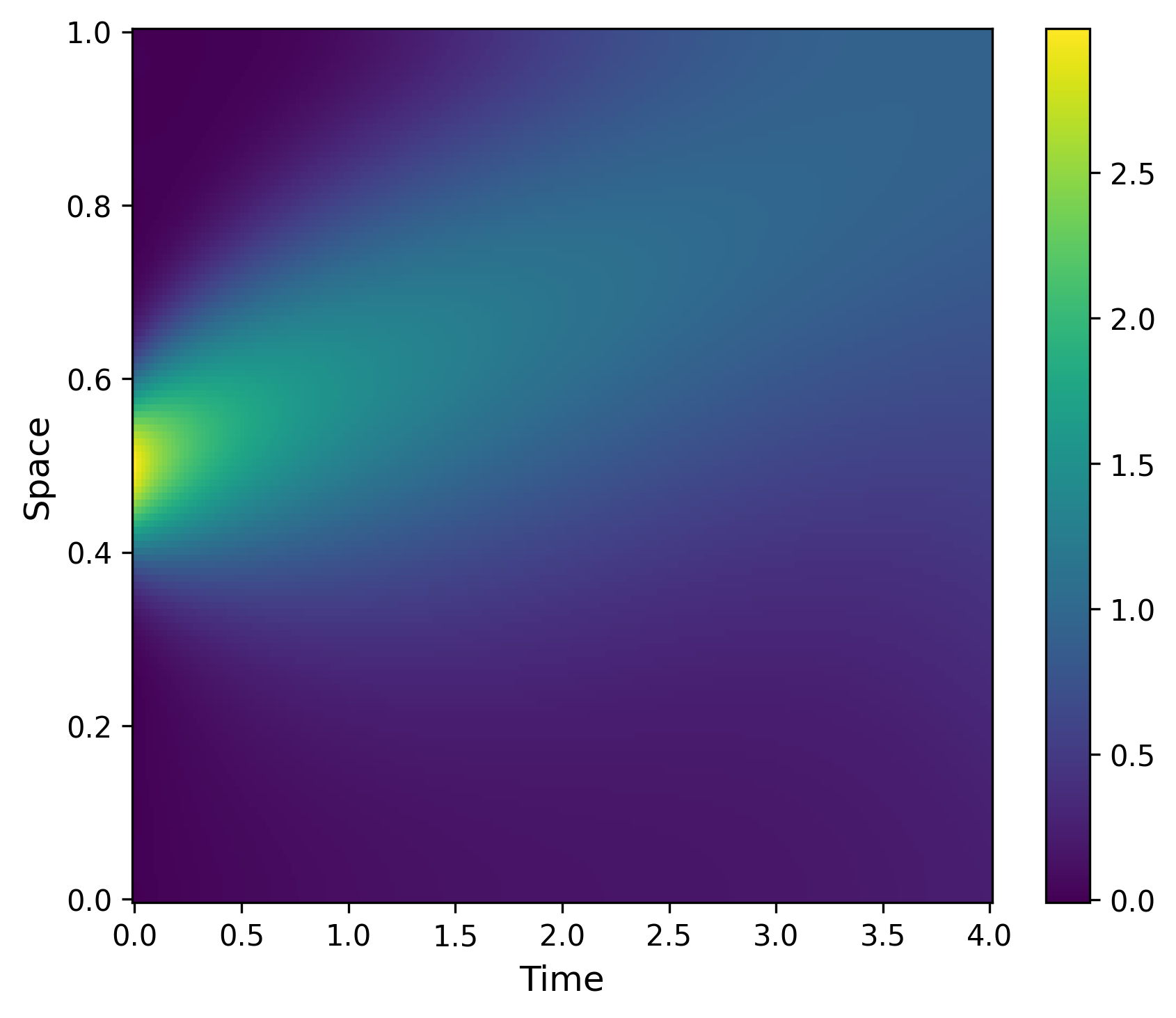}
        \caption{Prediction}
        \label{fig:pde_convection_diffusion_result}
    \end{subfigure}
    \hfill
    \begin{subfigure}{0.2\textwidth}
        \centering
        \includegraphics[width=\linewidth]{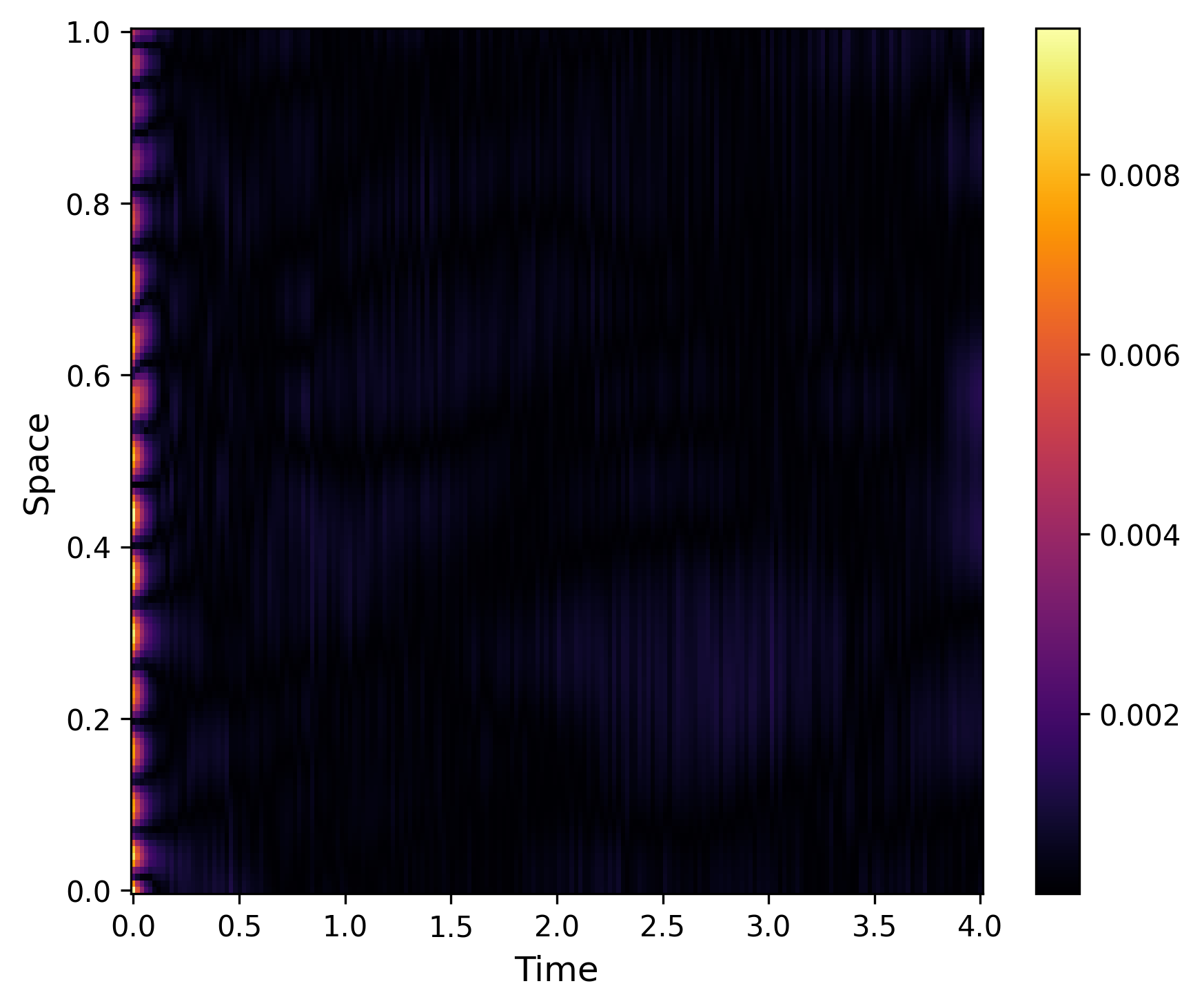}
        \caption{Absolute error}
        \label{fig:pde_convection_diffusion_error}
    \end{subfigure}
    \caption{Example output of SPOD-TrTINO for learning the solution of convection-diffusion equation (\Cref{subsubsect_conv_diff}).}
    \label{fig:pde_convection_diffusion_four_subfigures}
\end{figure}

\begin{figure}[htpt]
	\centering
	\setlength{\tabcolsep}{1.5pt} 
	
	\begin{tabular}{>{\centering\arraybackslash}m{0.1\textwidth} 
			*{6}{>{\centering\arraybackslash}m{0.14\textwidth}}}

		& $t=0.0$ & $t=0.4$ & $t=0.8$ & $t=1.2$ & $t=1.6$ & $t=2.0$ \\
		\noalign{\vspace{4pt}} 
		
		Input & 
		\includegraphics[width=\linewidth]{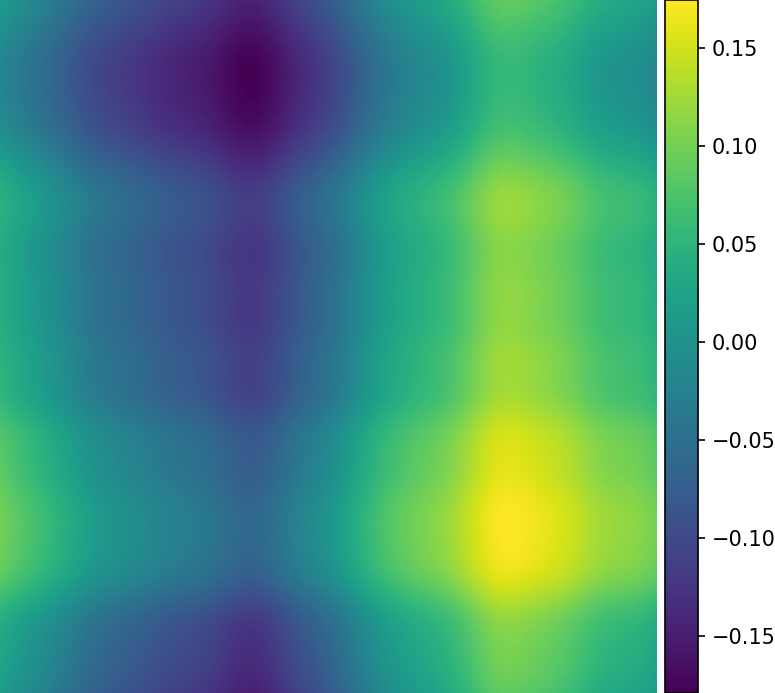} &
		\includegraphics[width=\linewidth]{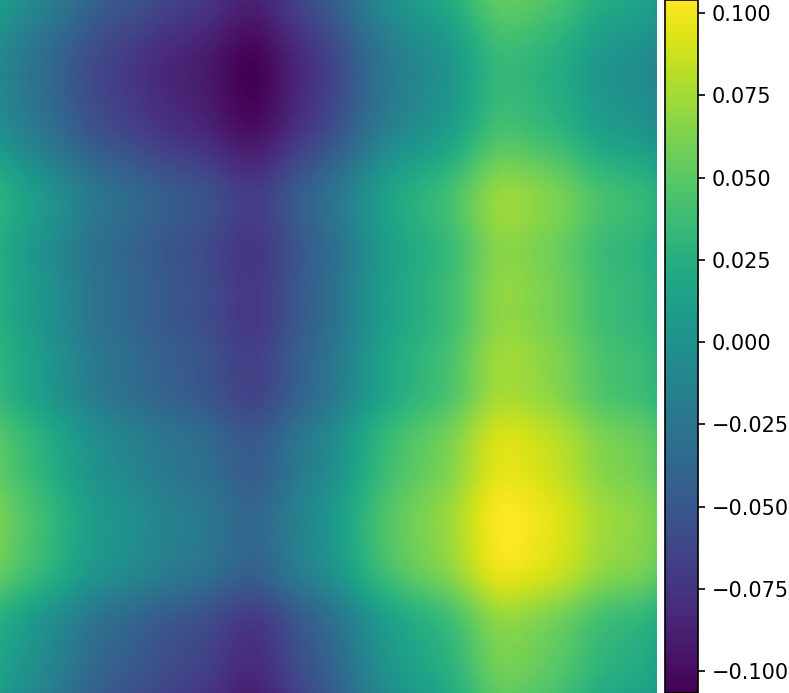} &
		\includegraphics[width=\linewidth]{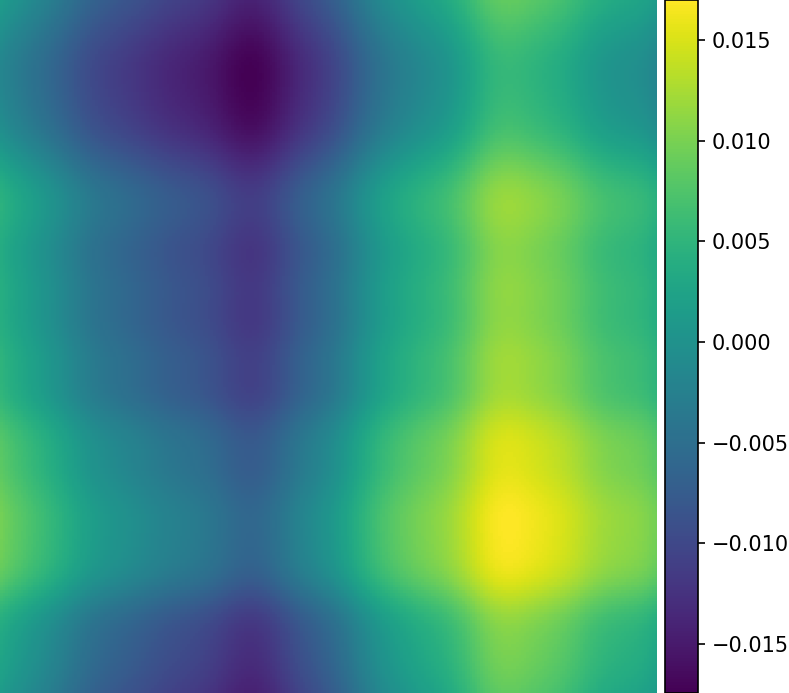} &
		\includegraphics[width=\linewidth]{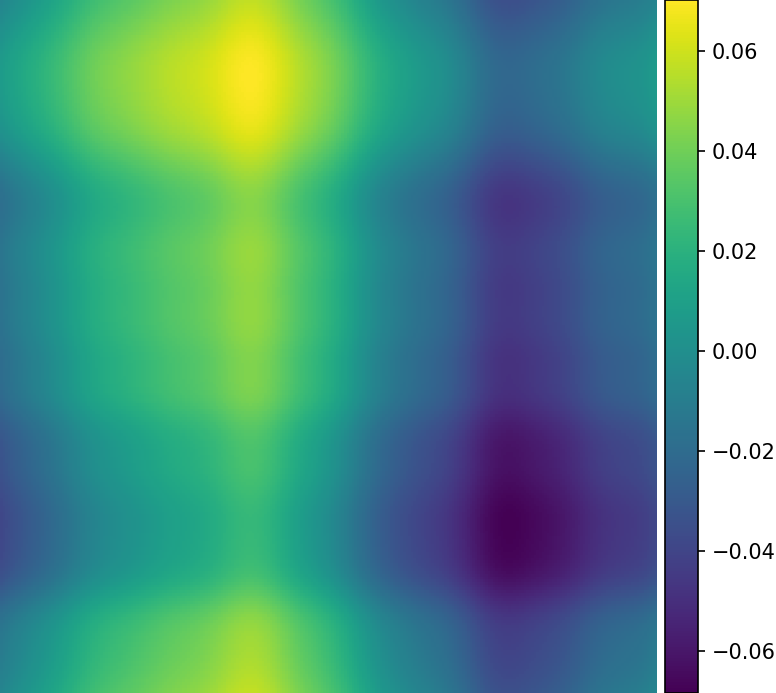} &
		\includegraphics[width=\linewidth]{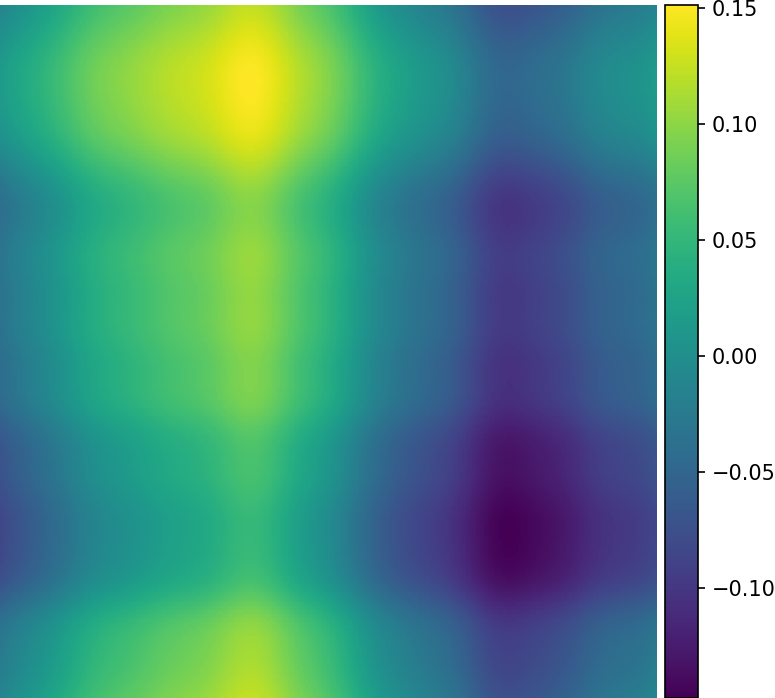} &
		\includegraphics[width=\linewidth]{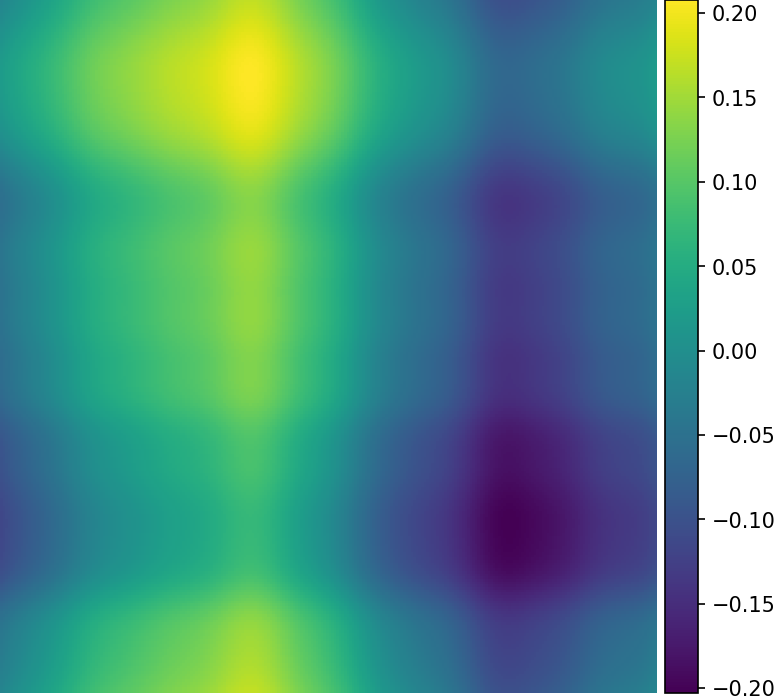} \\
		\noalign{\vspace{6pt}}

		Ground Truth & 
		\includegraphics[width=\linewidth]{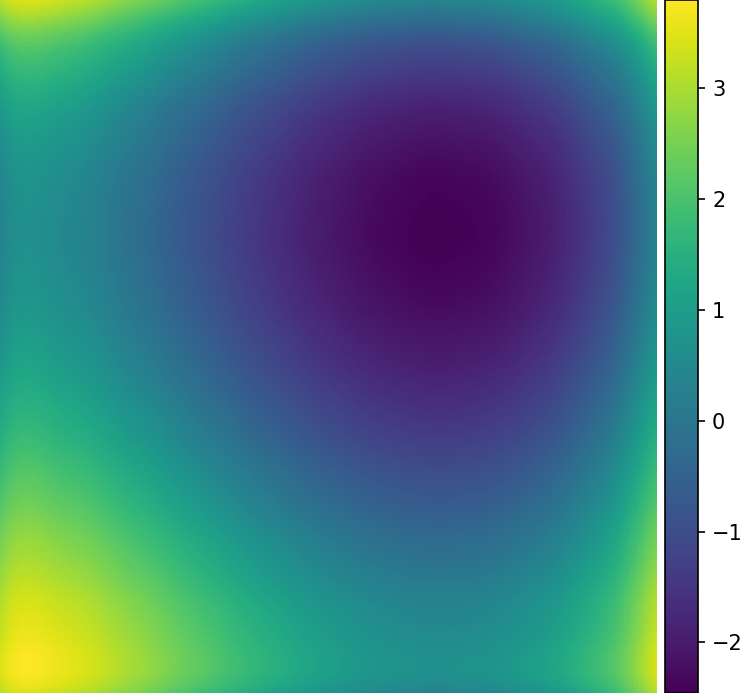} &
		\includegraphics[width=\linewidth]{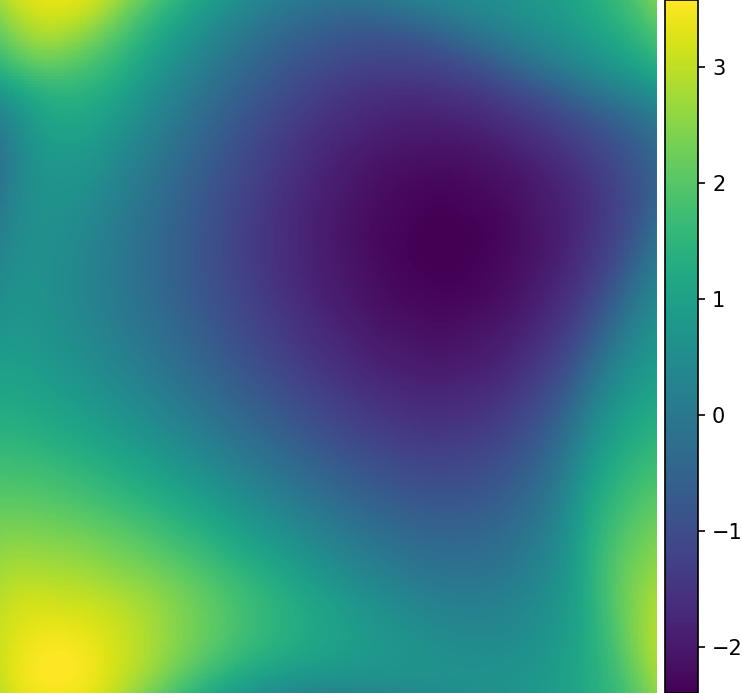} &
		\includegraphics[width=\linewidth]{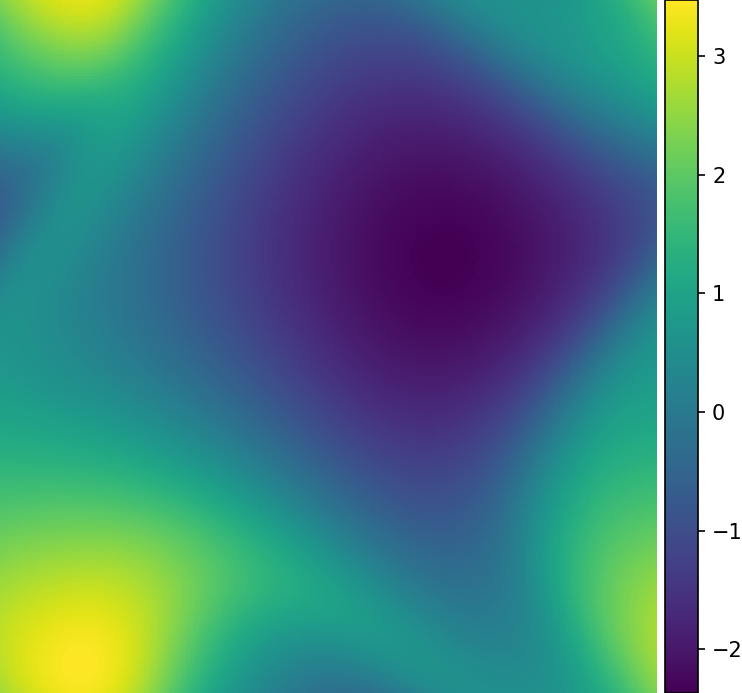} &
		\includegraphics[width=\linewidth]{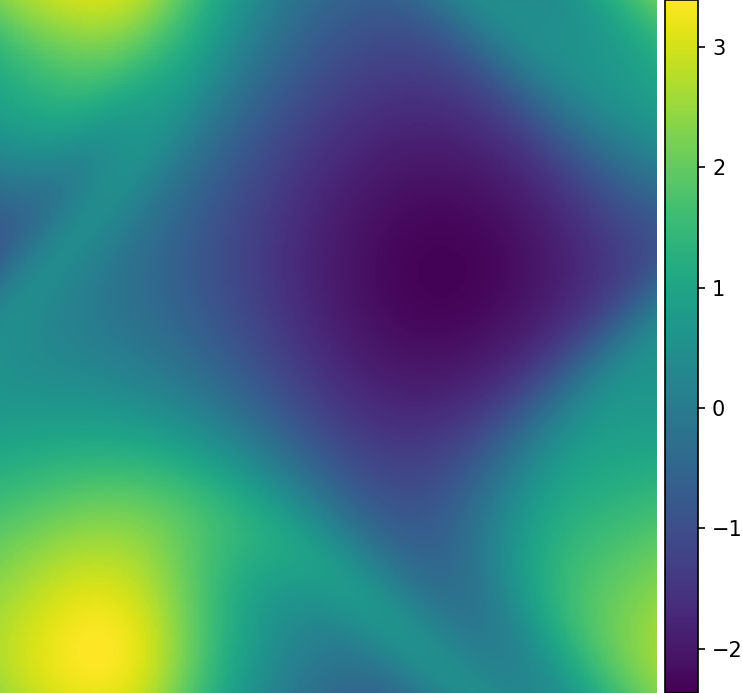} &
		\includegraphics[width=\linewidth]{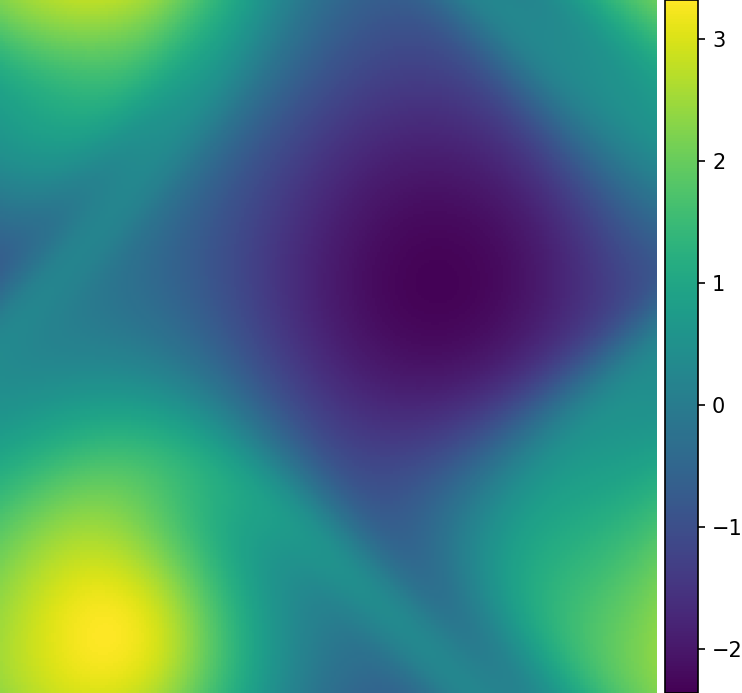} &
		\includegraphics[width=\linewidth]{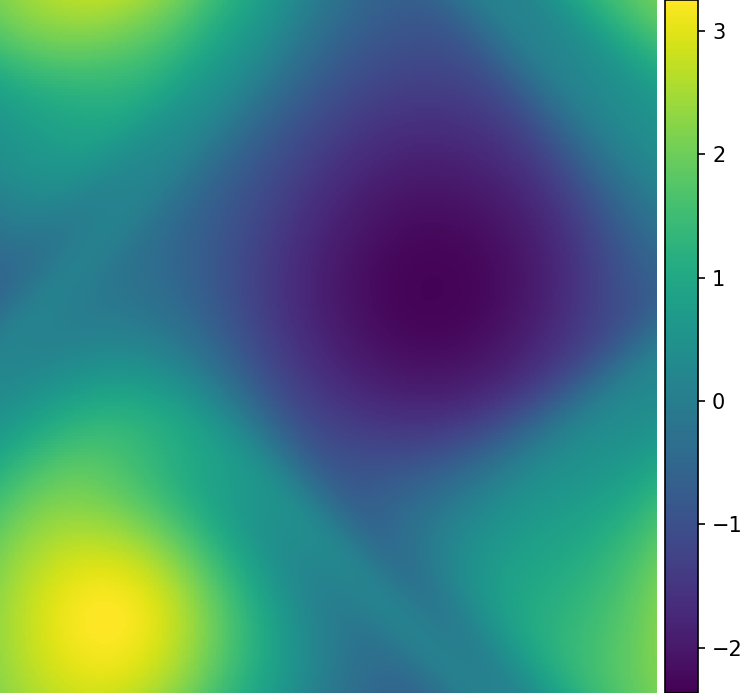} \\
		\noalign{\vspace{6pt}}

		SPOD-TrTINO Prediction & 
		\includegraphics[width=\linewidth]{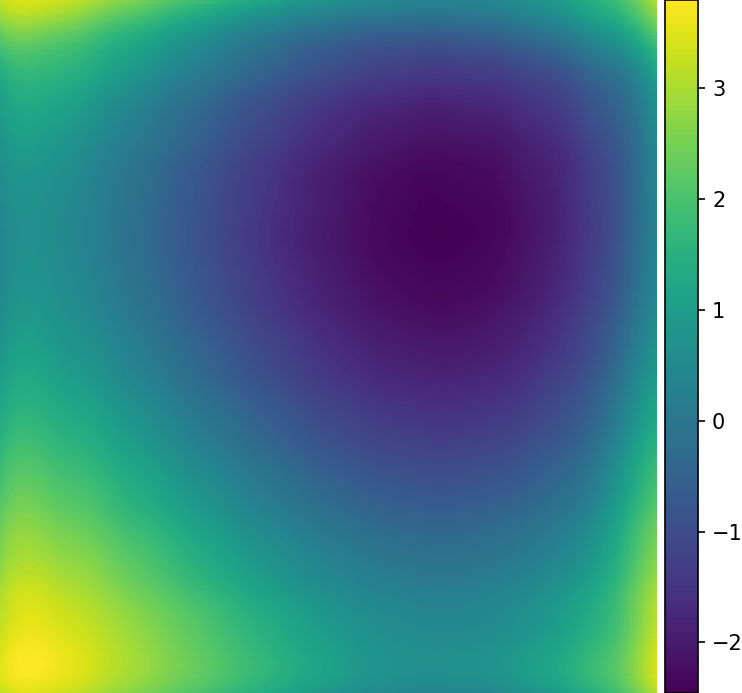} &
		\includegraphics[width=\linewidth]{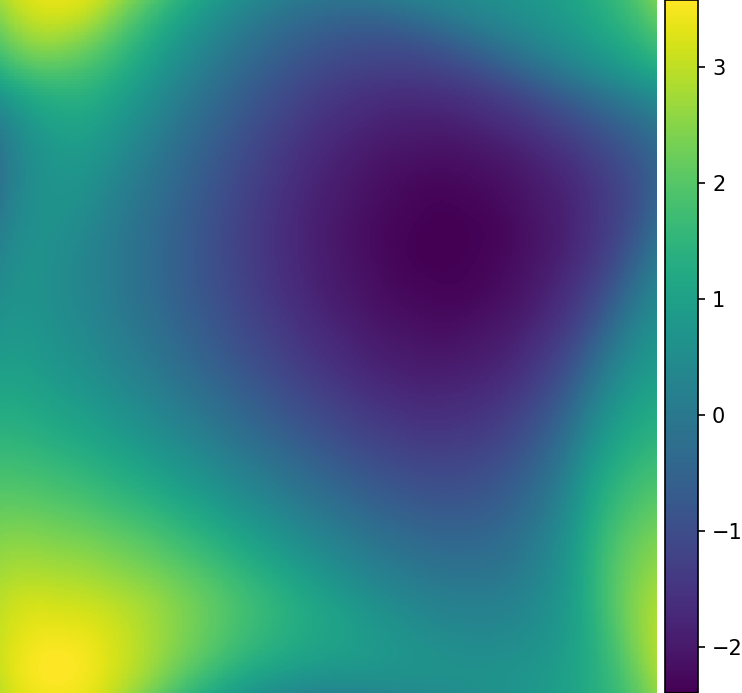} &
		\includegraphics[width=\linewidth]{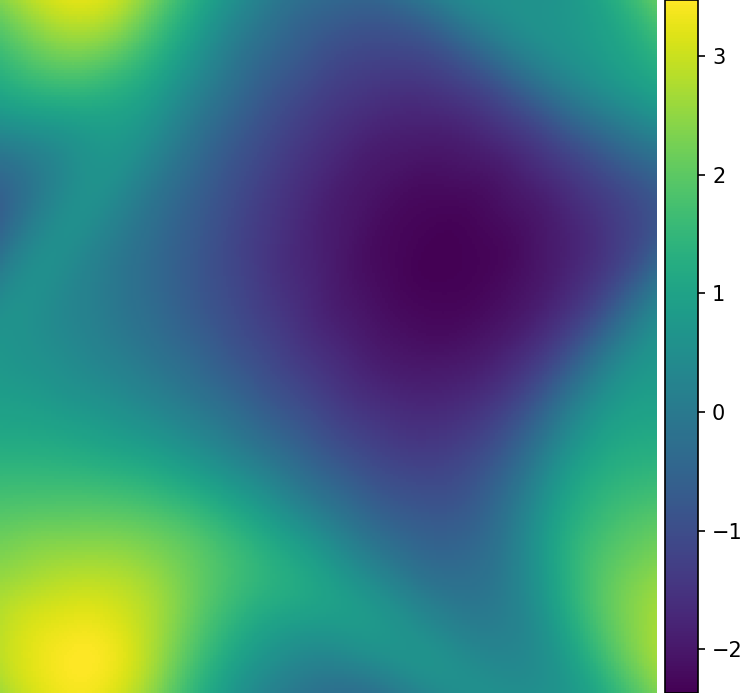} &
		\includegraphics[width=\linewidth]{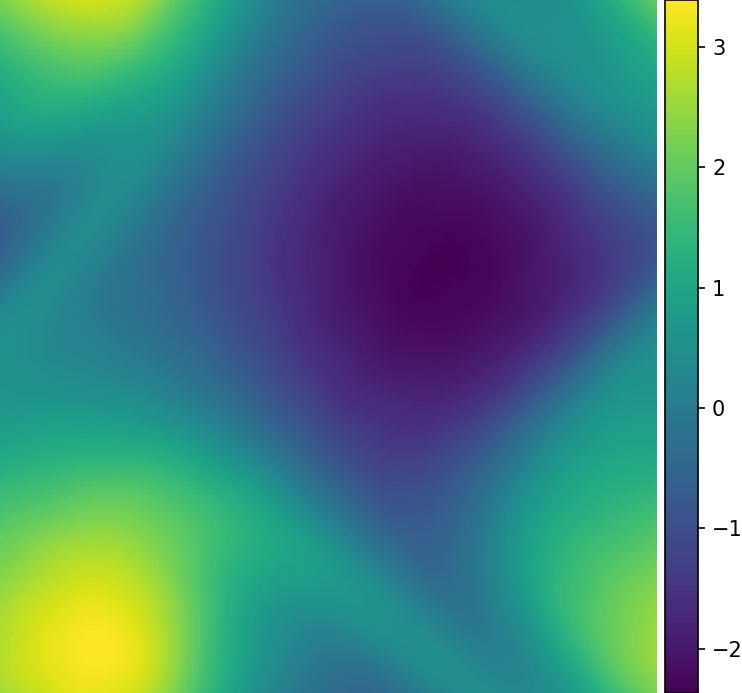} &
		\includegraphics[width=\linewidth]{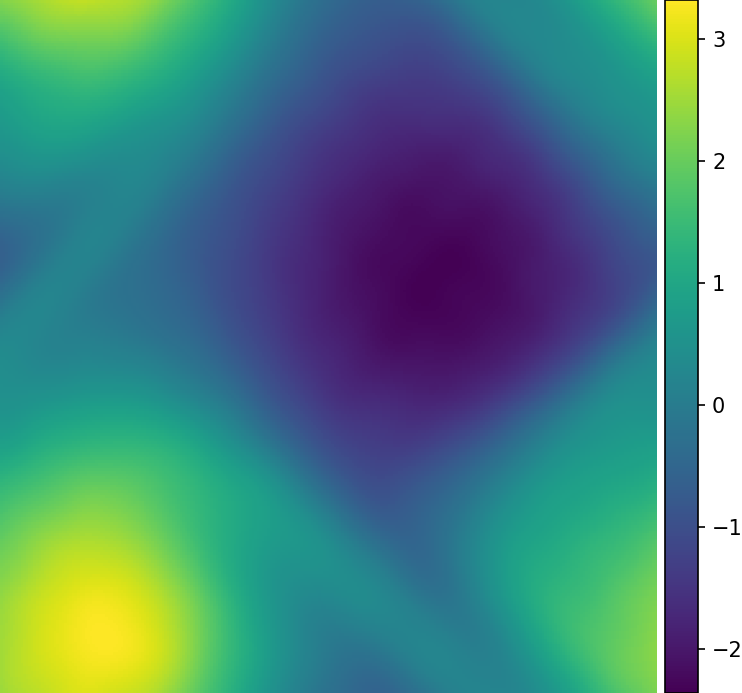} &
		\includegraphics[width=\linewidth]{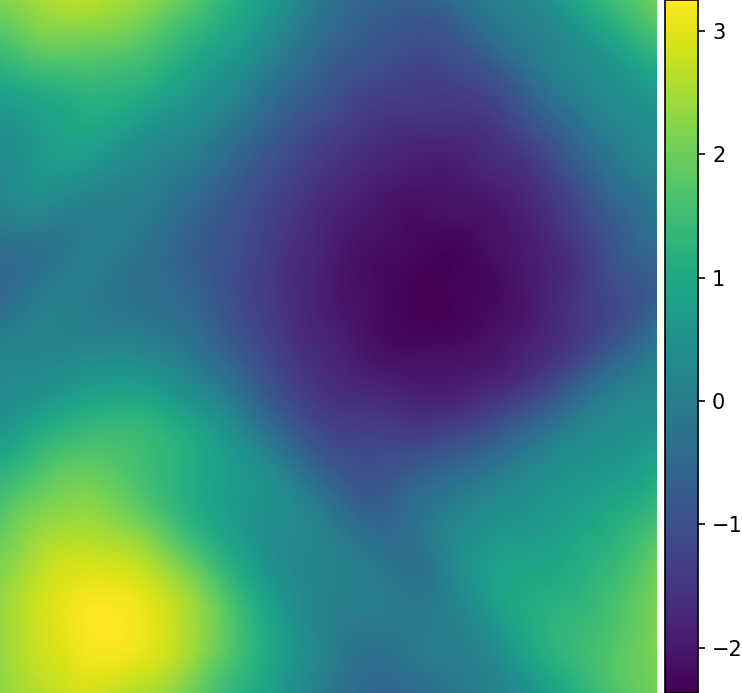} \\
		\noalign{\vspace{6pt}}
		
		Absolute Error & 
		\includegraphics[width=\linewidth]{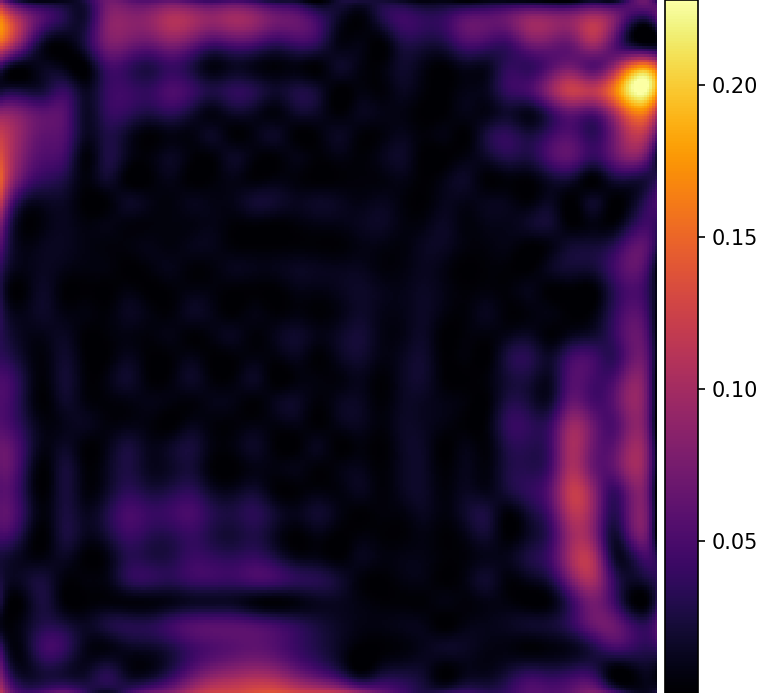} &
		\includegraphics[width=\linewidth]{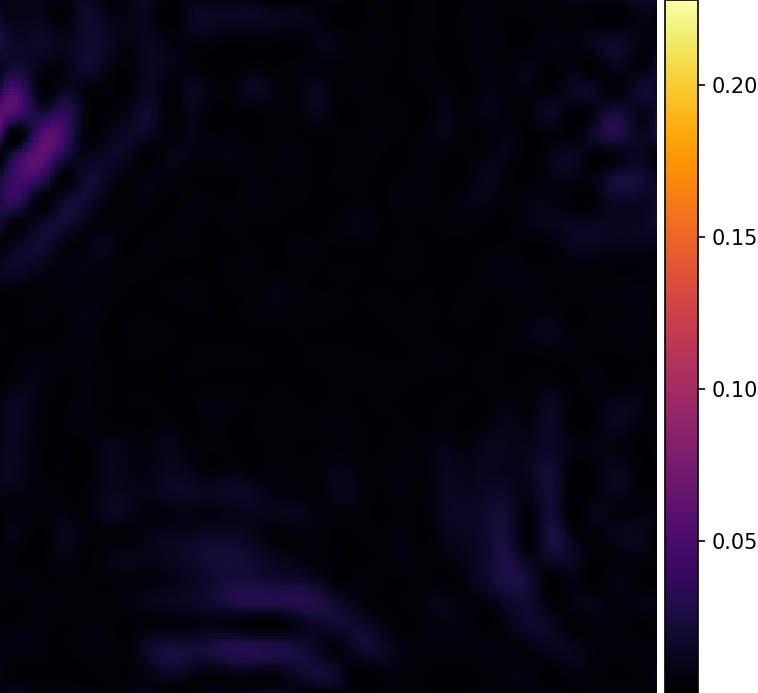} &
		\includegraphics[width=\linewidth]{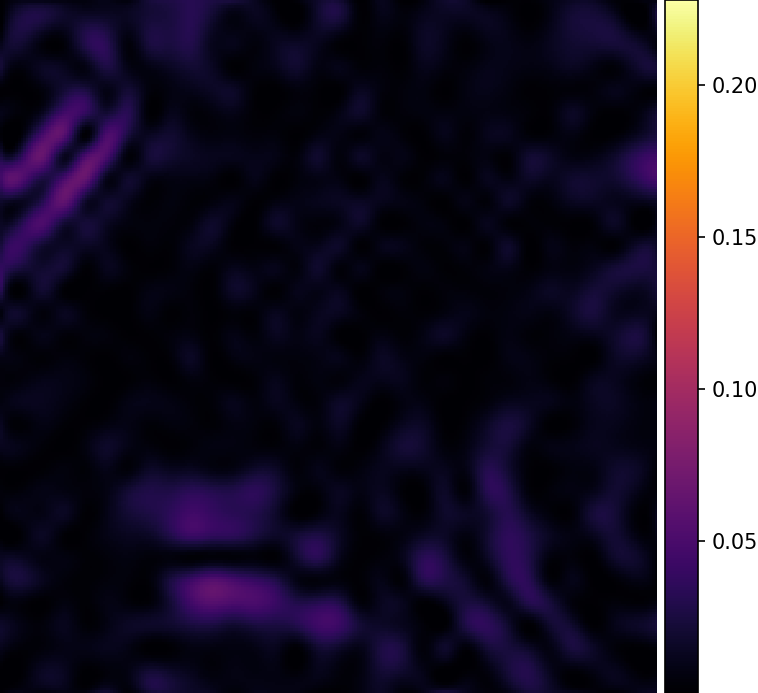} &
		\includegraphics[width=\linewidth]{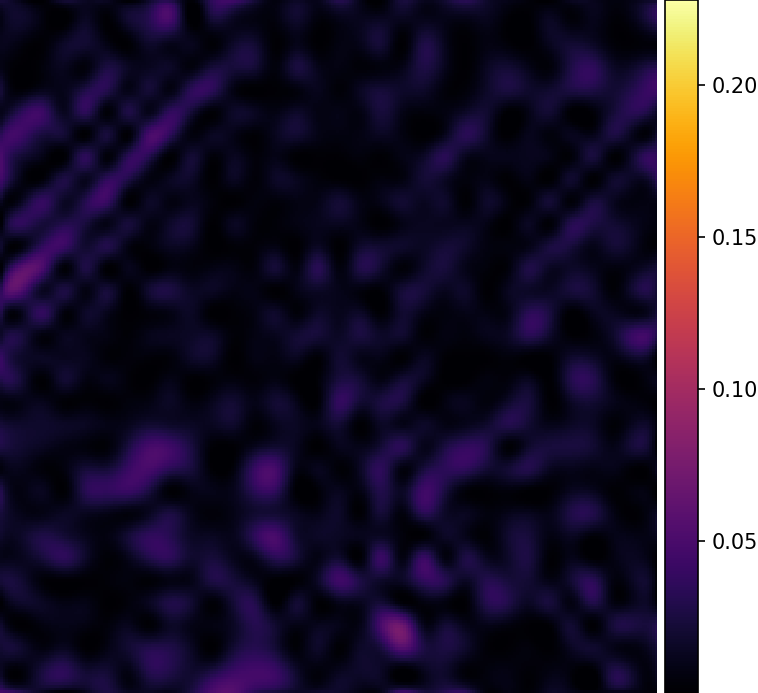} &
		\includegraphics[width=\linewidth]{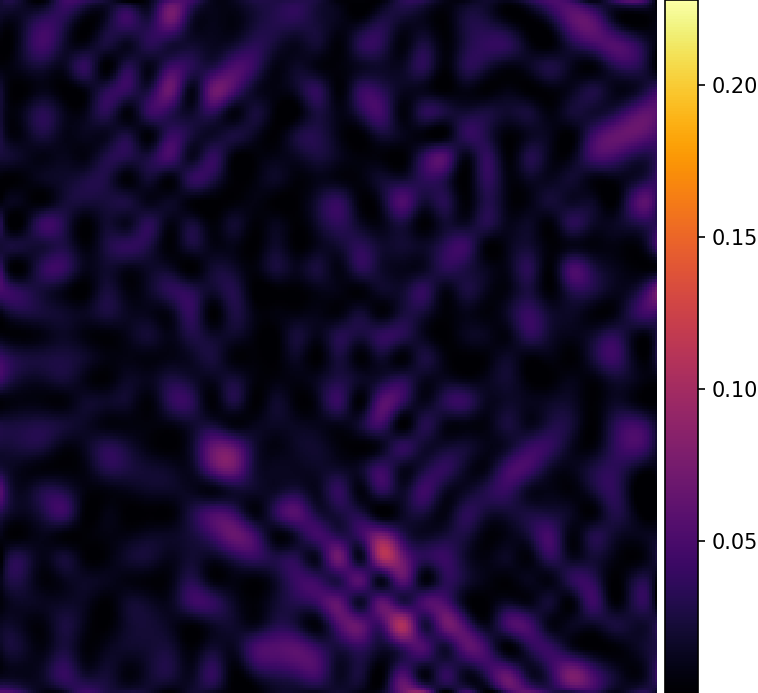} &
		\includegraphics[width=\linewidth]{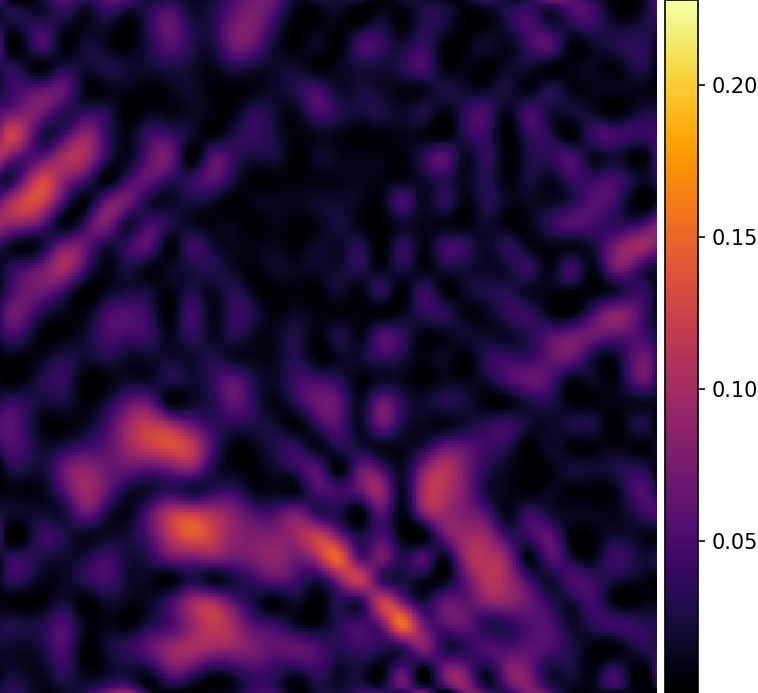} \\
		
	\end{tabular}
	
	\caption{Example output of SPOD-TrTINO for learning the solution Navier--Stokes equation from $t=0.0$ to $t=2.0$ (\Cref{subsubsect_navier_stokes}).}
	\label{fig:ns_results_matrix}
\end{figure}

\begin{figure}[h]
	\centering
	\setlength{\tabcolsep}{1.5pt} 

	\begin{tabular}{>{\centering\arraybackslash}m{0.1\textwidth} 
			*{6}{>{\centering\arraybackslash}m{0.14\textwidth}}}

		& $t=0.00$ & $t=0.25$ & $t=0.50$ & $t=0.75$ & $t=1.00$  \\
		\noalign{\vspace{4pt}} 
		Input & 
		\includegraphics[width=\linewidth]{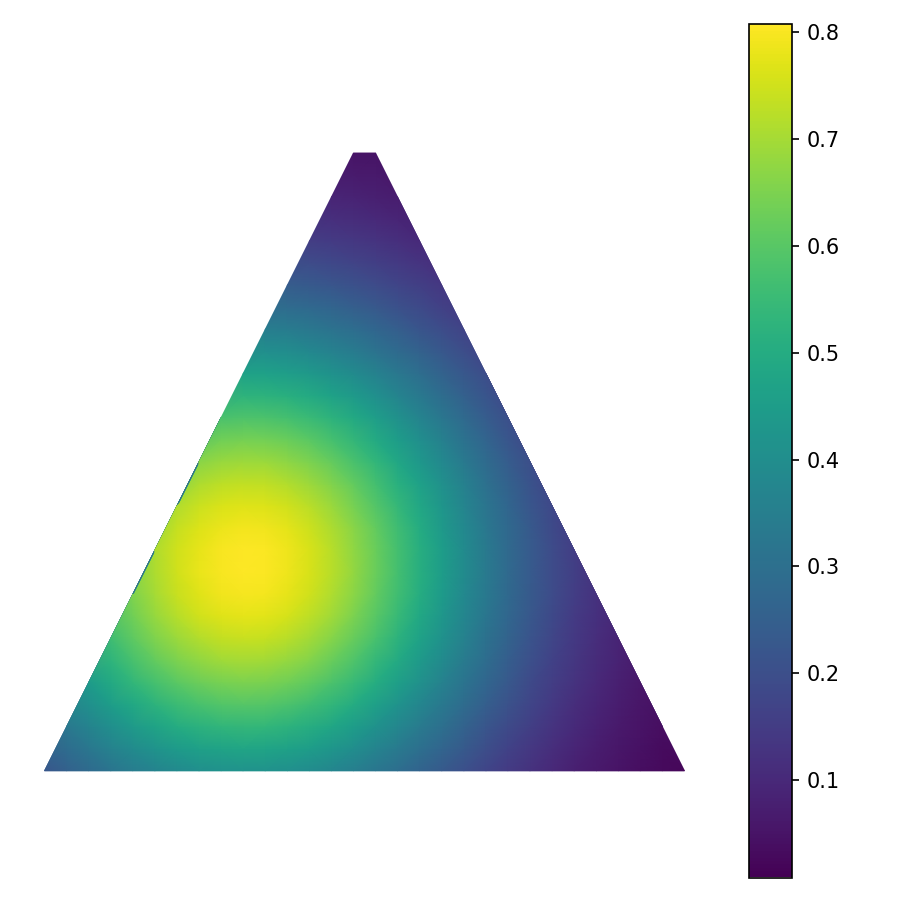} &
		\includegraphics[width=\linewidth]{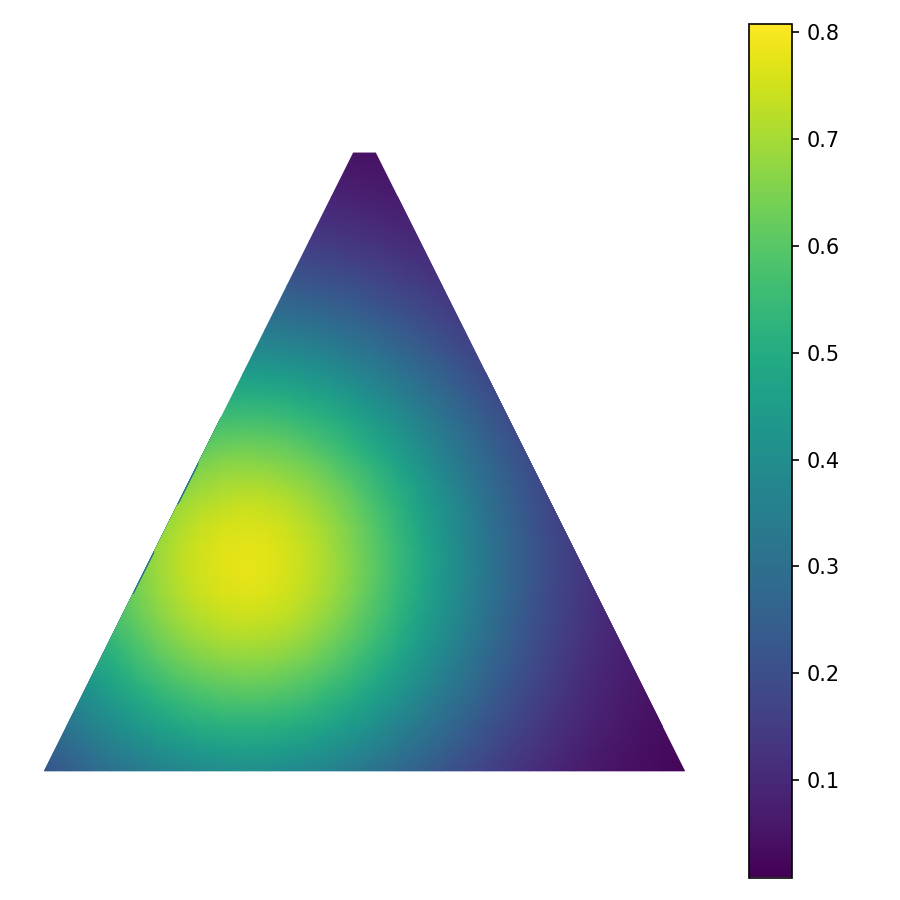} &
		\includegraphics[width=\linewidth]{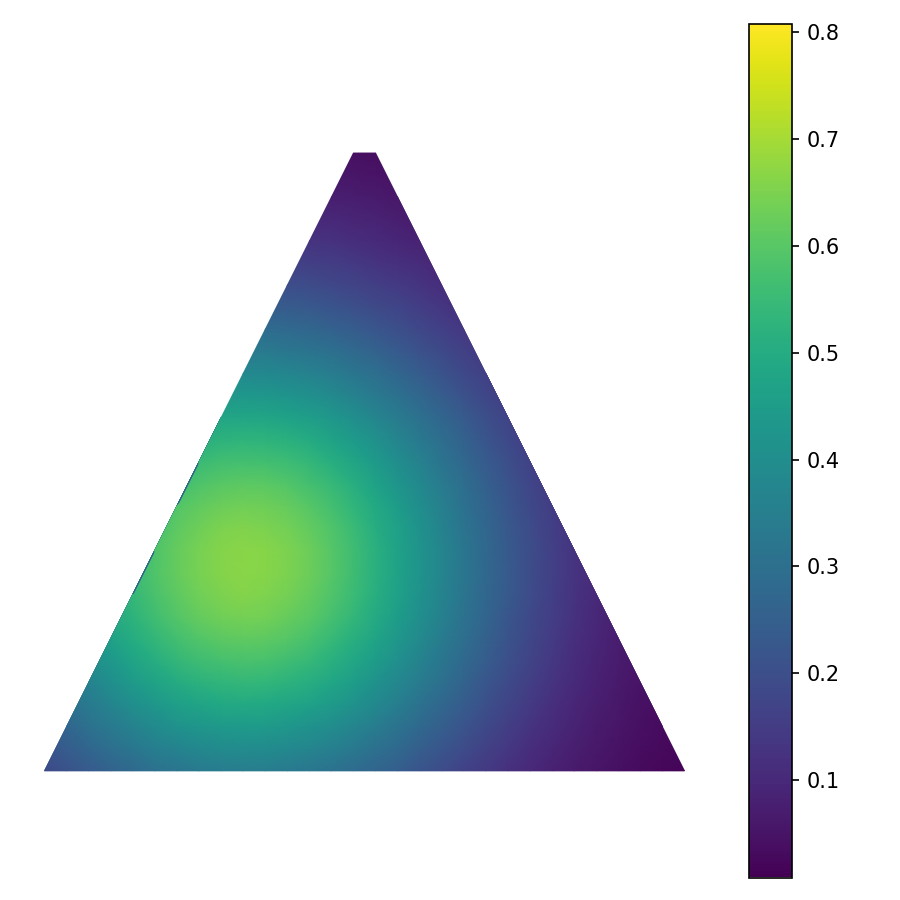} &
		\includegraphics[width=\linewidth]{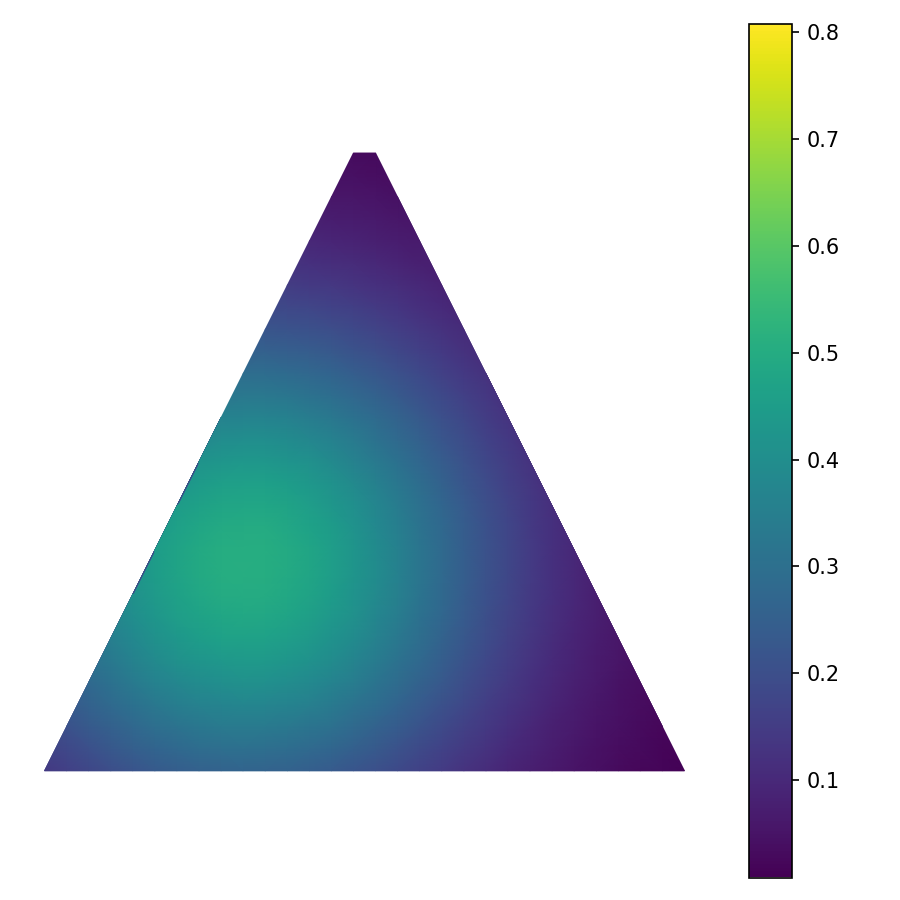} &
		\includegraphics[width=\linewidth]{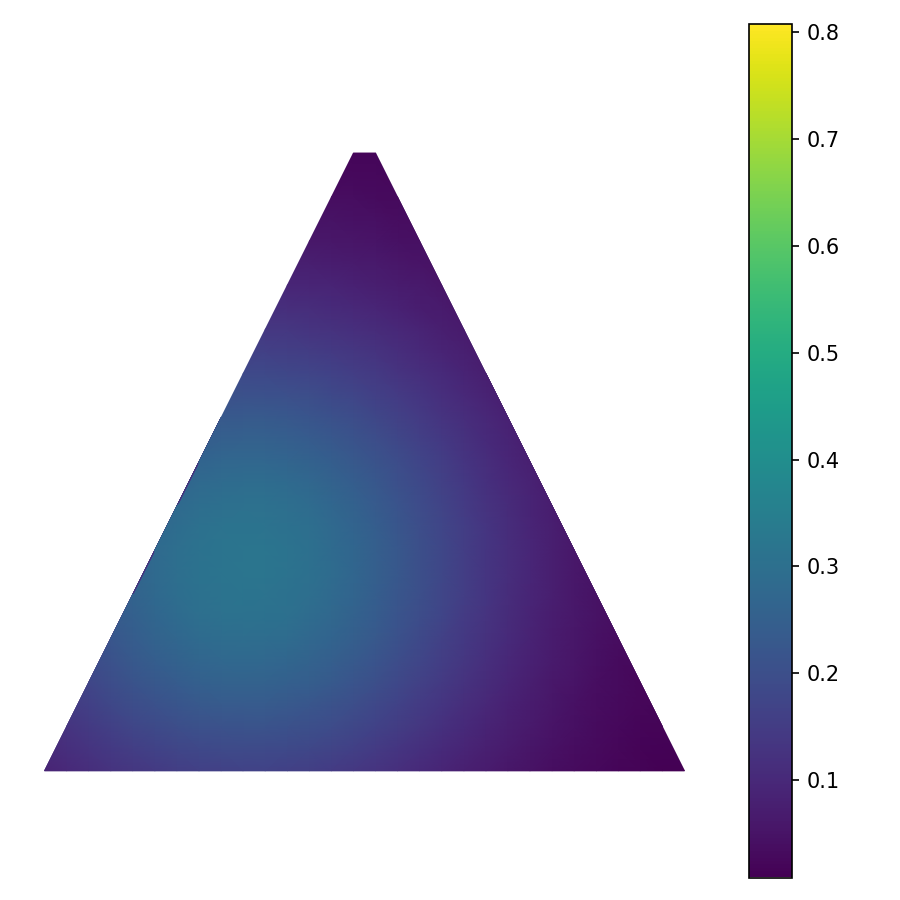} \\
		\noalign{\vspace{6pt}}
		
		Ground Truth & 
		\includegraphics[width=\linewidth]{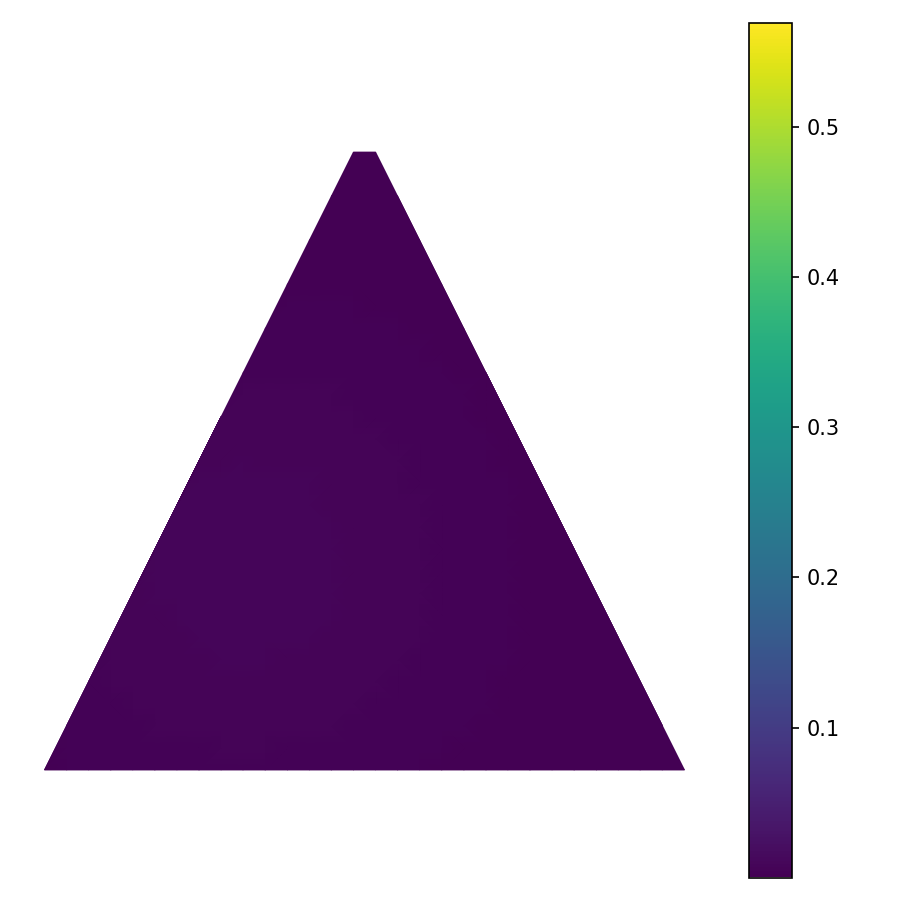} &
		\includegraphics[width=\linewidth]{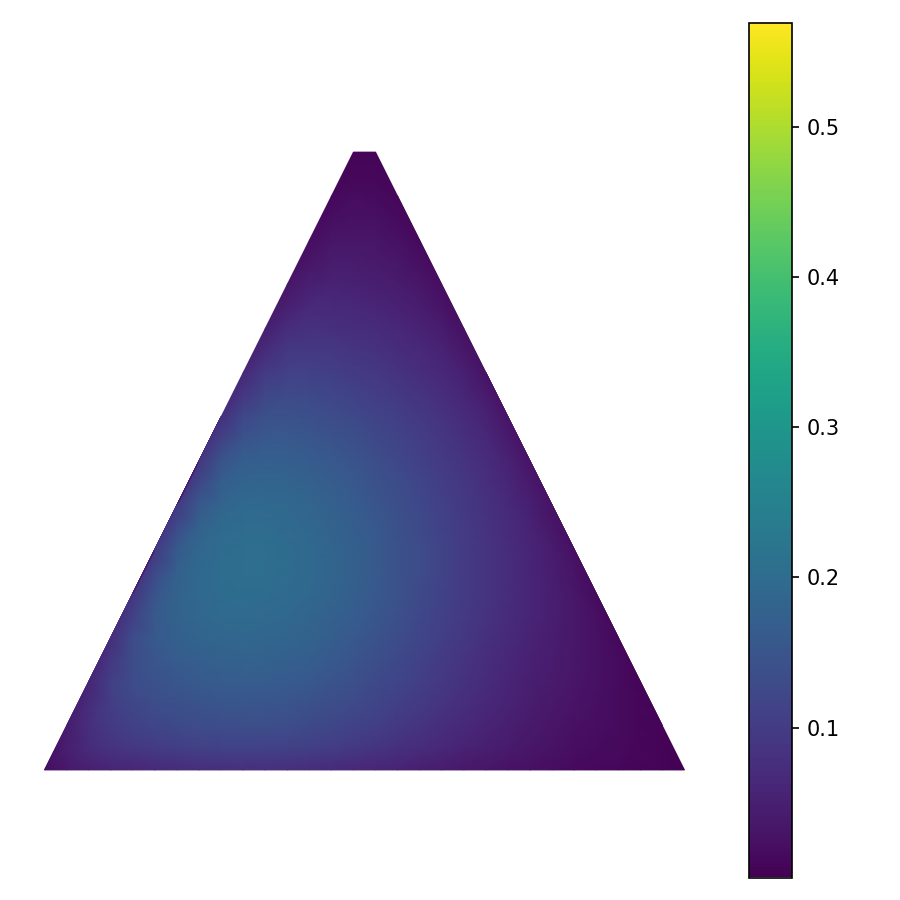} &
		\includegraphics[width=\linewidth]{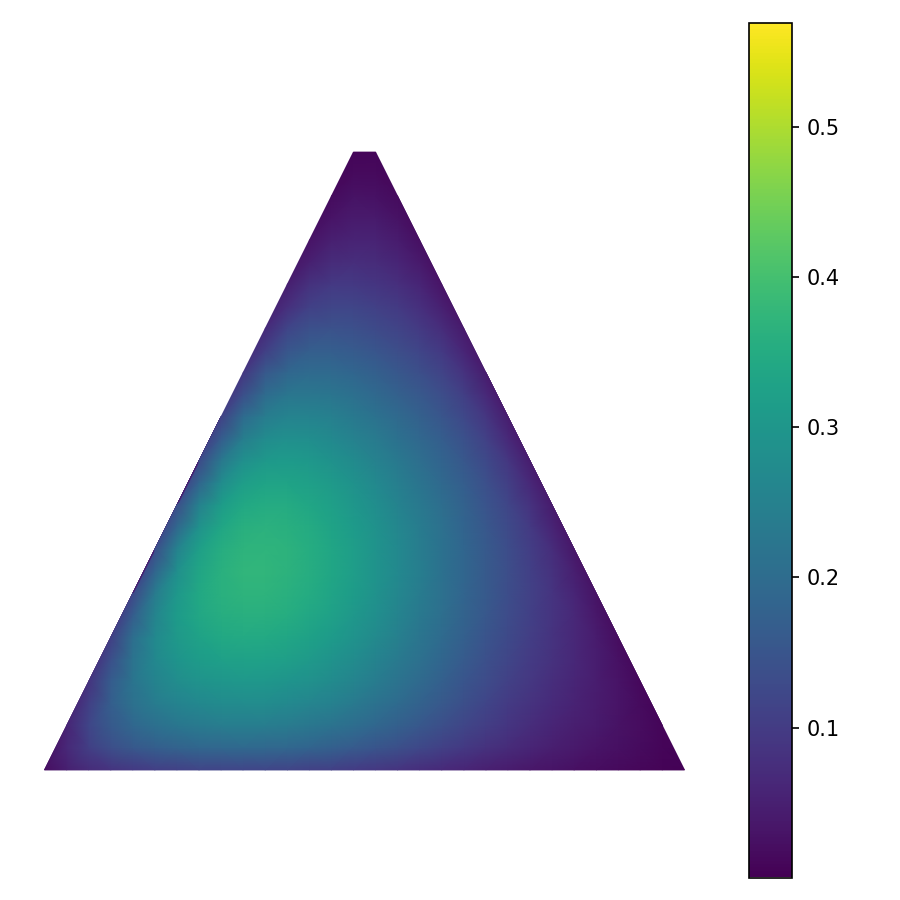} &
		\includegraphics[width=\linewidth]{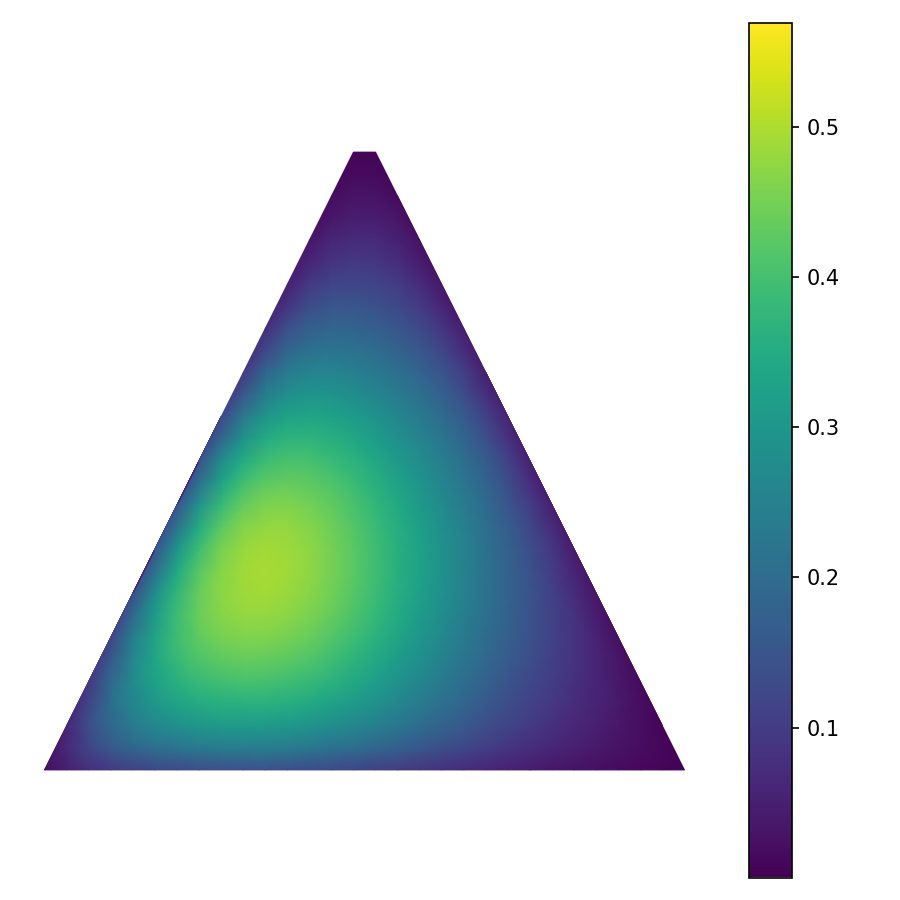} &
		\includegraphics[width=\linewidth]{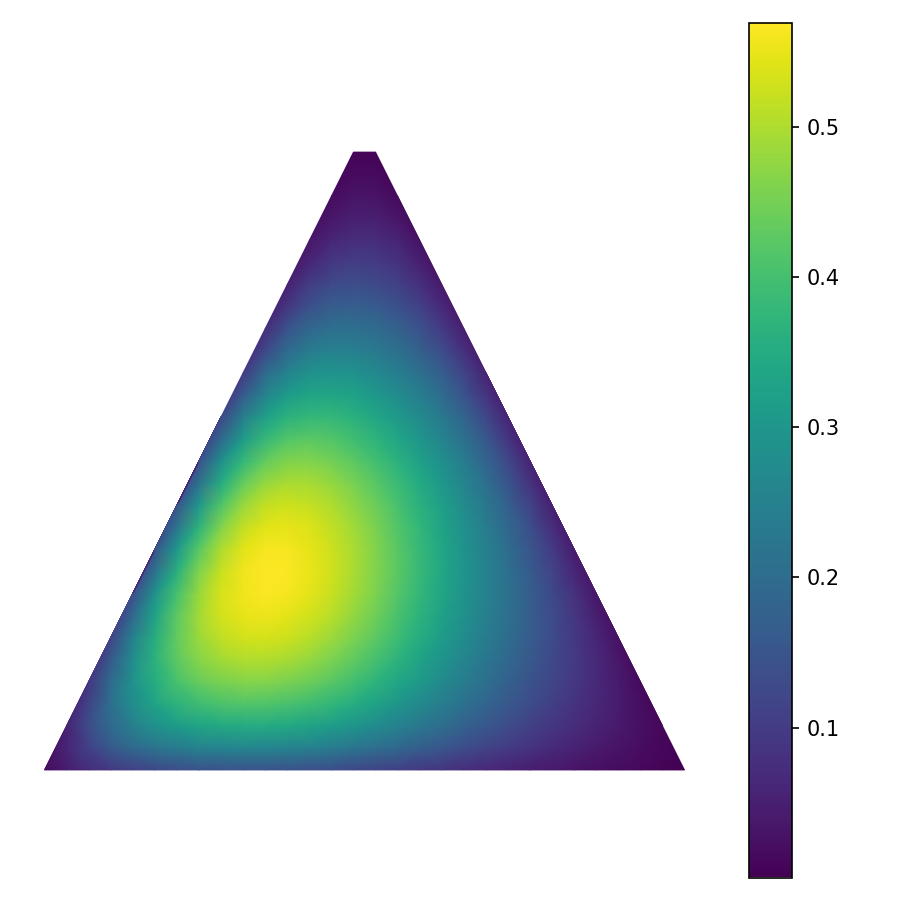} \\
		\noalign{\vspace{6pt}}

		SPOD-TrTINO Prediction & 
		\includegraphics[width=\linewidth]{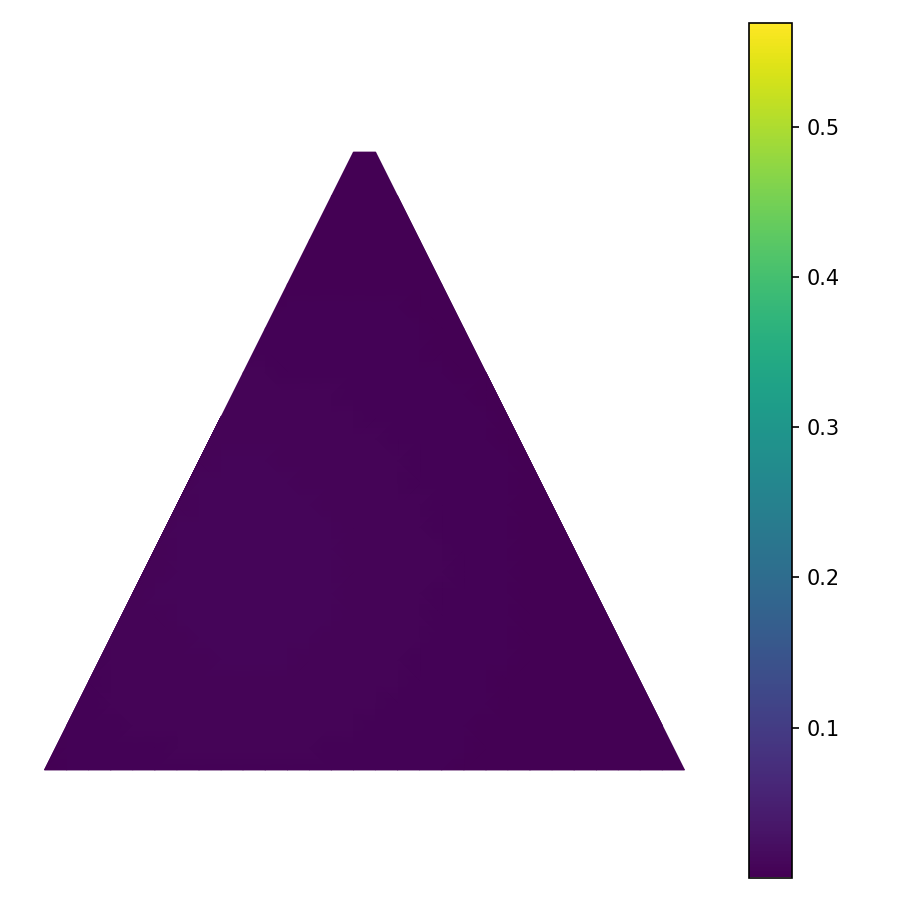} &
		\includegraphics[width=\linewidth]{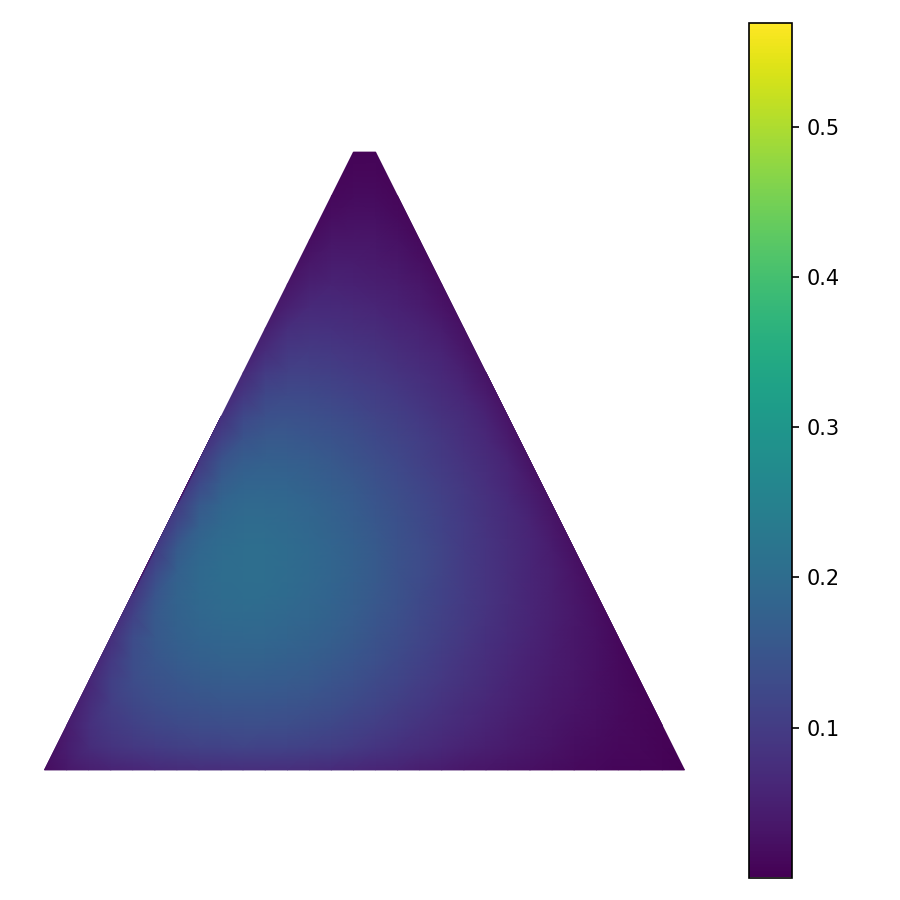} &
		\includegraphics[width=\linewidth]{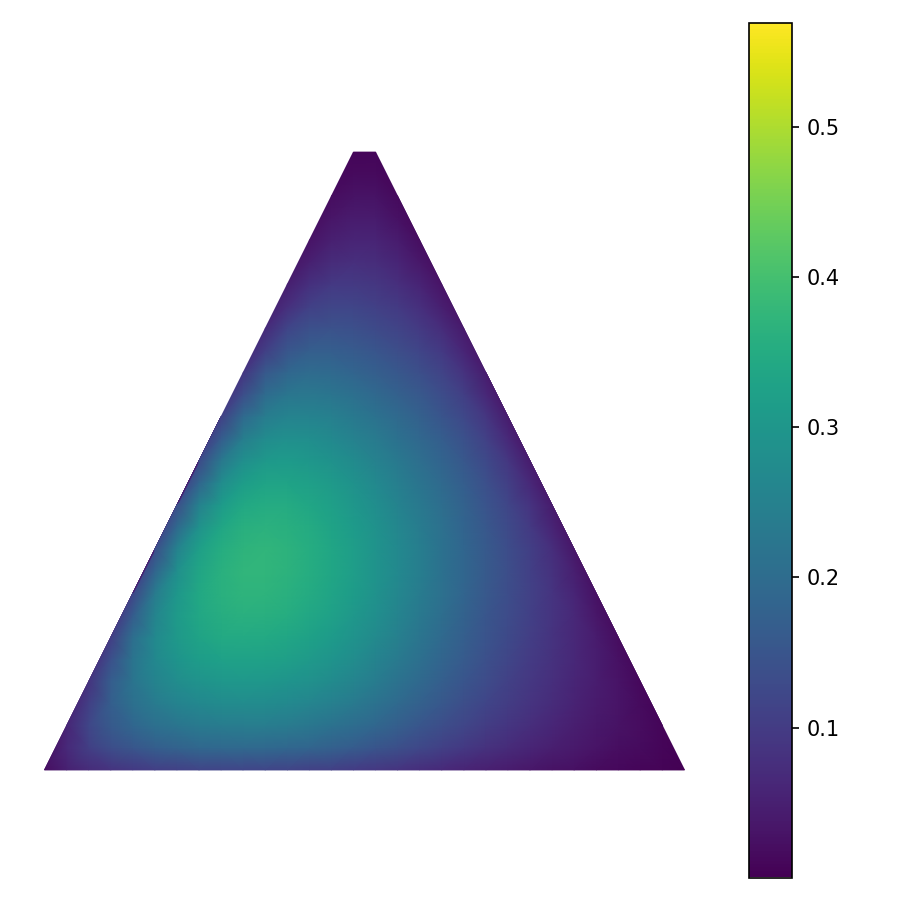} &
		\includegraphics[width=\linewidth]{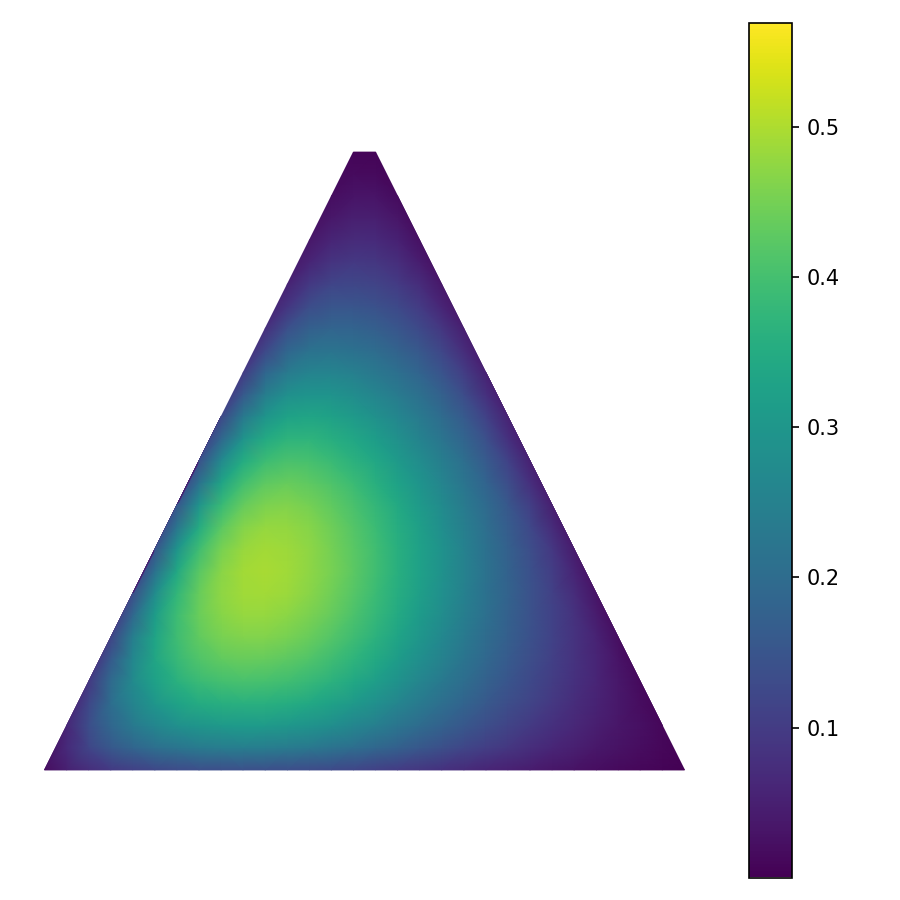} &
		\includegraphics[width=\linewidth]{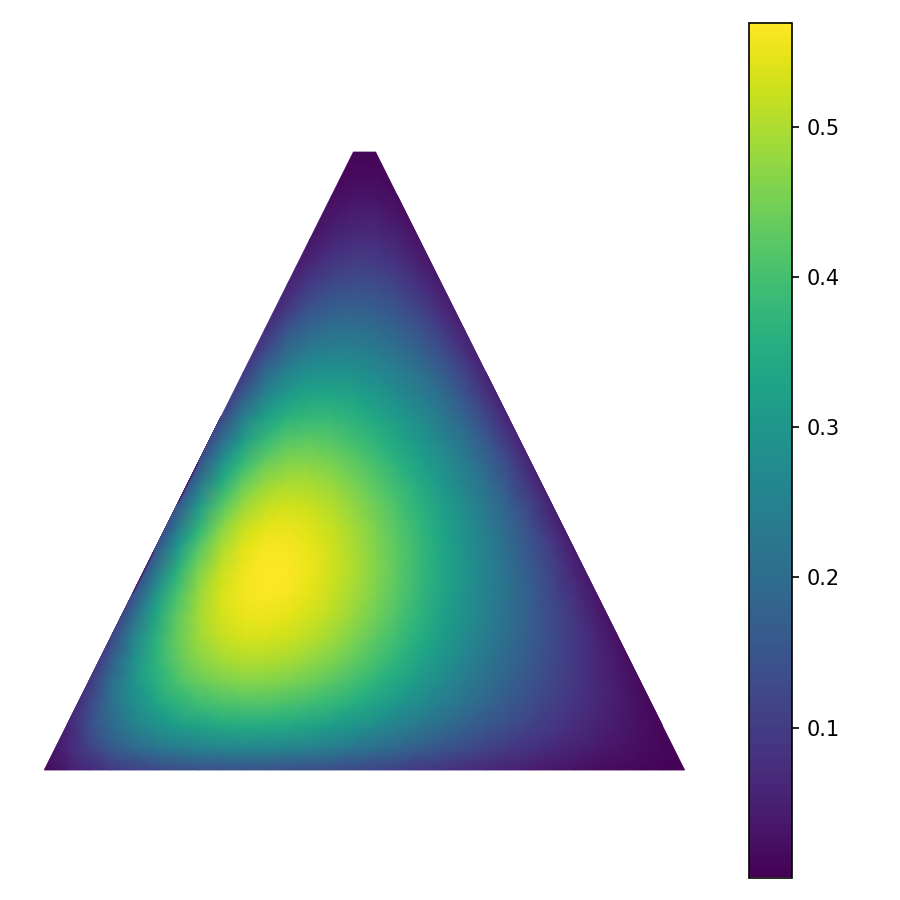} \\
		\noalign{\vspace{6pt}}

		Absolute Error & 
		\includegraphics[width=\linewidth]{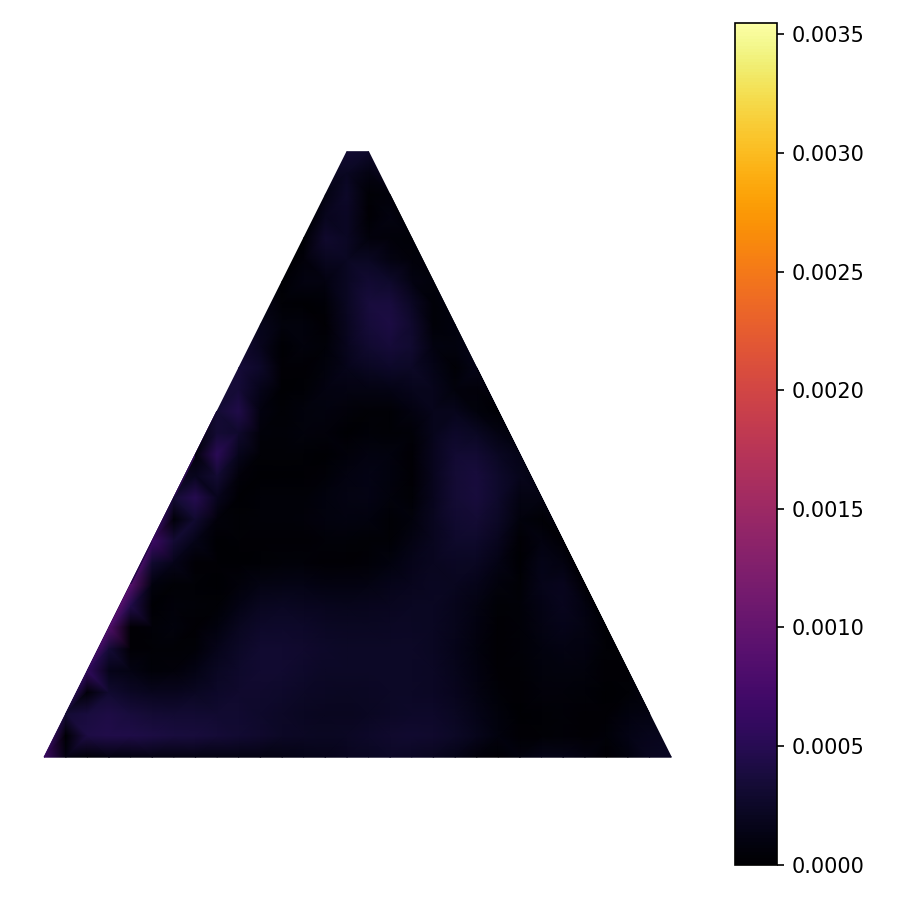} &
		\includegraphics[width=\linewidth]{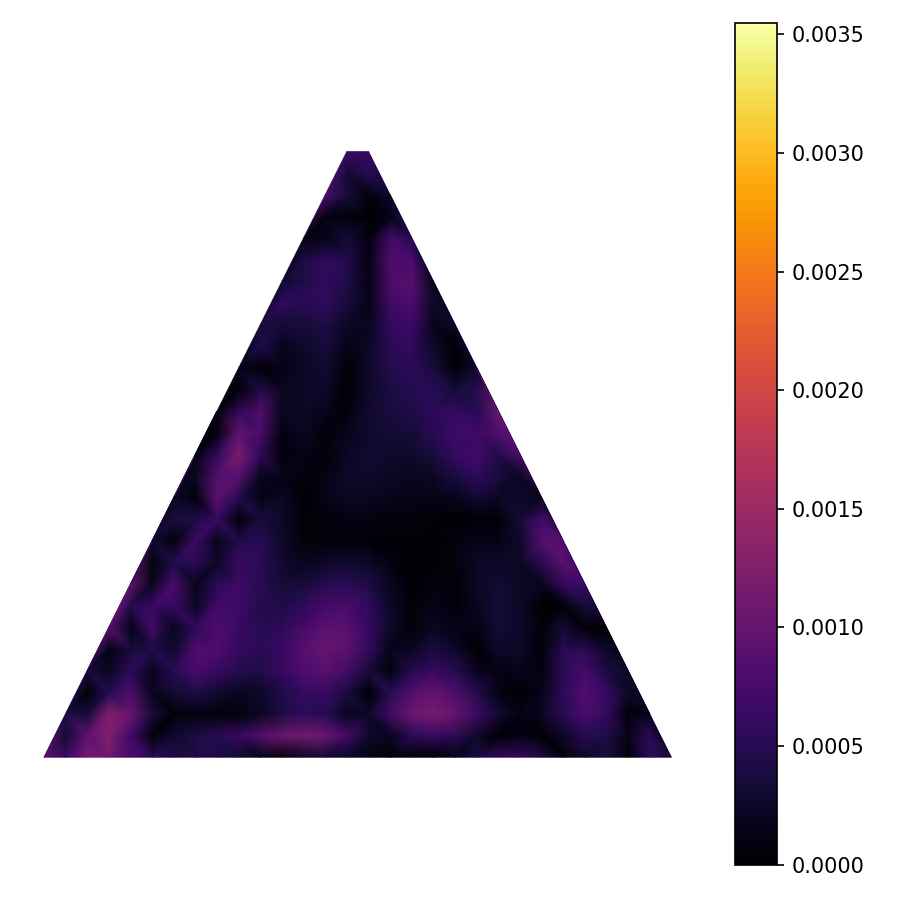} &
		\includegraphics[width=\linewidth]{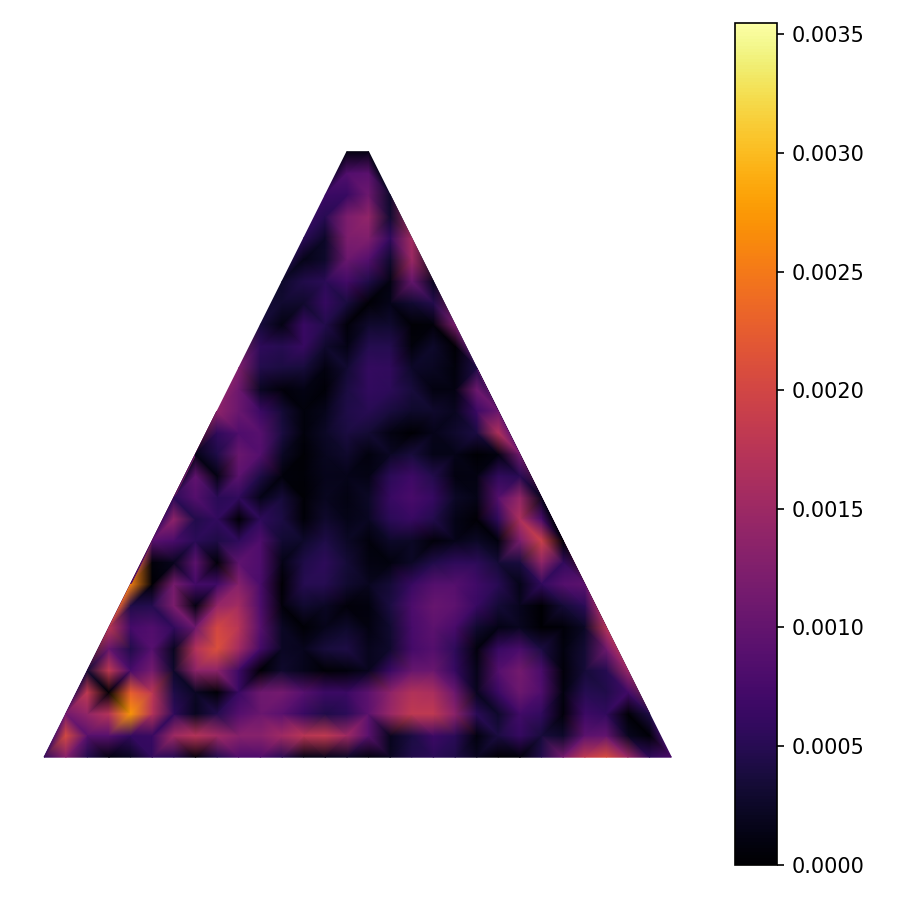} &
		\includegraphics[width=\linewidth]{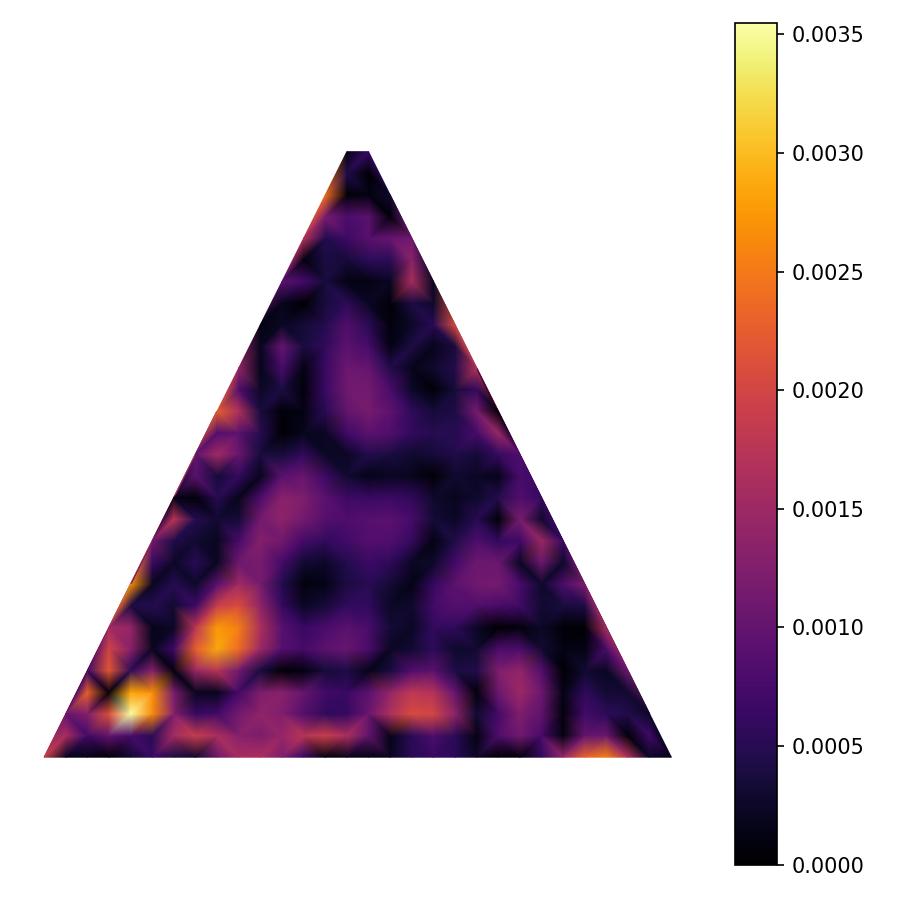} &
		\includegraphics[width=\linewidth]{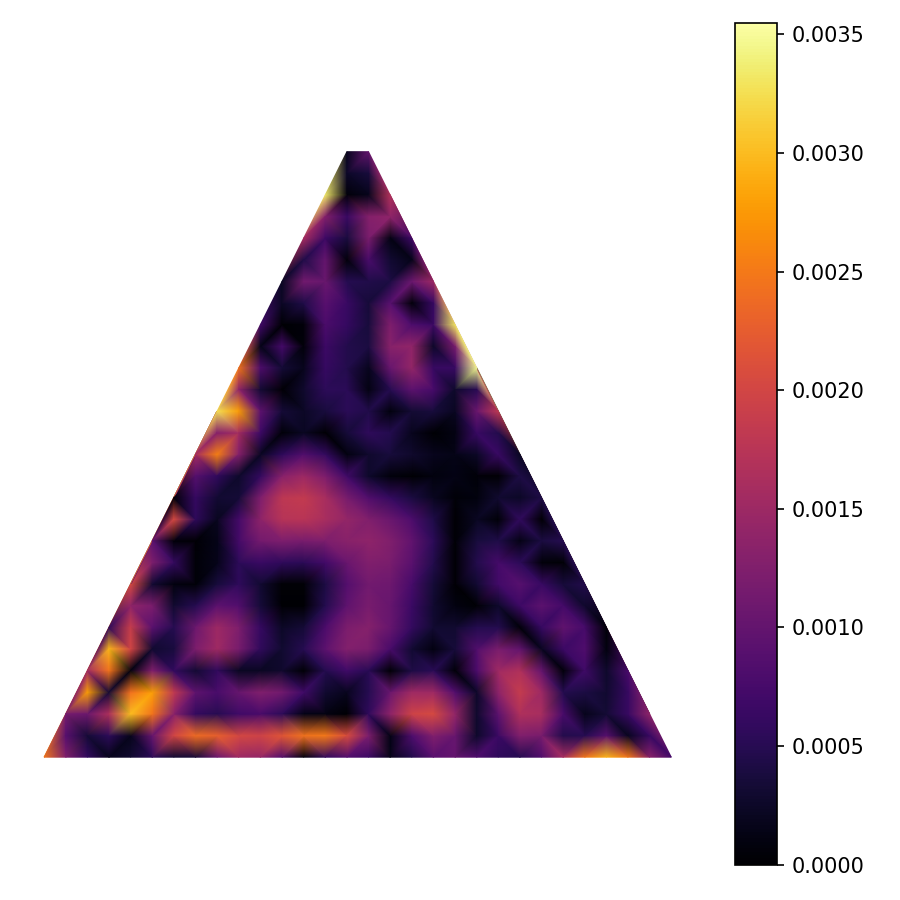} \\
		
	\end{tabular}
	
	\caption{Example output of SPOD-TrTINO for learning the solution to heat equation from $t=0.00$ to $t=1.00$ on a triangular spatial domain (\Cref{subsubsect_heat_equ}).}
	\label{fig:heat_results_matrix}
\end{figure}

\begin{figure}[h]
	\includegraphics[width=\textwidth]{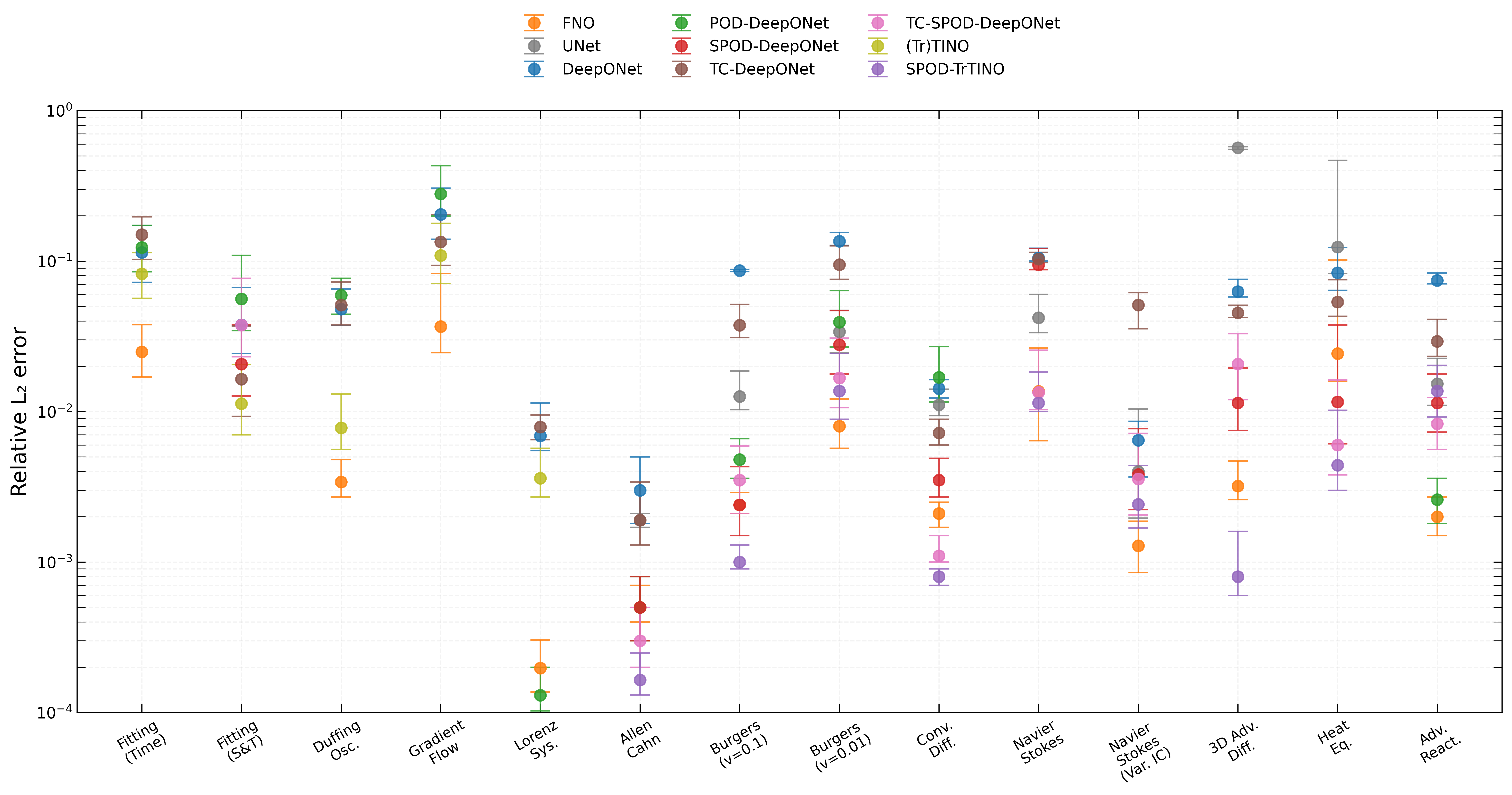}
	  \caption{Box plots of relative $L^2$ errors across all numerical experiments.
      }
	  \label{figure_numer_result}
\end{figure}

\section{Conclusion}\label{sect_conclusion}

For the deep operator learning problem for dynamic response of physical systems, we reviewe the existing deep learning models of DeepONet and TDNN, then propose TINO, the time-invariant neural operator, for learning nonlinear continuous spatiotemporal time-invariant operators.
TINO handles time information in a convolutional layer, while approximating the spatial output function space similar to the trunk net of DeepONet models.
For time-dependent PDEs with initial values that can be treated as a truncated time-invariant system, we propose TrTINO under the same framework, with a separation of time, input space, and output space variables.
The POD version of TrTINO is implemented with spatial POD for keeping a clear separation of space and time.
Numerical tests including 1-D to 3-D time-dependent PDEs justify the superior accuracy of TrTINO compared to baseline DeepONet implementations with time causality or without any time properties considered, provided the underlying PDE can be interpreted as a truncation of a time-invariant system.
We leave the embedding of selective memory in TINO to a future work.

\section*{Acknowledgment}
The author of WZ is partially supported by NSFC grants (Nos. 12571444, 12201603, 92270205) and by National Key R\&D Program of China (No. 2022YFA1005203).

% elsarticle_harv_custom: modified for keeping CAPITAL letters in title
%\bibliographystyle{elsarticle_harv_custom}
\bibliographystyle{plainnat}

\bibliography{reference.bib}

\appendix

\section{Parameter settings and data generation in numerical tests} \label{app_para_set}

\subsection{Operator networks}
For DeepONet-based models, the branch network is a 3-layer MLP (multilayer perceptron) with a width of $128$.
For variants with a learned trunk net (TINO, TrTINO, DON, TC-DON), the trunk net is a 3-layer MLP with a width of $512$.
In (Tr)TINO, the linear time convolution layer has an output dimension of $256$.
U-Net is configured with $32$ initial features. FNO hyperparameter settings (modes and widths) are summarized in Table \ref{tab:fno_settings}. Our implementation and parameter selection of FNO and U-Net referenced the settings described in \cite{lifourier}.

\begin{table}[h]
\centering
\caption{FNO Hyperparameter Settings.} \label{tab:fno_settings}
\begin{tabular}{llcc}
\toprule
Problem Dimension & Benchmark Cases & Modes & Width \\
\midrule
1D (Time) & ODEs, Space-independent operators & 16 & 64 \\
2D (Space-Time) & 1D PDEs, Space-dependent operators & $(12, 12)$ & 20 \\
3D (Space-Time) & 2D Navier-Stokes, Heat equation & $(12, 12, 12)$ & 20 \\
4D (Space-Time) & 3D Advection-Diffusion & $(12, 12, 12, 12)$ & 20 \\
\bottomrule
\end{tabular}
\end{table}

\subsection{Data generation}

The input functions and their parameter ranges have been detailed in the main text. For the sake of clarity, we summarize them again in \cref{tab:data_generation_details}. The parameters of the input function are randomly sampled in a uniform distribution within their respective domains. Both the training and testing sets are generated using the same data generation procedures. 

All ODE data were processed using the 4/5-th order Runge-Kutta scheme, with the help of MATLAB's ode45 toolkit. 1-D Allen-Cahn equation and 1-D convection-diffusion equation are numerically solved using the Method of Lines (MOL), which mainly relies on MATLAB's built-in pdepe solver. Burgers' equation uses the Fourier pseudo-spectral method for spatial discretization and the Exponential Time Differencing Runge-Kutta 4th order (ETDRK4) scheme for time integration. 2-D Navier--Stokes equation and 3-D advection-diffusion system use the Fourier pseudo-spectral method for spatial discretization and an IMEX time integration scheme that combines explicit forward Euler scheme for the nonlinear convection and forcing terms with an implicit Crank--Nicolson method for the linear diffusion, all implemented from scratch using PyTorch's FFT functions for time integration. 2-D heat equation uses a second-order central finite difference method. 1-D advection-reaction equation uses a first-order upwind finite difference method for spatial discretization and an explicit forward Euler scheme for time integration.

\begin{table}[h]
\centering
\caption{Detailed configurations of input functions and numerical methods for each benchmark.}
\label{tab:data_generation_details}
\resizebox{\columnwidth}{!}{
\begin{tabular}{
    >{\raggedright\arraybackslash}p{2.6cm} 
    >{\centering\arraybackslash}p{2.0cm} 
    >{\raggedright\arraybackslash}p{6.2cm} 
    >{\raggedright\arraybackslash}p{4.2cm} 
}
\toprule
\textbf{} & \textbf{Numerical Method} & \textbf{Input Function} & \textbf{Parameter Ranges} \\
\midrule

% Fitting tasks
\begin{tabular}{@{}l@{}} Time-only operator \\ (\Cref{subsubsect_fitting_time}) \end{tabular} & 
Analytical solution & 
$f(t) = A \sin(b t + c)$ & 
$A \in [0.5, 2], b \in [0.1, 2], c \in [0, 2\pi]$ \\

\begin{tabular}{@{}l@{}} Time-dependent op. \\ in 1-D space \\ (\Cref{subsubsect_fitting_time_space}) \end{tabular} & 
Analytical solution & 
$f(t, x) = A\sin(bt+c) e^{-(x-\mu)^2/(2\sigma^2)}$ & 
\begin{tabular}{@{}l@{}} $A \in [0.1, 1], b \in [0.2, 1], c \in [0, 2\pi],$ \\ $\mu \in [L/4, 3L/4], \sigma \in [0.1L, 0.4L]$ \end{tabular} \\

% ODEs
\begin{tabular}{@{}l@{}} Duffing oscillator \\ (\Cref{subsubsect_duffing_osc}) \end{tabular} & 
\begin{tabular}{@{}c@{}} 4/5-th order \\ Runge-Kutta \end{tabular} & 
$f(t) = A \sin(b t + c)$ & 
$A \in [0.5, 2], b \in [0.1, 2], c \in [0, 2\pi]$ \\

\begin{tabular}{@{}l@{}} Gradient flow \\ (\Cref{subsubsect_gradient_flow}) \end{tabular} & 
\begin{tabular}{@{}c@{}} 4/5-th order \\ Runge-Kutta \end{tabular} & 
$f(t) = A \sin(b t + c)$ & 
$A \in [0.5, 2], b \in [0.1, 2], c \in [0, 2\pi]$ \\

\begin{tabular}{@{}l@{}} Lorenz system \\ (\Cref{subsubsect_lorenz_sys}) \end{tabular} & 
\begin{tabular}{@{}c@{}} 4/5-th order \\ Runge-Kutta \end{tabular} & 
$f(t) = A \sin(b t + c)$ & 
$A \in [0.5, 2], b \in [0.1, 2], c \in [0, 2\pi]$ \\

% 1D PDEs
\begin{tabular}{@{}l@{}} 1-D Allen-Cahn \\ (\Cref{subsubsect_allen_cahn}) \end{tabular} & 
Method of Lines & 
$f(t,x) = A\sin(bt+c) e^{-(x-\mu)^2/(2\sigma^2)}$ & 
\begin{tabular}{@{}l@{}} $A \in [0.1, 1], b \in [0.2, 1], c \in [0, 2\pi],$ \\ $\mu \in [L/4, 3L/4], \sigma \in [0.1L, 0.4L]$ \end{tabular} \\

\begin{tabular}{@{}l@{}} 1-D Burgers' \\ (\Cref{subsubsect_burgers}) \end{tabular} & 
Pseudo-spectral ETDRK4 & 
\begin{tabular}{@{}l@{}} $f(t,x) = A\sin(bt+c) \sum_{n=1}^{10} \frac{1}{n^2} \big[$ \\ $\quad a_n \sin(2\pi n x) + b_n \cos(2\pi n x) \big]$ \end{tabular} & 
\begin{tabular}{@{}l@{}} $A \in [0.1, 1], b \in [0.2, 1], c \in [0, 2\pi],$ \\ $a_n, b_n \in [-1, 1]$ \end{tabular} \\

\begin{tabular}{@{}l@{}} 1-D conv. diff. \\ (\Cref{subsubsect_conv_diff}) \end{tabular} & 
Method of Lines & 
$f(t,x) = (A\sin(bt+c))^2 e^{-(x-\mu)^2/(2\sigma^2)}$ & 
\begin{tabular}{@{}l@{}} $A \in [0.1, 1], b \in [0.2, 1], c \in [0, 2\pi],$ \\ $\mu \in [L/4, 3L/4], \sigma \in [0.1L, 0.4L]$ \end{tabular} \\

% 2D/3D PDEs
\begin{tabular}{@{}l@{}} 2-D NS \\ (\Cref{subsubsect_navier_stokes}) \end{tabular} & 
Pseudo-spectral IMEX & 
\begin{tabular}{@{}l@{}} $f(t,\vec{x}) = A\sin(bt+c) \sum_{n=1}^{10} \frac{1}{n^2} \big[$ \\ $\quad a_{n,1} \sin(2\pi n x_1) + b_{n,1} \cos(2\pi n x_1) +$ \\ $\quad a_{n,2} \sin(2\pi n x_2) + b_{n,2} \cos(2\pi n x_2) \big]$ \end{tabular} & 
\begin{tabular}{@{}l@{}} $A \in [0.1, 1], b \in [0.2, 1], c \in [0, 2\pi],$ \\ $a_{n,i}, b_{n,i} \in [-1, 1]$ \end{tabular} \\

\begin{tabular}{@{}l@{}} 2-D NS (var. IC) \\ (\Cref{subsubsect_navier_stokes_icvar}) \end{tabular} & 
Adapted from 2-D NS & 
Adapted from 2-D NS & 
Adapted from 2-D NS \\

\begin{tabular}{@{}l@{}} 3-D adv. diff. \\ (\Cref{subsubsect_3d_adc_diff}) \end{tabular} & 
3-D Pseudo-spectral IMEX & 
\begin{tabular}{@{}l@{}} $f(t,\vec{x}) = A\sin(bt+c) \sum_{n=1}^{10} \sum_{i=1}^{3} \frac{1}{n^2} \big[$ \\ $\quad a_{n,i} \sin(2\pi n x_i) + b_{n,i} \cos(2\pi n x_i) \big]$ \end{tabular} & 
\begin{tabular}{@{}l@{}} $A \in [0.1, 1], b \in [0.2, 1], c \in [0, 2\pi],$ \\ $a_{n,i}, b_{n,i} \in [-1, 1]$ \end{tabular} \\

\begin{tabular}{@{}l@{}} 2-D heat \\ (\Cref{subsubsect_heat_equ}) \end{tabular} & 
Finite Difference Method & 
\begin{tabular}{@{}l@{}} $f(t,\vec{x}) = (A \sin(bt+c))^2 \cdot$ \\ $\quad \exp\left( -\frac{(x_1-x_{1\mathrm{c}})^2 + (x_2-x_{2\mathrm{c}})^2}{2\sigma^2} \right)$ \end{tabular} & 
\begin{tabular}{@{}l@{}} $A \in [0.1, 1], b \in [0.2, 1], c \in [0, 2\pi],$ \\ $x_{1\mathrm{c}}, x_{2\mathrm{c}} \in [0.3, 0.7], \sigma \in [0.2, 0.3]$ \end{tabular} \\

\begin{tabular}{@{}l@{}} 1-D adv. react. \\ (\Cref{subsubsect_adv_react}) \end{tabular} & 
Finite Difference Method & 
\begin{tabular}{@{}l@{}} $f(t,x) = A\sin(bt+c) \sum_{n=1}^{10} \frac{1}{n^2} \big[$ \\ $\quad a_n \sin(2\pi n x) + b_n \cos(2\pi n x) \big]$ \end{tabular} & 
\begin{tabular}{@{}l@{}} $A \in [0.1, 1], b \in [0.2, 1], c \in [0, 2\pi],$ \\ $a_n, b_n \in [-1, 1]$ \end{tabular} \\

\bottomrule
\end{tabular}
}
\end{table}

\end{document}